\documentclass[11pt,leqno]{article}
\usepackage{a4wide,amsmath}
\usepackage{amsfonts,amsbsy,amscd}
\usepackage{amssymb}
\usepackage{color}
\usepackage{tikz}
\usepackage[colorlinks]{hyperref}
\newtheorem{theorem}{Theorem}[section]

\newtheorem{lemma}[theorem]{Lemma}
\newtheorem{proposition}[theorem]{Proposition}

\newtheorem{definition}[theorem]{Definition}
\newtheorem{remark}[theorem]{Remark}

\numberwithin{equation}{section}

\newcommand{\RR}{\mathbb{R}}
\newcommand{\CC}{\mathbb{C}}
\newcommand{\NN}{\mathbb{N}}

\newcommand{\HH}{\mathbb{H}}

\newcommand{\esssup}{\mathop{\rm {ess\,sup}}\limits}   

\newcommand{\RE}{\mathop{\Re\mathfrak{e}}}  
\newcommand{\IM}{\mathop{\Im\mathfrak{m}}}  

\newcommand{\ee}{\mathrm{e}}
\newcommand{\ii}{\mathrm{i}}
\renewcommand{\colon}{:\,}

\newcommand{\eqdef}{\stackrel{{\rm {def}}}{=}}   

\newcommand{\Sgn}{\mathop{\rm {sign}}}

\newcommand{\Square}{$\sqcap$\hskip -1.5ex $\sqcup$}
\newcommand{\Blacksquare}{\vrule height 1.7ex width 1.7ex depth 0.2ex }
\newcommand{\proof}{{\em Proof. }}
\newcommand{\qed}{$\;$\Blacksquare}

\definecolor{violet}{rgb}{0.5,0,0.5}
\definecolor{orange}{cmyk}{0,0.3,0.7,0}

\newcommand{\td}{\,\mathrm{d}}
 
\title{\null\vspace*{-2.0cm}
      On the Heston Model with Stochastic Volatility:\\
      Analytic Solutions and Complete Markets%
\footnote{The work both authors was supported in part by
       le Minist\`ere des Affaires \'Etrang\`eres (France)
       and
       the German Academic Exchange Service (DAAD, Germany)
       within the exchange program ``PROCOPE''.\vspace{2.0mm}}
\vspace*{0.5cm}
}
 
\author{%
        B\'en\'edicte Alziary%
\thanks{{\it e-mail:} {\tt benedicte.alziary@ut-capitole.fr}}
\vspace*{0.3cm}
\\
        Toulouse School of Economics, I.M.T.,
        Universit\'e de Toulouse -- Capitole \\
        21 All\'ees de Brienne,
        F--31000 Toulouse Cedex, France \\
\vspace*{0.3cm}
\\
        and \\
\vspace*{0.3cm}
\and
        Peter Tak\'a\v{c}%
\thanks{{\it e-mail:} {\tt peter.takac@uni-rostock.de}.$\;$
        A part of this research was performed while
        this author was a visiting professor at
        Toulouse School of Economics, I.M.T.,
        Universit\'e de Toulouse -- Capitole, Toulouse, France.}
\vspace*{0.3cm}
\\
        Universit\"at Rostock,
        Institut f\"ur Mathematik \\
        Ulmenstra{\ss}e~69, Haus~3,
        D--18057 Rostock, Germany
\vspace*{0.5cm}
}
 
\date\today

\begin{document}
\baselineskip=16pt plus 1pt minus 1pt

\maketitle
\baselineskip=14pt plus 1pt minus 1pt
\baselineskip=16pt plus 1pt minus 1pt
 
\newpage

\begin{abstract}
We study the Heston model for pricing European options on stocks
with stochastic volatility.
This is a Black\--Scholes\--type equation whose spatial domain
for the logarithmic stock price $x\in \RR$
and the variance $v\in (0,\infty)$
is the half\--plane $\HH = \RR\times (0,\infty)$.
The {\it volatility\/} is then given by $\sqrt{v}$.
The diffusion equation for the price of the European call option
$p = p(x,v,t)$ at time $t\leq T$
is parabolic and degenerates at the boundary
$\partial \HH = \RR\times \{ 0\}$ as $v\to 0+$.
The goal is to hedge with this option against
volatility fluctuations, i.e., the function
$v\mapsto p(x,v,t)\colon (0,\infty)\to \RR$
and its (local) inverse are of particular interest.
We prove that
$\frac{\partial p}{\partial v}(x,v,t) \not= 0$
holds almost everywhere in $\HH\times (-\infty,T)$
by establishing the analyticity of $p$ in both,
space $(x,v)$ and time $t$ variables.
To this end, we are able to show that
the Black\--Scholes\--type operator,
which appears in the diffusion equation,
generates a holomorphic $C^0$-semigroup in a suitable weighted
$L^2$-space over $\HH$.
We show that the $C^0$-semigroup solution can be extended to
a holomorphic function in a complex domain in $\CC^2\times \CC$,
by establishing some new a~priori weighted $L^2$-estimates over
certain complex ``shifts'' of $\HH$
for the unique holomorphic extension.
These estimates depend only on
the weighted $L^2$-norm of the terminal data over $\HH$ (at $t=T$).
\end{abstract}

\par\vfill
\vspace*{0.5cm}
\noindent
\begin{tabular}{lll}
{\bf 2010 Mathematics Subject Classification:}
& Primary   & 35B65, 91G80;\\
& Secondary & 35K65, 35K15.\\
\end{tabular}

\par\vspace*{0.5cm}
\noindent
\begin{tabular}{ll}
{\bf Key words:}
& Heston model; stochastic volatility; \\
& Black\--Scholes equation; European call option; \\
& degenerate parabolic equation; terminal value problem; \\
& holomorphic extension; analytic solution \\
\end{tabular}
 
\baselineskip=16pt plus 1pt minus 1pt
\parskip=2mm plus .5mm minus .5mm
\newpage
 

\section{Introduction}
\label{s:Intro}

For several decades,
simple market models have been very important and useful products of
numerous mathematical studies of financial markets.
Several of them have become very popular and are extensively used by
the financial industry
({\sc F.\ Black} and {\sc M.\ Scholes} \cite{Black-Scholes},
 {\sc S.~L.\ Heston} \cite{Heston}, and
 {\sc J.-P.\ Fouque}, {\sc G.\ Papanicolaou}, and {\sc K.~R.\ Sircar}
 \cite{FouqPapaSir} to mention only a few).
These models are usually concerned with asset pricing
in a volatile market under clearly specified rules that are supposed
to guarantee ``fair pricing''
(e.g., arbitrage\--free prices in {\sc T.\ Bj\"ork} \cite{Bjoerk-3rd}).

Assets are typically represented by securities (e.g., bonds, stocks)
and their derivatives (such as options on stocks and similar contracts).
An important role of a derivative is to assess
the volatile behavior of a particular asset and replace it by
a suitable portfolio containing both,
the asset itself and its derivatives, in such a way that
the entire portfolio is less volatile than the asset itself.
A common way to achieve this objective is to add a derivative
on the volatile asset to the portfolio containing this asset.
This procedure, called hedging, is closely connected with
the problem of {\em\bfseries market completion\/}
({\sc M.~H.~A.\ Davis} \cite{Davis-royal}),
 {\sc M.\ Romano} and {\sc N.\ Touzi} \cite{RomanoTouzi}).
There have been a number of successful attempts to obtain
a market completion by (call or put) options on stocks.
The pricing of such options involves various kinds of
the Black\--Scholes\--type equations.
These attempts are typically based on
probabilistic, analytic, and numerical techniques,
some of them including even explicit formulas, cf.\
{\sc Y.\ Achdou} and {\sc O.\ Pironneau} \cite[Chapt.~2]{Achdou-Pironn}.
The basic principle behind all Black\--Scholes\--type models is that
the model must be {\em\bfseries arbitrage\--free\/}, that is,
any {\em arbitrage opportunity\/} must be excluded
which is possible only if
there exists an equivalent probability measure such that
the option price is a stochastic process that is a martingale
under this measure
(in which case it is called a {\it martingale measure\/}, cf.\
 {\sc T.\ Bj\"ork} \cite[{\S}3.3, pp.\ 32--33]{Bjoerk-3rd}).
It\^{o}'s formula then yields an equivalent linear parabolic equation
which will be the object of our investigation, cf.\
{\sc M.~H.~A.\ Davis} \cite{Davis-royal}.
Throughout our present work we study
the {\it\bfseries Heston model\/} of pricing for
{\it European call options on stocks\/}
with {\it\bfseries stochastic volatility\/}
({\sc S.~L.\ Heston} \cite{Heston})
by abstract analytic methods coming from partial differential equations
(PDEs, for short) and functional analysis.
Without any option, derivative, or other contingent claim added to
the Heston model, this model represents an incomplete market.
In probabilistic terms, this means that the martingale measure
mentioned above is possibly not unique.
We use a PDE to give a rigorous analytic formulation of
Heston's model in the next section (Section~\ref{s:reform}).
Our main results are presented in a functional\--analytic setting
in Section~\ref{s:Main}.

In our simple market, described by
the {\it\bfseries Heston stochastic volatility model\/}
(Heston model, for short),
market completion by a European call option on the stock has
the following meaning:
The basic quantities are
the {\it maturity time\/} $T$
(called also the {\it exercise time\/}), $0 < T$ $< \infty$,
at which the stock option matures;
the {\it real time\/} $t$, $-\infty < t\leq T$;
the {\it time to maturity\/} $\tau = T-t\geq 0$,
$0\leq \tau < \infty$;
the {\it spot price\/} of stock $S_t$ ($S_t > 0$) and
the (stochastic) {\it variance\/} of the stock market
$V_t$ ($V_t > 0$) at time $t\leq T$;
$\sqrt{V_t}$ is associated with
the (stochastic) {\it volatility\/} of the stock market;
the {\it strike price\/} (exercise price)
$K\equiv \mathrm{const} > 0$
of the stock option at maturity, a European call or put option;
a given (non\-negative) {\it payoff function\/}
$\hat{h}(S_T,V_T) = (S_T - K)^{+}$
at time $t = T$ (i.e., $\tau = 0$) for a European call option;
and the (call or put) {\it option price\/}
$P_t = U(S_t,V_t,t) > 0$ at time $t$,
given the stock price $S_t$ and the variance $V_t$.
In the derivation of {\sc S.~L.\ Heston}'s model \cite{Heston},
which is a system of two stochastic differential equations
for the pair $(S_t,V_t)$,
It\^{o}'s formula yields a diffusion equation for
the unknown option price $P_t =$ $U(S_t,V_t,t) > 0$ at time $t$
which depends only on the stock price $S_t$ and
the variance $V_t$ at time $t$.
This allows us to replace the {\it relative logarithmic stock price\/}
$X_t = \ln (S_t/K)$, a stochastic process valued in
$\RR = (-\infty, \infty)$, and the variance $V_t$,
another stochastic process valued in $(0,\infty)$, respectively,
by a pair of (independent) space variables $(x,v)$
valued in the open half\--plane
$\mathbb{H}\eqdef \mathbb{R}\times (0,\infty) \subset \RR^2$.
Consequently, the option price
\begin{math}
  P_t = p(X_t,V_t,t)\eqdef U\left( K\ee^{X_t} ,\, V_t, t\right)
\end{math}
is a stochastic process whose values at time $t$ ($t\leq T$)
are determined by the values of $(X_t,V_t)$.
Its terminal value, $P_T$ at maturity time $t=T$, is given by
\begin{equation*}
  P_T = p(X_T,V_T,T)
      = K\left( \mathrm{e}^{X_T} - 1\right)^{+} = (S_T - K)^{+}
  \quad\mbox{ for }\, (X_T,V_T)\in \HH \,.
\end{equation*}
The well\--known arbitrage\--free option pricing
({\sc T.\ Bj\"ork} \cite[Chapt.~7, pp.\ 92--108]{Bjoerk-3rd})
then yields the expectation formula
\begin{equation}
\label{e:expect}
  p(x,v,t) = K\cdot
    \mathrm{e}^{-r(T-t)}\, \mathbb{E}_{\mathbb{P}}
    \left[ \left( \mathrm{e}^{X_T} - 1\right)^{+}
           \mid X_t = x ,\ V_t = v \right]
\end{equation}
for $(x,v)\in \HH$ and $t\in (-\infty,T]$; see, e.g.,
{\sc J.-P.\ Fouque}, {\sc G.\ Papanicolaou}, and {\sc K.~R.\ Sircar}
\cite[{\S}2.4--2.5, pp.\ 42--48]{FouqPapaSir}.
In particular, the terminal condition at $t=T$ is fulfilled,
\begin{equation}
\label{e:mature}
  p(x,v,T) = K\left( \mathrm{e}^x - 1\right)^{+}
  \quad\mbox{ for }\, (x,v)\in \HH \,.
\end{equation}

The option price $p = p(x,v,t)\equiv p_{\tau}(x,v)$,
where $\tau = T-t\geq 0$, is determined by
an equivalent, risk neutral martingale measure
(\cite{Davis-royal, RomanoTouzi}),
which yields the stochastic process $(P_t)_{0\leq t\leq T}$.
This measure is unique if and only if
every contingent claim can be replicated by
a self\--financed trading strategy using bond, stock, and option;
that is to say, if and only if the option price $(P_t)_{0\leq t\leq T}$
{\it\bfseries completes the market\/}
({\sc J.~M.\ Harrison} and {\sc S.~R. Pliska}
 \cite{HarrisPliska-1, HarrisPliska-2}).
Applying It\^{o}'s formula to this process, one concludes that,
equivalently to the probabilistic expectation formula
\eqref{e:expect} for $p(x,v,t)$,
this option price can be calculated directly from
a partial differential equation of parabolic type with
the terminal value condition \eqref{e:mature}.
Thus, given the (relative logarithmic) stock price $x\in \RR$
at a fixed time $t\in (-\infty, T]$, the function
$\tilde{p}_{x,t}\colon v\mapsto p(x,v,t)$
yields the (unique) option price for every $v\in (0,+\infty)$.
According to
{\sc I.\ Bajeux\--Besnainou} and {\sc J.-Ch.\ Rochet}
\cite[p.~12]{Bajeux-Rochet},
the characteristic property of a complete market is that
$\tilde{p}_{x,t}\colon (0,+\infty)\to \RR_+$
is injective (i.e., one\--to\--one), which means that
any particular option value $p = \tilde{p}_{x,t}(v)$
cannot be attained at two different values of the variance
$v\in (0,+\infty)$.
We take advantage of this property to give
an alternative definition of a {\it\bfseries complete market\/}
using differential calculus rather than probability theory, see
our Definition~\ref{def-complete} in Section~\ref{s:Appl}.
This is a purely mathematical problem that we solve in this article
for the Heston model by analytic methods, with a help from
\cite[Sect.~5]{Bajeux-Rochet} and the work by
{\sc M.~H.~A.\ Davis} and {\sc J.\ Ob{\l}{\'o}j} \cite{Davis-Obloj};
see Section~\ref{s:Appl} below, Theorem~\ref{thm-complete}.
We refer the reader to the monograph by
{\sc F.\ Delbaen} and {\sc W.\ Schachermayer} \cite{Delbaen-Schacher}
for an up\--to\--date treatment of
complete markets with no arbitrage opportunity
(particularly in Chapter~9, pp.\ 149--205).

There are several other stochastic volatility models, see, e.g.,
those listed in \cite[Table 2.1, p.~42]{FouqPapaSir}
and those treated in
\cite{FouqPapaSir, Hull-White, Lewis, Stein-Stein, Wiggins},
that are already known to allow or may allow market completion
by a European call or put option.
However, the rigorous proofs of market completeness (and their methods)
vary from model to model; cf.\
{\sc T.\ Bj\"ork} \cite{Bjoerk-3rd}.
Some of them are more probabilistic
({\sc R.~M.\ Anderson} and {\sc R.~C.\ Raimondo} \cite{Anders-Raimond}
with {\it ``endogenous completeness''\/}
of a diffusion driven equilibrium market,
{\sc I.\ Bajeux\--Besnainou} and {\sc J.-Ch.\ Rochet} \cite{Bajeux-Rochet},
{\sc J.\ Hugonnier}, {\sc S.\ Malamud}, and {\sc E.\ Trubowitz}
\cite{HugoMalaTrubo},
{\sc D.\ Kramkov} \cite{Kramkov-1},
{\sc D.\ Kramkov} and {\sc S.\ Predoiu} \cite{Kramkov-Pred}, and
{\sc M.\ Romano} and {\sc N.\ Touzi} \cite{RomanoTouzi}),
others more analytic (PDEs), e.g., in
{\sc M.~H.~A.\ Davis} \cite{Davis-royal},
{\sc M.~H.~A.\ Davis} and {\sc J.\ Ob{\l}{\'o}j} \cite{Davis-Obloj},
and
{\sc P.\ Tak\'a\v{c}} \cite{Takac-12}.

In the derivation of {\sc S.~L.\ Heston}'s model \cite{Heston},
It\^{o}'s formula yields the following diffusion equation
(in {\sc Heston}'s original notation)
\begin{equation}
\label{eq:Heston}
  \left( \frac{\partial}{\partial t} + \mathbf{A} \right)
  U(s,v,t) = 0 \quad\mbox{ for }\, s > 0,\ v > 0,\ t < T \,.
\end{equation}
The variables $s$ and $v$, respectively, stand for the values of
the stochastic processes
$(S_t)_{t\geqslant 0}$ and $(V_t)_{t\geqslant 0}$
at a time $t\geq 0$ on a (continuous) path
$\omega\colon [0,\infty)\to (0,\infty)^2$
(that belongs to the underlying probabilistic space $\Omega$), i.e.,
$s = S_t(\omega) > 0$ and $v = V_t(\omega) > 0$.
We call $\mathbf{A}$ the {\em Black\--Scholes\--It\^{o} operator\/}
for the Heston model; it is defined by
\begin{equation}
\label{eq:Ito-oper}
\begin{aligned}
& (\mathbf{A}U)(s,v,t)\eqdef
\\
& {} v\cdot
  \left(
    \frac{1}{2}\, s^2\, \frac{\partial^2 U}{\partial s^2}(s,v,t)
  + \rho\sigma\, s\, \frac{\partial^2 U}{\partial s\;\partial v}(s,v,t)
  + \frac{1}{2}\, \sigma^2\, \frac{\partial^2 U}{\partial v^2}(s,v,t)
  \right)
\\
& {}
  + (r-q)\, s\, \frac{\partial U}{\partial s}(s,v,t)
  + \left[ \kappa (\theta - v) - \lambda(s,v,t)\right]
    \frac{\partial U}{\partial v}(s,v,t)
\\
& {}
  - r\, U(s,v,t)
  \quad\mbox{ for $s>0$, $v>0$, and $t < T$, }
\end{aligned}
\end{equation}
with the following additional quantities (constants) as given data:
the {\it risk free rate of interest\/} $r\in \RR$;
the {\it dividend yield\/} $q\in \RR$;
the {\it instantaneous drift\/} of the stock price returns
$r-q\equiv {}- q_r\in \RR$
(when interpreted under
 the original, ``real\--world'' probability measure);
the {\it volatility\/} $\sigma > 0$ of
the stochastic volatility $\sqrt{v}$;
the {\it correlation\/} $\rho\in (-1,1)$
between the Brownian motions for the stock price and the volatility;
the {\it rate of mean reversion\/} $\kappa > 0$ of
the stochastic volatility $\sqrt{v}$;
the {\it long term variance\/} $\theta > 0$
(called also {\it long\--run variance\/} or {\it long\--run mean level\/})
of the {\it stochastic variance\/} $v$; and
the {\it price of volatility risk\/} $\lambda(s,v,t)\geq 0$,
in \cite{Heston} chosen to be linear,
$\lambda(s,v,t)\equiv \lambda v$ with a constant
$\lambda\equiv \mathrm{const}\geq 0$.
We assume a constant {\it risk free rate of interest\/} $r$ and
a constant {\it dividend yield\/} $q$; hence,
$r-q = {}- q_r$ is
the {\it instantaneous drift of the stock price returns\/}
(under the original probability measure);
All three quantities, $r$, $q$, and $q_r$, may take any real values;
but, typically, one has $0 < r\leq q < \infty$ whence also $q_r\geq 0$.
We refer the reader to the monograph by
{\sc J.~C.\ Hull} \cite[Chapt.~26, pp.\ 599--607]{Hull-book}
and to {\sc S.~L.\ Heston}'s original article \cite{Heston}
for further description of all these quantities.

The diffusion equation \eqref{eq:Heston}
is supplemented first by the following
{\it dynamic boundary condition\/} as $v\to 0+$,
\begin{equation}
\label{bc:Heston}
  \left( \frac{\partial}{\partial t} + \mathbf{B} \right)
  U(s,0,t) = 0 \quad\mbox{ for }\, s > 0,\ t < T \,.
\end{equation}
The {\em boundary operator\/} $\mathbf{B}$ is the trace of
the Black\--Scholes\--It\^{o} operator $\mathbf{A}$ as $v\to 0+$;
it corresponds to the Black\--Scholes operator with zero volatility:
\begin{equation}
\label{eq:Ito-oper_BC}
\begin{aligned}
& (\mathbf{B}U)(s,0,t)\eqdef
\\
&   (r-q)\, s\, \frac{\partial U}{\partial s}(s,0,t)
  + \kappa\theta\, \frac{\partial U}{\partial v}(s,0,t)
  - r\, U(s,0,t)
\\
& \quad\mbox{ for $s>0$, $v=0$, and $-\infty < t < T$. }
\end{aligned}
\end{equation}
The original {\em Heston boundary conditions\/} in \cite{Heston},
\begin{equation}
\label{bc:Ito-oper}
\left\{
\begin{alignedat}{2}
  U(0,v,t)
& = 0 \quad\mbox{ for $v>0$; }
\\
  \lim_{s\to \infty}\,
  \frac{\partial}{\partial s}\, ( U(s,v,t) - s )
& = 0 \quad\mbox{ for $v>0$; }
\\
  \lim_{v\to \infty}\, ( U(s,v,t) - s )
& = 0 \quad\mbox{ for $s>0$, }
\end{alignedat}
\right.
\end{equation}
at all times $t\in (-\infty,T)$,
seem to be ``economically'' motivated.
Mathematically, one may attempt to motivate them by
the asymptotic behavior of the solution
\begin{math}
  U_{\mathrm{BS}}(s,t)\equiv U_{\mathrm{BS}}(s,v_0,t)
\end{math}
to the Black\--Scholes equation, for $s > 0$ and $t\leq T$, where
the variance $v_0 = \sigma_0^2 > 0$ is a given constant
determined from the constant volatility $\sigma_0 > 0$.
What we mean are the following {\em boundary conditions\/},
\begin{equation}
\label{bc_BS:Ito-oper}
\left\{
\begin{alignedat}{2}
  U_{\mathrm{BS}}(0,v,t)
& = 0 \quad\mbox{ for $v>0$; }
\\
  \lim_{s\to \infty}\,
  \frac{\partial}{\partial s}\, ( U_{\mathrm{BS}}(s,v,t) - s )
& = 0 \quad\mbox{ for $v>0$; }
\\
  \lim_{v\to \infty}\, ( U_{\mathrm{BS}}(s,v,t) - s )
& = 0 \quad\mbox{ for $s>0$, }
\end{alignedat}
\right.
\end{equation}
at all times $t\in (-\infty,T)$.
Roughly speaking, the difference
$U(s,v,t) - U_{\mathrm{BS}}(s,v,t)$
becomes asymptotically small near the boundary,
and so does its $s$\--partial derivative as $s\to \infty$.
The terminal condition as $t\to T-$ for both solutions,
$U$ and $U_{\mathrm{BS}}$, is the {\it payoff function\/}
$\hat{h}(s,v) = (s-K)^{+}$ for $s>0$,
\begin{equation*}
  U(s,v,T) = U_{\mathrm{BS}}(s,v,T) = (s-K)^{+} \,.
\end{equation*}
The solution $U_{\mathrm{BS}}(s,t)$ of the Black\--Scholes equation
has been calculated explicitly in the original article by
{\sc F.\ Black} and {\sc M.\ Scholes} \cite{Black-Scholes};
see also
{\sc J.-P.\ Fouque}, {\sc G.\ Papanicolaou}, and {\sc K.~R.\ Sircar}
\cite[{\S}1.3.4, p.~16]{FouqPapaSir}.

Finally, the diffusion equation \eqref{eq:Heston}
is supplemented also by the following
{\it terminal condition\/} as $t\to T-$,
which is given by the payoff function
$\hat{h}(s,v) = (s-K)^{+}$,
\begin{equation}
\label{term:Heston}
  U(s,v,T) = (s-K)^{+} \quad\mbox{ for }\, s > 0,\ v > 0 \,.
\end{equation}
We would like to point out that, by our mathematical approach,
we are able to treat much more general terminal conditions
\begin{math}
  U(s,v,T) = u_0(s,v) \,\mbox{ for }\, s > 0,\ v > 0 \,;
\end{math}
see Proposition~\ref{prop-Lions} and Theorem~\ref{thm-Main}
in Section~\ref{s:Main} below, where
$u_0\in H$ -- a weighted $L^2$\--type Lebesgue space.
Hence, we are not restricted to
European call options \eqref{term:Heston}.
The terminal\--boundary value problem for eq.~\eqref{eq:Heston}
with the boundary conditions
\eqref{bc:Heston} and \eqref{bc:Ito-oper}, as it stands,
poses a {\it mathematically\/} challenging problem, in particular,
due to the degeneracies in the diffusion part of the operator
$\mathbf{A}$:
Some or all of the coefficients of the second partial derivatives
tend to zero as $s\to 0+$ and/or $v\to 0+$,
making the diffusion effects disappear on the boundary
$\{ (s,0)\colon s>0\}$, cf.\ eq.\ \eqref{eq:Ito-oper_BC}.
Similar questions concerned with terminal and boundary conditions
are addressed in
{\sc E.\ Ekstr\"om} and {\sc J. Tysk} \cite{Ekstrom-Tysk}.
However, their work treats only smooth solutions with
only smooth terminal data and, thus, excludes
the (very basic) European call and put options.

This article is organized as follows.
We begin with a rigorous mathematical formulation of
the Heston model in Section~\ref{s:reform}.
We make use of weighted Lebesgue and Sobolev spaces
originally introduced in
{\sc P.\ Daskalopoulos} and {\sc P.~M.~N.\ Feehan} \cite{Daska-Feehan-14}
and \cite[Sect.~2, p.~5048]{Daska-Feehan-16} and
{\sc P.~M.~N.\ Feehan} and {\sc C.~A.\ Pop} \cite{Feehan-Pop-14}.
An extension of the problem from the real to a complex domain
is formulated in Section~\ref{s:prelim}.
Our main results,
Proposition~\ref{prop-Lions} and Theorem~\ref{thm-Main},
are stated in Section~\ref{s:Main}.
Before giving the proofs of these two results,
in Section~\ref{s:Appl} we present an application of them
to {\sc S.~L.\ Heston}'s model \cite{Heston}
for {\em European call options\/} in Mathematical Finance.
There we also provide an affirmative answer
(Theorem~\ref{thm-complete})
to the problem of {\em market completeness\/} as described in
{\sc M.~H.~A.\ Davis} and {\sc J.\ Ob{\l}{\'o}j} \cite{Davis-Obloj}.
Our contribution to market completeness is also
an alternative definition for a market to be complete
(Definition~\ref{def-complete})
which is based on classical concepts of differential calculus
({\sc I.\ Bajeux\--Besnainou} and {\sc J.-Ch.\ Rochet}
 \cite[p.~12]{Bajeux-Rochet})
rather than on probability theory.
In addition, we discuss the important {\it Feller condition\/}
in Remark~\ref{rem-complete}
and also mention another application to a related model
in Remark~\ref{rem-compl_other}.
The proofs of our main results from Section~\ref{s:Main}
are gradually developed in
Sections \ref{s:Heston-real} through \ref{s:L^2-bound}
and completed in Section~\ref{s:End}.
Finally, Appendix~\ref{s:Trace,Sobolev} contains some technical
asymptotic results for functions from our weighted Sobolev spaces,
whereas Appendix~\ref{s:density} is concerned with
the density of certain analytic functions in these spaces.


\section{Formulation of the mathematical problem}
\label{s:reform}

In this section we introduce {S.~L.\ Heston}'s model
\cite[Sect.~1, pp.\ 328--332]{Heston}
and formulate the associated Cauchy problem as
an evolutionary equation of (degenerate) parabolic type.

\subsection{Heston's stochastic volatility model}
\label{ss:Heston}

We consider the {\sc Heston\/} model given under
a \emph{risk neutral measure\/} via
equations $(1)-(4)$ in \cite[pp.\ 328--329]{Heston}.
The model is defined on a filtered probability space
$(\Omega, \mathcal{F}, (\mathcal{F}_t)_{t\geqslant 0}, \mathbb{P})$,
where $\mathbb{P}$ is a risk neutral probability measure, and
the filtration $(\mathcal{F}_t)_{t\geqslant 0}$
satisfies the usual conditions.
Recalling that $S_t$ denotes the {\em stock price\/} and
$V_t$ the (stochastic) {\em variance\/} of the stock market
at (the real) time $t\geq 0$,
the unknown pair $(S_t,V_t)_{t\geqslant 0}$ satisfies
the following system of {\em stochastic\/} differential equations,
\begin{equation}
\label{e:SVmodel}
\left\{
\begin{aligned}
  \frac{\td S_t}{S_t} &= {}- q_r \td t + \sqrt{V_t} \td W_t \,,
\\
  \td V_t &= \kappa\, (\theta - V_t) \td t
           + \sigma\, \sqrt{V_t}\, \td Z_t \,,
\end{aligned}
\right.
\end{equation}
where $(W_t)_{t\geqslant 0}$ and $(Z_t)_{t\geqslant 0}$
are two Brownian motions with the {\it correlation coefficient\/}
$\rho\in (-1,1)$, a constant given by
$\td \langle W,Z\rangle_t = \rho\td t$.
This is the original {\sc Heston} system in \cite{Heston}.

If $X_t = \ln (S_t/K)$ denotes the (natural) logarithm of
the {\em scaled stock price\/} $S_t / K$ at time $t\geq 0$,
relative to the strike price $K>0$,
then the pair $(X_t,V_t)_{t\geqslant 0}$ satisfies
the following system of stochastic differential equations,
\begin{equation}
\label{e:XVmodel}
\left\{
\begin{aligned}
  \td X_t &= {}- \left( q_r + \genfrac{}{}{}1{1}{2} V_t \right) \td t
           + \sqrt{V_t} \td W_t \,,
\\
  \td V_t &= \kappa\, (\theta - V_t) \td t
           + \sigma\, \sqrt{V_t}\, \td Z_t \,.
\end{aligned}
\right.
\end{equation}

Following \cite[Sect.~4]{Davis-Obloj}, let us consider
a European call option written in this market with {\em payoff\/}
$\hat{h}(S_T,V_T)\equiv \hat{h}(S_T)\geq 0$ at maturity $T>0$, where
$\hat{h}(s) = (s-K)^{+}$ for all $s>0$.
As usual, for $x\in \RR$ we abbreviate
$x^{+}\eqdef \max\{ x,\,0\}$ and $x^{-}\eqdef \max\{ -x,\,0\}$.
Recalling {\sc Heston}'s notation in eqs.\
\eqref{eq:Heston} and \eqref{eq:Ito-oper}, we denote
$x = X_t(\omega)\in \RR$.
We set
$h(x,v)\equiv h(x) = K\, (\ee^x - 1)^{+}$
for all $x = \ln (s/K)\in \mathbb{R}$, so that
$h(x) = \hat{h}(s) = \hat{h}(K\ee^x)$ for $x\in \RR$.
Hence, if the instant values
$(X_t(\omega) ,\, V_t(\omega)) = (x,v)\in \HH$
are known at time $t\in (0,T)$, where
$\mathbb{H} = \mathbb{R}\times (0,\infty) \subset \RR^2$,
the \emph{arbitrage\--free price\/}
$P^h_t$ of the European call option at this time is given by
the following expectation formula
(with respect to the risk neutral probability measure $\mathbb{P}$)
which is justified in \cite{Davis-Obloj} and \cite{Takac-12}:
$P^h_t = p(X_t,V_t,t)$ where
\begin{equation}
\label{e:optionprice}
\begin{aligned}
  p(x,v,t)
& = \mathrm{e}^{-r(T-t)}\, \mathbb{E}_{\mathbb{P}}
    \left[ \hat{h}(S_T)\mid \mathcal{F}_t \right]
  = \mathrm{e}^{-r(T-t)}\, \mathbb{E}_{\mathbb{P}}
    \left[ h(X_T)\mid \mathcal{F}_t \right]
\\
& = \mathrm{e}^{-r(T-t)}\, \mathbb{E}_{\mathbb{P}}
    \left[ h(X_T)\mid X_t = x ,\ V_t = v \right] \,.
\end{aligned}
\end{equation}
Furthermore, $p$ solves the ({\it terminal\/} value) Cauchy problem
\begin{equation}
\label{e:IVP.p}
\left\{
\begin{alignedat}{2}
  \frac{\partial p}{\partial t} + \mathcal{G}_t\, p - rp &= 0 \,,\quad
&&(x,v,t)\in \HH\times (0,T) \,;
\\
  p(x,v,T) &= h(x) \,,\quad &&(x,v)\in \HH \,,
\end{alignedat}
\right.
\end{equation}
with $\mathcal{G}_t$ being
the (time\--independent) infinitesimal generator of
the time\--homogeneous Markov process $(X_t,V_t)$;
cf.\ {\sc A.\ Friedman} \cite[Chapt.~6]{Friedman-75}
or {\sc B.\ {\O}ksendal} \cite[Chapt.~8]{Oksendal}.
Indeed, first, eq.~\eqref{eq:Heston} is derived from
eqs.\ \eqref{e:XVmodel} and \eqref{e:optionprice}
by It\^{o}'s formula, then
the diffusion equation \eqref{e:IVP.p} is obtained from
eq.~\eqref{eq:Heston} using
\begin{align*}
  s
& = K\mathrm{e}^x \,,\quad
  \frac{\mathrm{d}s}{\mathrm{d}x}
  = s \,,
\\
    p(x,v,t)
& = U(s,v,t) \,,\quad
    \frac{\partial p}{\partial x}(x,v,t)
  = s\, \frac{\partial U}{\partial s}(s,v,t) \,,
\\
    \frac{\partial^2 p}{\partial x^2}(x,v,t)
& = s\, \frac{\partial U}{\partial s}(s,v,t)
  + s^2\, \frac{\partial^2 U}{\partial s^2}(s,v,t)
\\
& = \frac{\partial p}{\partial x}(x,v,t)
  + s^2\, \frac{\partial^2 U}{\partial s^2}(s,v,t) \,.
\end{align*}
Hence, the function
$\overline{p}\colon (x,v,t)\mapsto p(x,v,T-t)$
verifies a linear Cauchy problem of the following type, with the notation
$\mathbf{x} = (x_1,x_2)\equiv (x,v)\in \mathbb{H}$,
\begin{equation}
\label{e:div_Cauchy}
\left\{
\begin{aligned}
    \frac{ \partial\overline{p} }{\partial t}
  - \sum_{i,j=1}^2 a_{ij}(\mathbf{x}, t)\,
      \frac{ \partial^2\overline{p} }{\partial x_i\, \partial x_j}
  - \sum_{j=1}^2 b_j(\mathbf{x}, t)\,
      \frac{ \partial\overline{p} }{\partial x_j}
  - c(\mathbf{x}, t)\, \overline{p}
\\
\begin{alignedat}{2}
  =&\; f(\mathbf{x}, t)
&&  \quad\mbox{ for } (\mathbf{x}, t)\in \mathbb{H}\times (0,T) \,;
\\
  \overline{p}(\mathbf{x}, 0) =&\; u_0(\mathbf{x})
&&  \quad\mbox{ for } \mathbf{x}\in \mathbb{H} \,,
\end{alignedat}
\end{aligned}
\right.
\end{equation}
with the function $f(\mathbf{x}, t)\equiv 0$ on the right\--hand side
(which may become nontrivial in related Cauchy problems later on),
the initial data
$u_0(\mathbf{x}) = u_0(x,v) = p(x,v,T) = h(x)$
at $t=0$, and the coefficients
\begin{align*}
&
\begin{aligned}
  a(x,v,t) &= \frac{v}{2}\left(
\begin{array}{cc}
  1          & \rho\sigma\\
  \rho\sigma & \sigma^2
\end{array}
  \right) \in \RR^{2\times 2}_{\mathrm{sym}} \,,
\end{aligned}
\\
&
\begin{aligned}
  b(x,v,t) &= \left(
\begin{array}{c}
  {}- q_r - \genfrac{}{}{}1{1}{2} v\\
  \kappa\, (\theta - v) - \lambda(x,v,T-t)
\end{array}
  \right) \in \RR^2 \,,
\quad
  c(x,v,t) &= - r\in \RR \,,
\end{aligned}
\end{align*}
where the variable $\mathbf{x} = (x_1,x_2)\in \RR^2$
has been replaced by $(x,v)\in \HH\subset \RR^2$.
We have also replaced the meaning of the temporal variable $t$
as real time ($t\leq T$)
by the {\em time to maturity\/} $t$
($t\geq 0$), so that the real time has become $\tau = T-t$.
According to
{\sc S.~L.\ Heston} \cite[eq.~(6), p.~329]{Heston},
the unspecified term $\lambda(x,v,T-t)$ in the vector $b(x,v,t)$
represents the {\it\bfseries price of volatility risk\/}
and is specifically chosen to be
$\lambda(x,v,T-t)\equiv \lambda v$ with a constant $\lambda\geq 0$.
As we have already pointed out in the Introduction
(Section~\ref{s:Intro}),
we can treat much more general terminal conditions
$u_0(\mathbf{x}) = u_0(x,v) = p(x,v,T) = h(x,v)$
than just those corresponding to a European call option,
$p(x,v,T) = h(x) =  K\, (\ee^x - 1)^{+}$ for $(x,v)\in \HH$;
see Section~\ref{s:Main} below.
In particular, we do not need the convexity of the function
$h(x) =  K\, (\ee^x - 1)^{+}$ of $x\in \RR$ used heavily in
{\sc M.\ Romano} and {\sc N.\ Touzi} \cite{RomanoTouzi}.

Next, we eliminate the constants $r\in \RR$ and $\lambda\geq 0$,
respectively, from eq.~\eqref{e:div_Cauchy} by substituting
\begin{equation}
\label{e:U-->u}
  p^{*}(x,v,t)\eqdef \ee^{-r(T-t)}\, \overline{p}(x,v,t)
                   = \ee^{-r(T-t)}\, p(x,v,T-t)
  \quad\mbox{ for }\quad \overline{p} (x,v,t) \,,
\end{equation}
which is the {\it discounted option price\/}, and replacing
$\kappa$ by $\kappa^{*} = \kappa + \lambda > 0$ and
$\theta$ by
\begin{math}
  \theta^{*} = \frac{\kappa\theta}{\kappa + \lambda} > 0 .
\end{math}
Hence, we may set $r = \lambda = 0$.
Finally, we introduce also the re\--scaled variance
$\xi = v / \sigma > 0$ for $v\in (0,\infty)$
and abbreviate
$\theta_{\sigma}\eqdef \theta / \sigma\in \RR$.
These substitutions will have a simplifying effect
on our calculations later.
Eq.~\eqref{e:div_Cauchy} then yields
the following initial value problem for the unknown function
$u(x,\xi,t) = p^{*}(x, \sigma\xi ,t)$:
\begin{equation}
\label{e:Cauchy}
\left\{
\begin{alignedat}{2}
  \frac{\partial u}{\partial t} + \mathcal{A} u &= f(x,\xi,t)
  &&\quad\mbox{ in }\, \HH\times (0,T) \,;
\\
  u(x,\xi,0) &= u_0(x,\xi)
  &&\quad\mbox{ for }\, (x,\xi)\in \HH \,,
\end{alignedat}
\right.
\end{equation}
with the function $f(x,\xi,t)\equiv 0$ on the right\--hand side
and the initial data $u_0(x,\xi)\equiv h(x)$ at $t=0$,
where the (autonomous linear) {\it\bfseries Heston operator\/}
$\mathcal{A}$, derived from eq.~\eqref{e:div_Cauchy},
takes the following form,
\begin{align}
\label{e:Heston-oper}
&
\begin{aligned}
  (\mathcal{A}u)(x,\xi) \eqdef
& {} - \frac{1}{2}\, \sigma\xi\cdot
  \left(
    \frac{\partial^2 u}{\partial x^2}(x,\xi)
  + 2\rho\, \frac{\partial^2 u}{\partial x\;\partial\xi}(x,\xi)
  + \frac{\partial^2 u}{\partial\xi^2}(x,\xi)
  \right)
\\
& {}
  + \left( q_r + \genfrac{}{}{}1{1}{2} \sigma\xi\right)
    \cdot \frac{\partial u}{\partial x}(x,\xi)
  - \kappa (\theta_{\sigma} - \xi)
    \cdot \frac{\partial u}{\partial\xi}(x,\xi)
\end{aligned}
\\
\nonumber
&
\begin{alignedat}{2}
& \equiv\;
&&{} - \frac{1}{2}\, \sigma\xi\cdot
  \left( u_{xx} + 2\rho\, u_{x\xi} + u_{\xi\xi}\right)
\\
&
&&{} + \left( q_r + \genfrac{}{}{}1{1}{2} \sigma\xi\right) \cdot u_x
  {} - \kappa (\theta_{\sigma} - \xi) \cdot u_{\xi}
  \quad\mbox{ for $(x,\xi)\in \HH$. }
\end{alignedat}
\end{align}
Recall $\theta_{\sigma} = \theta / \sigma$.
We prefer to use
the following asymmetric ``divergence'' form of $\mathcal{A}$,
\begin{align}
\label{eq:Heston-oper}
&
\begin{aligned}
  (\mathcal{A}u)(x,\xi) =
& {} - \frac{1}{2}\, \sigma\xi\cdot
  \left[ \frac{\partial}{\partial x}
  \left(
    \frac{\partial u}{\partial x}(x,\xi)
  + 2\rho\, \frac{\partial u}{\partial\xi}(x,\xi)
  \right)
  + \frac{\partial^2 u}{\partial\xi^2}(x,\xi)
  \right]
\\
& {}
  + \left( q_r + \genfrac{}{}{}1{1}{2} \sigma\xi\right)
    \cdot \frac{\partial u}{\partial x}(x,\xi)
  - \kappa (\theta_{\sigma} - \xi)
    \cdot \frac{\partial u}{\partial\xi}(x,\xi)
\end{aligned}
\\
\nonumber
&
\begin{alignedat}{2}
& \equiv\;
&&{} - \frac{1}{2}\, \sigma\xi\cdot
  \left[
  \left( u_x + 2\rho\, u_{\xi}\right)_x + u_{\xi\xi}
  \right]
\\
&
&&{} + \left( q_r + \genfrac{}{}{}1{1}{2} \sigma\xi\right) \cdot u_x
  {} - \kappa (\theta_{\sigma} - \xi) \cdot u_{\xi}
  \quad\mbox{ for $(x,\xi)\in \HH$. }
\end{alignedat}
\end{align}

The {\it\bfseries boundary operator\/}
defined in eq.~\eqref{eq:Ito-oper_BC} transforms the left\--hand side of
eq.~\eqref{bc:Heston} into the following (logarithmic) form
on the boundary $\partial\HH = \RR\times \{ 0\}$ of $\HH$:
\begin{equation}
\label{log:Ito-oper_BC}
\begin{aligned}
& \ee^{-r\tau}\,
  \left( \frac{\partial}{\partial\tau} + \mathbf{B} \right)
  U(s,0,\tau) \Big\vert_{\tau = T-t}
  = {}- \left( \frac{\partial}{\partial t} + \mathcal{B} \right)
  u(x,0,t)
\\
& = {}
  - \frac{\partial u}{\partial t}(x,0,t)
  - q_r\, \frac{\partial u}{\partial x}(x,0,t)
  + \kappa\theta_{\sigma}\,
    \frac{\partial u}{\partial\xi}(x,0,t)
\\
& \quad\mbox{ for $x\in \RR$ and $0 < t < \infty$. }
\end{aligned}
\end{equation}
The remaining {\it\bfseries boundary conditions\/}
\eqref{bc:Ito-oper} become
\begin{equation}
\label{log-bc:Ito-oper}
\left\{
\begin{alignedat}{2}
  u(-\infty,\xi,t)\eqdef
  \lim_{x\to -\infty}\,
  \left( u(x,\xi,t) - K\mathrm{e}^{x-r(T-t)} \right)
& = 0 \quad\mbox{ for $\xi > 0$; }
\\
  \lim_{x\to +\infty}\,
  \left[
  \ee^{-x}\cdot
  \frac{\partial}{\partial x}\,
  \left( u(x,\xi,t) - K\mathrm{e}^{x-r(T-t)} \right)
  \right]
& = 0 \quad\mbox{ for $\xi > 0$; }
\\
  \lim_{\xi\to \infty}\,
  \left( u(x,\xi,t) - K\mathrm{e}^{x-r(T-t)} \right)
& = 0 \quad\mbox{ for $x\in \RR$, }
\end{alignedat}
\right.
\end{equation}
at all times $t\in (0,\infty)$.

In the next paragraph we give a definition of $\mathcal{A}$
as a densely defined, closed linear operator in a Hilbert space.

\subsection{Weak formulation in a weighted $L^2$\--space}
\label{ss:weak_form}

Now we formulate the initial\--boundary value problem
for eq.~\eqref{eq:Heston} with the boundary conditions
\eqref{bc:Heston} and \eqref{bc:Ito-oper}
in a weighted $L^2$\--space.
In the context of the Heston model,
similar weighted Lebesgue and Sobolev spaces were used earlier in
{\sc P.\ Daskalopoulos} and {\sc P.~M.~N.\ Feehan} \cite{Daska-Feehan-14}
and \cite[Sect.~2, p.~5048]{Daska-Feehan-16} and
{\sc P.~M.~N.\ Feehan} and {\sc C.~A.\ Pop} \cite{Feehan-Pop-14}.
To this end, we wish to consider the Heston operator $\mathcal{A}$,
defined in eq.~\eqref{eq:Heston-oper} above,
as a densely defined, closed linear operator in
the weighted Lebesgue space
$H\equiv L^2(\mathbb{H};\mathfrak{w})$,
where the weight
$\mathfrak{w}\colon \mathbb{H}\to (0,\infty)$ is defined by
\begin{equation}
\label{def:w}
  \mathfrak{w}(x,\xi)\eqdef
  \xi^{\beta - 1}\, \ee^{ - \gamma |x| - \mu\xi }
  \quad\mbox{ for $(x,\xi)\in \HH$, }
\end{equation}
and $H = L^2(\mathbb{H};\mathfrak{w})$
is the {\em complex\/} Hilbert space endowed with
the inner product
\begin{equation}
\label{def:w_prod}
  (u,w)_{H}\equiv
  (u,w)_{ L^2(\mathbb{H};\mathfrak{w}) } \eqdef
  \int_{\mathbb{H}} u\, \bar{w}\cdot \mathfrak{w}(x,\xi)
    \,\mathrm{d}x \,\mathrm{d}\xi
    \quad\mbox{ for }\, u,w\in H \,.
\end{equation}
Here, $\beta, \gamma, \mu\in (0,\infty)$
are suitable positive constants that will be specified later,
in Section~\ref{s:Heston-real}
(see also Appendix~\ref{s:Trace,Sobolev}).
However, it is already clear that if we want that the weight
$\mathfrak{w}(x,\xi)$ tends to zero as $\xi\to 0+$,
we have to assume $\beta > 1$.
Similarly, if we want that the initial condition
$u_0(x,\xi) = K (\mathrm{e}^x - 1)^{+}$ for $(x,\xi)\in \HH$
belongs to $H$, we must require $\gamma > 2$.

We prove in Section~\ref{s:Heston-real}, {\S}\ref{ss:bound-R},
that the sesquilinear form associated to $\mathcal{A}$,
\begin{equation*}
  (u,w)\mapsto (\mathcal{A}u, w)_H
       \equiv  (\mathcal{A}u, w)_{ L^2(\mathbb{H};\mathfrak{w}) } \,,
\end{equation*}
is {\it bounded\/} on $V\times V$, where $V$ denotes
the {\em complex\/} Hilbert space
$H^1(\mathbb{H};\mathfrak{w})$ endowed with the inner product
\begin{equation}
\label{def:w_prod-H^1}
\begin{aligned}
  (u,w)_{V}\equiv
  (u,w)_{ H^1(\mathbb{H};\mathfrak{w}) }
& \eqdef \int_{\mathbb{H}}
  \left( u_x\, \bar{w}_x + u_{\xi}\, \bar{w}_{\xi} \right)
  \cdot \xi\cdot \mathfrak{w}(x,\xi) \,\mathrm{d}x \,\mathrm{d}\xi
\\
& + \int_{\mathbb{H}} u\, \bar{w} \cdot \mathfrak{w}(x,\xi)
    \,\mathrm{d}x \,\mathrm{d}\xi
    \quad\mbox{ for }\, u,w\in H^1(\mathbb{H};\mathfrak{w}) \,.
\end{aligned}
\end{equation}
In particular, by Lemmas \ref{lem-Trace} and~\ref{lem-Trace_x}
in the Appendix (Appendix~\ref{s:Trace,Sobolev}),
every function
$u\in V = H^1(\mathbb{H};\mathfrak{w})$
satisfies also the following (natural) {\it zero boundary conditions\/},
\begin{align}
\label{e:trace:v=0}
  \xi^{\beta}\cdot \int_{-\infty}^{+\infty}
    |u(x,\xi)|^2\cdot \ee^{- \gamma |x|} \,\mathrm{d}x
& \,\longrightarrow\, 0 \quad\mbox{ as }\, \xi\to 0+ \,,
\\
\label{e:trace:v=infty}
  \xi^{\beta}\, \ee^{- \mu\xi}\cdot \int_{-\infty}^{+\infty}
    |u(x,\xi)|^2\cdot \ee^{- \gamma |x|} \,\mathrm{d}x
& \,\longrightarrow\, 0 \quad\mbox{ as }\, \xi\to \infty \,,
\end{align}
and
\begin{align}
\label{e:trace:x=+-infty}
  \ee^{- \gamma |x|}\cdot \int_0^{\infty}
    |u(x,\xi)|^2\cdot \xi^{\beta}\, \ee^{- \mu\xi} \,\mathrm{d}\xi
& \,\longrightarrow\, 0 \quad\mbox{ as }\, x\to \pm\infty \,.
\end{align}
(We are no longer using the letter $v = V_t(\omega) > 0$ for variance;
 it has been replaced by the re\--scaled variance
 $\xi = v / \sigma > 0$.)
The following additional {\it vanishing boundary conditions\/}
are determined by our particular {\em realization\/}
of the Heston operator $\mathcal{A}$ with the domain
$V = H^1(\mathbb{H};\mathfrak{w})$, cf.~\eqref{e:Heston-bilin} below:
\begin{align}
\label{bc_xi:Heston-bilin}
\left\{
\begin{aligned}
  \xi^{\beta}\cdot \int_{-\infty}^{+\infty}
    u_{\xi}(x,\xi)\cdot \bar{w}(x,\xi)
    \cdot \ee^{- \gamma |x|} \,\mathrm{d}x
& \,\longrightarrow\, 0 \quad\mbox{ as }\, \xi\to 0+ \,;
\\
  \xi^{\beta}\, \ee^{- \mu\xi}\cdot \int_{-\infty}^{+\infty}
    u_{\xi}(x,\xi)\cdot \bar{w}(x,\xi)
    \cdot \ee^{- \gamma |x|} \,\mathrm{d}x
& \,\longrightarrow\, 0 \quad\mbox{ as }\, \xi\to \infty \,,
\end{aligned}
\right.
\\
\label{bc_x:Heston-bilin}
    \ee^{- \gamma |x|}\cdot \int_0^{\infty}
  ( u_x + 2\rho\, u_{\xi} )\, \bar{w}(x,\xi)\cdot
    \xi^{\beta}\, \ee^{- \mu\xi} \,\mathrm{d}\xi
  \,\longrightarrow\, 0 \quad\mbox{ as }\, x\to \pm\infty \,,
\end{align}
for every function $w\in V$.
The validity of these boundary conditions on the boundary
$\partial\HH = \RR\times \{ 0\}$ of the half\--plane
$\HH = \RR\times (0,\infty)\subset \RR^2$
(i.e., as $\xi\to 0+$)
and as $\xi\to \infty$ is discussed below,
in {\S}\ref{ss:Heston_BC}.
They guarantee that $\mathcal{A}$ is a closed, densely defined
linear operator in the Hilbert space $H$
which possesses a unique extension to a bounded linear operator
$V\to V'$, denoted by $\mathcal{A}\colon V\to V'$ again,
with the property that there is a constant $c\in \RR$ such that
$\mathcal{A} + c\, I$ is {\it coercive\/} on $V$.
Consequently, every function $v\in V$ from the domain
$\mathcal{D}(\mathcal{A}) \subset H$ of $\mathcal{A}$,
\begin{math}
  \mathcal{D}(\mathcal{A}) = \{ v\in V\colon \mathcal{A}v\in H\} ,
\end{math}
must satisfy not only
\eqref{e:trace:v=0}, \eqref{e:trace:v=infty}, and
\eqref{e:trace:x=+-infty} (thanks to $v\in V$),
but also the boundary conditions
\eqref{bc_xi:Heston-bilin} and \eqref{bc_x:Heston-bilin}
(owing to $v\in \mathcal{D}(\mathcal{A})$).
A detailed discussion of all boundary conditions is provided
in {\S}\ref{ss:Heston_BC} below.
The coercivity of $\mathcal{A} + c\, I$ on $V$ will be proved
in Section~\ref{s:Heston-real}, {\S}\ref{ss:coerce-R}.

The sesquilinear form
\begin{math}
  (u,w)\mapsto (\mathcal{A}u, w)_H
\end{math}
is used in the {\em Hilbert space definition\/} of
the linear operator $\mathcal{A}$ by the following procedure.
For any given
$u,w\in H^1(\mathbb{H};\mathfrak{w}) \cap W^{2,\infty}(\HH)$,
we use eq.~\eqref{eq:Heston-oper} to calculate the inner product
\begin{align}
\label{e:Heston-bilin}
\begin{aligned}
& (\mathcal{A}u, w)_H\equiv
  (\mathcal{A}u, w)_{ L^2(\mathbb{H};\mathfrak{w}) } =
\\
& \frac{\sigma}{2}\int_{\HH}
  \left[
  ( u_x + 2\rho\, u_{\xi} )\cdot \bar{w}_x + u_{\xi}\cdot \bar{w}_{\xi}
  \right] \cdot \xi\cdot \mathfrak{w}(x,\xi)
    \,\mathrm{d}x \,\mathrm{d}\xi
\\
& {}+ \frac{\sigma}{2}\int_{\HH}
  \left[
  ( u_x + 2\rho\, u_{\xi} )\, \bar{w}\cdot
    \xi\cdot \partial_x \mathfrak{w}(x,\xi) + u_{\xi}\cdot \bar{w}\cdot
    \partial_{\xi} \bigl( \xi\cdot \mathfrak{w}(x,\xi) \bigr)
  \right]
    \,\mathrm{d}x \,\mathrm{d}\xi
\\
& {}- \frac{\sigma}{2}\int_0^{\infty}
  ( u_x + 2\rho\, u_{\xi} )\, \bar{w}\cdot
      \xi\cdot \mathfrak{w}(x,\xi) \,\mathrm{d}\xi
    \,\Big\vert_{x = -\infty}^{x = +\infty}
\\
& {}- \frac{\sigma}{2}\int_{-\infty}^{+\infty}
      u_{\xi}\cdot \bar{w}\cdot \xi\cdot \mathfrak{w}(x,\xi)
    \,\mathrm{d}x
    \,\Big\vert_{\xi = 0}^{\xi = \infty}
\end{aligned}
\\
\nonumber
\begin{aligned}
& {}- \int_{\HH}
  \left[
{}- \left( q_r + \genfrac{}{}{}1{1}{2} \sigma\xi\right) u_x
  + \kappa (\theta_{\sigma} - \xi)\, u_{\xi}
  \right] \cdot \bar{w}\cdot \mathfrak{w}(x,\xi)
    \,\mathrm{d}x \,\mathrm{d}\xi
\end{aligned}
\end{align}
\begin{align*}
\begin{aligned}
& {}= \frac{\sigma}{2}\int_{\HH}
  \left(
    u_x\cdot \bar{w}_x + 2\rho\, u_{\xi}\cdot \bar{w}_x
  + u_{\xi}\cdot \bar{w}_{\xi}
  \right) \cdot \xi\cdot \mathfrak{w}(x,\xi)
    \,\mathrm{d}x \,\mathrm{d}\xi
\\
& {}+ \frac{\sigma}{2}\int_{\HH}
  \bigl[ {}- \gamma\, \Sgn x\cdot
  ( u_x + 2\rho\, u_{\xi} )\, \bar{w}\cdot \xi
  + (\beta - \mu\xi)\, u_{\xi}\cdot \bar{w}
  \bigr] \mathfrak{w}(x,\xi)
    \,\mathrm{d}x \,\mathrm{d}\xi
\end{aligned}
\\
\begin{aligned}
  {}- \frac{\sigma}{2}
& \left[
      \lim_{x\to +\infty}
    \left( \ee^{- \gamma |x|}\cdot \int_0^{\infty}
  ( u_x + 2\rho\, u_{\xi} )\, \bar{w}\cdot
    \xi^{\beta}\, \ee^{- \mu\xi} \,\mathrm{d}\xi
    \right)
  \right.
\\
& \left.
{}  - \lim_{x\to -\infty}
    \left( \ee^{- \gamma |x|}\cdot \int_0^{\infty}
  ( u_x + 2\rho\, u_{\xi} )\, \bar{w}\cdot
    \xi^{\beta}\, \ee^{- \mu\xi} \,\mathrm{d}\xi
    \right)
  \right]
\end{aligned}
\\
\begin{aligned}
& {}+ \frac{\sigma}{2}
  \left[
      \lim_{\xi\to 0+}
    \left( \xi^{\beta}\cdot \int_{-\infty}^{+\infty}
      u_{\xi}\cdot \bar{w}\cdot \ee^{- \gamma |x|} \,\mathrm{d}x
    \right)
    - \lim_{\xi\to \infty}
    \left( \xi^{\beta}\, \ee^{- \mu\xi}\cdot \int_{-\infty}^{+\infty}
      u_{\xi}\cdot \bar{w}\cdot \ee^{- \gamma |x|} \,\mathrm{d}x
    \right)
  \right]
\end{aligned}
\\
\begin{aligned}
& {}- \int_{\HH}
  \left( {}- q_r\, u_x + \kappa\theta_{\sigma}\, u_{\xi}
  \right) \cdot \bar{w}\cdot \mathfrak{w}(x,\xi)
    \,\mathrm{d}x \,\mathrm{d}\xi
\\
& {}+ \int_{\HH}
  \left( \genfrac{}{}{}1{1}{2}\sigma\, u_x + \kappa\, u_{\xi}
  \right) \bar{w}\cdot \xi\cdot \mathfrak{w}(x,\xi)
    \,\mathrm{d}x \,\mathrm{d}\xi \,,
\end{aligned}
\end{align*}
where we now impose the vanishing boundary conditions
\eqref{bc_xi:Heston-bilin} and \eqref{bc_x:Heston-bilin}.

Hence, the sesquilinear form \eqref{e:Heston-bilin} becomes
\begin{align}
\label{e:Heston-diss:u=w}
&
\begin{aligned}
  (\mathcal{A}u, w)_H
& {}= \frac{\sigma}{2}\int_{\HH}
  \left(
    u_x\cdot \bar{w}_x + 2\rho\, u_{\xi}\cdot \bar{w}_x
  + u_{\xi}\cdot \bar{w}_{\xi}
  \right) \cdot \xi\cdot \mathfrak{w}(x,\xi)
    \,\mathrm{d}x \,\mathrm{d}\xi
\\
& {}+ \frac{\sigma}{2}\int_{\HH}
  (1 - \gamma\, \Sgn x)\, u_x\cdot \bar{w}
    \cdot \xi\cdot \mathfrak{w}(x,\xi)
    \,\mathrm{d}x \,\mathrm{d}\xi
\\
& {}+ \int_{\HH}
  \left( \kappa - \gamma\rho\sigma\, \Sgn x
       - \genfrac{}{}{}1{1}{2} \mu\sigma
  \right) u_{\xi}\cdot \bar{w}
    \cdot \xi\cdot \mathfrak{w}(x,\xi)
    \,\mathrm{d}x \,\mathrm{d}\xi
\end{aligned}
\\
\nonumber
&
\begin{aligned}
& {}+ q_r
    \int_{\HH} u_x\cdot \bar{w}\cdot \mathfrak{w}(x,\xi)
    \,\mathrm{d}x \,\mathrm{d}\xi
    + \left( \genfrac{}{}{}1{1}{2} \beta\sigma - \kappa\theta_{\sigma}
      \right)
    \int_{\HH} u_{\xi}\cdot \bar{w}\cdot \mathfrak{w}(x,\xi)
    \,\mathrm{d}x \,\mathrm{d}\xi \,.
\end{aligned}
\end{align}
All integrals on the right\--hand side converge absolutely for any pair
$u,w\in V$; see the proof of our Proposition~\ref{prop-bound-R} below.
In what follows we use the last formula, eq.~\eqref{e:Heston-diss:u=w},
to define the sesquilinear form \eqref{e:Heston-bilin} in $V\times V$.
Of course, in the calculations above we have assumed
the boundary conditions in
\eqref{bc_xi:Heston-bilin} and \eqref{bc_x:Heston-bilin}.

We make use of the {\em Gel'fand triple\/}
$V\hookrightarrow H = H'\hookrightarrow V'$, i.e.,
we first identify the Hilbert space $H$ with its dual space $H'$,
by the Riesz representation theorem,
then use the imbedding $V\hookrightarrow H$,
which is dense and continuous, to construct its adjoint mapping
$H'\hookrightarrow V'$,
a dense and continuous imbedding of $H'$ into the dual space $V'$
of $V$ as well.
The (complex) inner product on $H$ induces
a sesquilinear duality between $V$ and $V'$;
we keep the notation
\begin{math}
  ( \,\cdot\, , \,\cdot\, )_H
\end{math}
also for this duality.



\subsection{The Cauchy problem in the real domain}
\label{ss:Cauchy-real}

Let us return to the initial value problem \eqref{e:Cauchy}.
The letter $T$ stands for an arbitrary (finite) upper bound on time $t$.
The latter, $t$, can still be regarded as time to maturity.

\begin{definition}\label{def-weak_sol}\nopagebreak
\begingroup\rm
Let\/ $0 < T < \infty$, $f\in L^2((0,T)\to V')$, and\/ $u_0\in H$.
A function $u\colon \HH\times [0,T]$ $\to \RR$
is called a {\it weak solution\/} to the initial value problem
\eqref{e:Cauchy}
if it has the following properties:
\begin{itemize}
\item[{\rm (i)}]
the mapping
$t\mapsto u(t)\equiv u(\,\cdot\,, \,\cdot\,, t)\colon [0,T]\to H$
is a continuous function, i.e.,
$u\in C([0,T]\to H)$;
\item[{\rm (ii)}]
the initial value $u(0) = u_0$ in $H$;
\item[{\rm (iii)}]
the mapping
$t\mapsto u(t)\colon (0,T)\to V$
is a B\^ochner square\--integrable function, i.e.,
$u\in L^2((0,T)\to V)$; and
\item[{\rm (iv)}]
for every function
\begin{equation*}
  \phi\in L^2((0,T)\to V)\cap W^{1,2}((0,T)\to V')
  \hookrightarrow C([0,T]\to H) \,,
\end{equation*}
the following equation holds,
\begin{equation}
\label{def:weak_sol}
\begin{aligned}
&   ( u(T), \phi(T) )_H
  - \int_0^T
    \left( u(t), \genfrac{}{}{}1{\partial\phi}{\partial t}(t)
    \right)_H \,\mathrm{d}t
  + \int_0^T (\mathcal{A}u(t), \phi(t))_H \,\mathrm{d}t
\\
& = ( u_0, \phi(0) )_H
  + \int_0^T (f(t), \phi(t))_H \,\mathrm{d}t \,.
\end{aligned}
\end{equation}
\end{itemize}
\endgroup
\end{definition}
\par\vskip 10pt

The following remarks are in order:

First, our definition of a weak solution is equivalent with that given in
{\sc L.~C.\ Evans} \cite[{\S}7.1]{Evans-98}, p.~352.
We are particularly interested in the solution with the initial value
$u_0(x,\xi) = K\, (\mathrm{e}^x - 1)^{+}$ for $(x,\xi)\in \HH$, cf.\
eq.~\eqref{term:Heston}.
Clearly, we have $u_0\in H$
if and only if $\gamma > 2$, $\beta > 0$, and $\mu > 0$.
However, if the European put option with the initial value
$u_0(x,\xi) = K\, (1 - \mathrm{e}^x)^{+}$ for $(x,\xi)\in \HH$
is considered, any small constant $\gamma > 0$ will do.

$W^{1,2}((0,T)\to V')$ denotes the Sobolev space of all functions
$\phi\in L^2((0,T)\to V')$ that possess a distributional time\--derivative
$\phi'\in L^2((0,T)\to V')$.
The norm is defined in the usual way; cf.\
{\sc L.~C.\ Evans} \cite[{\S}5.9]{Evans-98}.
The properties of
$V\equiv H^1(\mathbb{H};\mathfrak{w})$
justify the notation
$V'= H^{-1}(\mathbb{H};\mathfrak{w})$.

The continuity of the imbedding
\begin{equation*}
  L^2((0,T)\to V)\cap W^{1,2}((0,T)\to V')
  \hookrightarrow C([0,T]\to H)
\end{equation*}
is proved, e.g.,
in {\sc L.~C.\ Evans} \cite[{\S}5.9]{Evans-98}, Theorem~3 on p.~287.


\subsection{The Heston operator and boundary conditions}
\label{ss:Heston_BC}

We have seen in our definition of
the sesquilinear form \eqref{e:Heston-diss:u=w}
in paragraph {\S}\ref{ss:weak_form}
that the boundary conditions
\eqref{bc_xi:Heston-bilin} and \eqref{bc_x:Heston-bilin}
are necessary for performing integration by parts to obtain
the sesquilinear form \eqref{e:Heston-diss:u=w}.
They should be valid for every weak solution
$u\colon \HH\times [0,T]\to \RR$ of the initial value problem
\eqref{e:Cauchy} at a.e.\ time $t\in (0,T)$, and for every test function
$w\in V$.
A natural way to satisfy these conditions is to estimate
the absolute value of the integrals from above by Cauchy's inequality
and then impose or verify the following boundary conditions,
\begin{align}
\label{bc_xi:bound_u}
\left\{
\begin{aligned}
  \xi^{\beta}\cdot \int_{-\infty}^{+\infty}
    |u_{\xi}(x,\xi)|^2\cdot \ee^{- \gamma |x|} \,\mathrm{d}x
& \leq \mathrm{const} < \infty \quad\mbox{ as }\, \xi\to 0+ \,;
\\
  \xi^{\beta}\, \ee^{- \mu\xi}\cdot \int_{-\infty}^{+\infty}
    |u_{\xi}(x,\xi)|^2\cdot \ee^{- \gamma |x|} \,\mathrm{d}x
& \leq \mathrm{const} < \infty \quad\mbox{ as }\, \xi\to \infty+ \,,
\end{aligned}
\right.
\\
\label{bc_x:bound_u}
    \ee^{- \gamma |x|}\cdot \int_0^{\infty}
  | u_x + 2\rho\, u_{\xi} |^2\cdot
    \xi^{\beta}\, \ee^{- \mu\xi} \,\mathrm{d}\xi
  \leq \mathrm{const} < \infty \quad\mbox{ as }\, x\to \pm\infty \,,
\end{align}
together with \eqref{e:trace:v=0}, \eqref{e:trace:v=infty}, i.e.,
\begin{align}
\label{bc:bound_w}
\left\{
\begin{aligned}
  \xi^{\beta}\cdot \int_{-\infty}^{+\infty}
    |w(x,\xi)|^2\cdot \ee^{- \gamma |x|} \,\mathrm{d}x
& \,\longrightarrow\, 0 \quad\mbox{ as }\, \xi\to 0+ \,;
\\
  \xi^{\beta}\, \ee^{- \mu\xi}\cdot \int_{-\infty}^{+\infty}
    |w(x,\xi)|^2\cdot \ee^{- \gamma |x|} \,\mathrm{d}x
& \,\longrightarrow\, 0 \quad\mbox{ as }\, \xi\to \infty \,,
\end{aligned}
\right.
\end{align}
and \eqref{e:trace:x=+-infty} for $w$ in place of $u$.
In other words, we have
\begin{itemize}
\item
$\;$
\eqref{bc_xi:bound_u} and \eqref{bc:bound_w} $\,\Rightarrow\,$
\eqref{bc_xi:Heston-bilin} $\quad$ whereas $\quad$
\eqref{bc_x:bound_u} and \eqref{e:trace:x=+-infty} $\,\Rightarrow\,$
\eqref{bc_x:Heston-bilin}.
\end{itemize}
Indeed, by Lemma~\ref{lem-Trace},
the latter boundary conditions, \eqref{bc:bound_w},
are satisfied for every test function $w\in V$.
Similarly, \eqref{e:trace:x=+-infty} holds by Lemma~\ref{lem-Trace_x}.
We stress that only the boundary conditions in
\eqref{bc_xi:bound_u} and \eqref{bc_x:bound_u} are {\it imposed\/};
they do {\it not\/} follow from $u\in V$.

Two of these boundary conditions on the boundary
$\partial\HH = \RR\times \{ 0\}$ of the half\--plane
$\HH = \RR\times (0,\infty)\subset \RR^2$
limit from above the growth of the solution $u(x,\xi)$
at an arbitrarily low volatility level
$\sqrt{\xi}$, i.e., as the variance $\xi\to 0+$.

From now on, we use exclusively formula \eqref{e:Heston-diss:u=w}
to define the linear operator $\mathcal{A}\colon V\to V'$
that appears in the sesquilinear form \eqref{e:Heston-bilin}
obtained directly for the Heston operator \eqref{eq:Heston-oper}.
This means that we no longer need the boundary conditions in
\eqref{bc_xi:bound_u} and \eqref{bc_x:bound_u}
(or in \eqref{bc_xi:Heston-bilin} and \eqref{bc_x:Heston-bilin})
imposed on $u\in V$.

We refer the reader to the recent work by
{\sc P.~M.~N.\ Feehan} \cite{Feehan-13},
Appendix~B, {\S}B.1, pp.\ 57--58,
for numerous interesting properties of $\mathcal{A}$.

\begin{remark}\label{rem-prop-Feller}\nopagebreak
{\rm (Coercivity conditions.)}$\;$
\begingroup\rm
It is important to remark at this stage of our investigation of
the Heston operator $\mathcal{A}$ that,
in order to ensure the coercivity of $\mathcal{A} + c\, I$ on $V$,
one has to assume the well\--known {\it\bfseries Feller condition\/}
(\cite{Feller, Guo-Grzel-Ooster}),
\begin{equation}
\label{e:Feller}
  \genfrac{}{}{}1{1}{2} \sigma^2 - \kappa\theta < 0 \,.
\end{equation}

However, {\it Feller's condition\/} \eqref{e:Feller}
is {\it not sufficient\/} for obtaining the desired coercivity.
We need to guarantee also
\begin{equation*}
  c_1' = \genfrac{}{}{}1{1}{2} \sigma
    \left[
    \left( \frac{\kappa}{\sigma} - \gamma\, |\rho| \right)^2
  - \gamma (1+\gamma)
    \right] \geq 0 \,;
\end{equation*}
cf.\ ineq.\ \eqref{e:Heston-diss:RE}
in the proof of Proposition~\ref{prop-coerc-R} below.
That is, we need to assume
\begin{equation}
\label{ineq:c_1'>0}
  \kappa\geq \sigma
    \left( \gamma\, |\rho| + \sqrt{ \gamma (1+\gamma) } \right)
  \quad \left(\, > \sigma\gamma (|\rho| + 1) \,\right) \,.
\end{equation}

The last inequality is an additional condition to
{\it Feller's condition\/},
$\genfrac{}{}{}1{1}{2} \sigma^2 - \kappa\theta < 0$,
both of them requiring the rate of mean reversion $\kappa > 0$
of the stochastic volatility in system~\eqref{e:SVmodel}
to be sufficiently large.
This additional condition is caused by the fact that
{\sc W.\ Feller} \cite{Feller} considers only
an analogous problem in one space dimension ($\xi\in \RR_+$),
so that the solution $u = u(\xi)$ is independent from $x\in \RR$.
In particular, if the initial condition
$u_0 = u(\,\cdot\,, \,\cdot\,, 0)\in H$ for $u(x,\xi,t)$
permits us to take $\gamma > 0$ arbitrarily small, then
inequality \eqref{ineq:c_1'>0} is easily satisfied,
provided {\it Feller's condition\/}
$\genfrac{}{}{}1{1}{2} \sigma^2 - \kappa\theta < 0$
is satisfied.
This is the case for the European put option with the initial condition
$u_0(x,\xi) = K\, (1 - \mathrm{e}^x)^{+}$ (${}\leq K$)
for $(x,\xi)\in \HH$.
However, if we wish to accommodate also initial conditions of type
$u_0(x,\xi) = K\, (\mathrm{e}^x - 1)^{+}$ for $(x,\xi)\in \HH$,
then we are forced to take $\gamma > 2$ to ensure that $u_0\in H$.

We refer the reader to the recent monograph by
{\sc G.~H. Meyer} \cite{Meyer-2015}
for a discussion of the role of Feller's condition
in the boundary conditions in Heston's model.
\hfill\Square
\endgroup
\end{remark}
\par\vskip 10pt

We will see in Section~\ref{s:Main} that the initial value problem
\eqref{e:Cauchy} has a unique weak solution
$u\colon \HH\times [0,T]\to \RR$.
Recall that, by eq.~\eqref{term:Heston},
we are particularly interested in the solution with the initial value
$u_0(x,\xi) = K\, (\mathrm{e}^x - 1)^{+}$ for $(x,\xi)\in \HH$.
We are not able to show that even this particular solution satisfies
{\sc Heston}'s boundary conditions
\eqref{bc:Heston} and \eqref{log-bc:Ito-oper}.
However, the asymptotic boundary conditions in \eqref{log-bc:Ito-oper}
are taken into account by the choice of function spaces $H$ and $V$.
{\sc Heston}'s boundary operator \eqref{log:Ito-oper_BC}
assumes the existence of traces of certain functions of
$(x,\xi)$ as $\xi\to 0+$ which have to satisfy
a partial differential equation derived from \eqref{bc:Heston}.
In conditions \eqref{e:trace:x=+-infty} and \eqref{bc:bound_w}
we assume only that some of the functions in
the boundary operator \eqref{log:Ito-oper_BC}
do not blow up too fast as $\xi\to 0+$.


\section{The complex domain: Preliminaries and notation}
\label{s:prelim}

We complexify the real space\--time domain
$\HH\times (0,\infty)$ as follows:

We denote by
\begin{equation}
\label{e:strip}
  \mathfrak{X}^{(r)} \eqdef \RR + \ii (-r,r)\subset \CC
\end{equation}
the {\em complex strip\/} of width $2r$, $r\in (0,\infty)$,
which consists of all (complex) numbers
$z = x + \ii y\in \CC$ whose imaginary part, $y = \IM z$,
is bounded by $|y| < r$, while the real part, $x = \RE z$,
may take any value $x\in \RR$
(see Figure~1).
This is the complexification of the variable $x\in \RR$.
The remaining two independent variables, $\xi,t\in (0,\infty)$,
will be complexified by angular domains with the vertex at zero.
We denote by
\begin{equation}
\label{e:angle}
  \Delta_{\vartheta} \eqdef
  \{ \zeta = \varrho\ee^{\ii\theta} \in \CC\colon
     \varrho > 0\mbox{ and } \theta\in (-\vartheta, \vartheta) \}
\end{equation}
the {\em complex angle\/} of angular width $2\vartheta$,
$\vartheta\in (0, \pi / 2)$
(Figure~2).
Notice that the standard logarithm
$\zeta\mapsto z = \log\zeta$ is a conformal mapping from the angle
$\Delta_{\vartheta}$ onto the strip $\mathfrak{X}^{(\vartheta)}$.
Now, given any
$\vartheta_{\xi}, \vartheta_t\in (0, \pi / 2)$,
we complexify $\xi$ as
$\zeta = \xi + \ii\eta\in \Delta_{\vartheta_{\xi}}$,
so that $\xi = \RE\zeta > 0$, and $t$ as
$t = \alpha + \ii\tau\in \Delta_{\vartheta_t}$,
whence $\alpha = \RE t > 0$, thus stressing that
we allow for complex time $t\in \Delta_{\vartheta_t}$
in accordance with the usual notation for
holomorphic $C^0$-semigroups.
The half\--plane $\HH = \RR\times (0,\infty)$
is naturally imbedded into the complex domain
\begin{equation}
\label{x,v:strip}
  \mathfrak{V}^{(r)} \eqdef
  \mathfrak{X}^{(r)} \times \Delta_{\arctan r}\subset \CC^2 \,,
    \quad r\in (0,\infty) \,.
\end{equation}
%


\null\hspace{1.0cm} 
\begin{tikzpicture}[scale=1]

    \draw[->] (-1,0) -- (5.5,0) node[right, below] {$\quad \quad x\in \RR$};
     \draw[->] (0,-1.8) -- (0,1.8) node[left,above] {$\ii y \in \ii\RR\quad\quad$};
           \draw[red] (-1,1) -- (5.5,1);
           \draw[red] (-1,-1) -- (5.5,-1);
           \draw[<->] (1.5,0.05)--(1.5,0.95) node[right,midway]
           {$r(\alpha)$};
           \draw[<->] (1.5,-0.05)--(1.5,-0.95) node[right,midway]
           {$r(\alpha)$};
           \draw[green](5.08,0.6) node {$z = x+\ii y\in \CC$};
           \draw[green,fill=green](3.5,0.6) circle (0.4mm);
           \draw (2.5,-2.5) node { Figure~1.$\;$ Strip
           $\mathfrak{X}^{(r)} = \RR+\ii (-r,r))$ for $r = r(\alpha)$,
           $\alpha > 0$. };
\end{tikzpicture}

\vspace{1cm}

\hspace{2.0cm} 
\begin{tikzpicture}[scale=1]

           \draw[->] (-1,0) -- (5.5,0) node[right, below]
           {$\quad \quad \xi\in (0,+\infty)$};
           \draw[->] (0,-1.8) -- (0,1.8) node[left,above]
           {$\ii\eta\in \ii\RR\quad\quad$};
                \draw[red] (0,0) -- (5.5,2);
                \draw[red] (0,0) -- (5.5,-2);
                \draw[green](5.6,0.6) node {$\zeta = \xi + \ii\eta\in \CC$};
                \draw[green,fill=green](4,0.6) circle (0.4mm);
                \draw [<->](2,0) arc (0:20:2) ;
\draw (12:2.5) node {$\vartheta (\alpha)$};
\draw [<->](2,0) arc (0:-20:2) ;
\draw (-12:2.5) node {$\vartheta (\alpha)$};
                \draw (2.5,-2.5) node { Figure~2.$\;$ Angle
                $\Delta_{\vartheta}$. };
\end{tikzpicture}

\vspace{1cm}

\hspace{1.5cm} 
\begin{tikzpicture}[scale=1]

\draw[->] (-1,0) -- (5.5,0) node[right, below] {$\quad\alpha$};
\draw(4.2,0) node[left,below]{$T$};
\draw[-](4.5,-3) -- (4.5,3);
\draw[->] (0,-1.8) -- (0,1.8) node[left,above] {$\ii \tau \quad$};
\draw[domain=0.001:1.5,samples=100,color=blue, ultra  thick] plot ({\x},{0.6*(\x)});
\draw[domain=1.5:5,samples=100,color=blue, ultra  thick] plot ({\x},{0.6*(\x)});
\draw[domain=0.001:1.5,samples=100,color=blue, ultra thick] plot ({\x},{-0.6*(\x)});
\draw[domain=1.5:5,samples=100,color=blue, ultra  thick] plot ({\x},{-0.6*(\x)});
\draw(2.4,0) node[right,below]{$T'$};
\draw[dashed](2.1,-1.7) -- (2.1,1.7);
\draw(-0.2,0) node[left,below] {$0$};
\draw[thick](2.03,0.5) -- (2.17,0.5) node[right] {$\tau$};
\draw (8.5,0) node { Figure~3.$\;$ $\Sigma^{(\alpha)} (\nu_0)$.};
\end{tikzpicture}

\vspace{1cm}

\hspace{0.0cm} 
\begin{tikzpicture}[scale=1]

\draw[->] (-1,0) -- (5.5,0) node[right, below] {$\quad\alpha$};            
\draw(4.2,0) node[left,below]{$T$};
\draw[-](4.5,-1.7) -- (4.5,1.7);
\draw[->] (0,-1.8) -- (0,1.8); 
\draw[-](-0.08,0.9) -- (0.08,0.9);
\draw[-](-0.08,-0.9) -- (0.08,-0.9);
\draw[->] (0,-1.8) -- (0,1.8) node[left,above] {$\ii y \quad$};
\draw(-1,0.53) node[left,above] {$\kappa_0\cdot \min\{ \alpha, T'\} \quad\quad$};
\draw(-1.15,-0.53) node[left,below] {$-\kappa_0\cdot \min\{ \alpha, T'\} \quad\quad$};
\draw[domain=0.001:2.1,samples=100,color=blue, ultra thick] plot ({\x},{0.6*(\x)});
\draw[domain=2.1:5,samples=100,dashed,color=blue] plot ({\x},{0.6*(\x)});
\draw[domain=0.001:2.1,samples=100,color=blue, ultra thick] plot ({\x},{-0.6*(\x)});
\draw[domain=2.1:5,samples=100,dashed,color=blue] plot ({\x},{-0.6*(\x)});
\draw(2.4,0) node[right,below]{$T'$};
\draw[dashed](2.1,-1.7) -- (2.1,1.7);
\draw(-0.2,0) node[left,below] {$0$};
\draw[blue, ultra thick](2.1,1.26)--(5,1.26);
\draw[blue, ultra thick](2.1,-1.26)--(5,-1.26);
\draw[thick](2.03,0.7) -- (2.17,0.7) node[right] {$y$};
\draw (8.5,0) node { Figure~4.$\;$ $\Gamma^{(T')}_T (\kappa_0,\nu_0)$.};
\end{tikzpicture}

\vspace{1cm}


In order to give a plausible lower estimate on
the space\--time domain of holomorphy
(i.e., the domain of complex analyticity)
of a weak solution $u$ to
the homogeneous initial value problem \eqref{e:Cauchy}
with $f\equiv 0$,
we introduce a few more subsets of $\CC^2\times \CC$
(cf.\ {\sc P.\ Tak\'a\v{c}} et al.\ \cite[p.~428]{TakacBoller}
 or   {\sc P.\ Tak\'a\v{c}} \cite[pp.\ 58--59]{Takac-12}):

The two constants $\kappa_0, \nu_0\in (0,\infty)$
used below will be specified later (in Theorem~\ref{thm-Main});
$0\leq \alpha < \infty$ is an arbitrary number.
First, we set
\begin{align}
\label{e:Pi}
  \mathfrak{V}^{(\kappa_0\alpha)}
& {}=
  \mathfrak{X}^{(\kappa_0\alpha)} \times
  \Delta_{\arctan (\kappa_0\alpha)}
\\
\nonumber
& = \bigl\{
    (z,\zeta) = (x + \ii y, \xi + \ii\eta)\in \CC^2\colon
\\
\nonumber
& \qquad
    |y| < \kappa_0\alpha \;\mbox{ and }\;
    |\arctan (\eta / \xi)| < \kappa_0\alpha ,\ \xi > 0
    \bigr\} \,,
\\
\label{e:Sigma}
  \Sigma^{(\alpha)} (\nu_0)
& = \left\{ t = \alpha + \ii\tau\in \CC\colon
            \nu_0 |\tau| < \alpha \right\}
  = \alpha
  + \ii\left( - \nu_0^{-1}\alpha \,, \nu_0^{-1}\alpha \right)
\end{align}
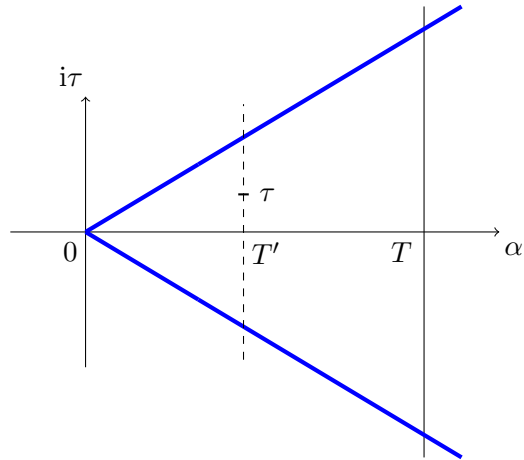
(Figure~3),
and for $0 < T'\leq T\leq \infty$,
we introduce the following complex parabolic domain,
\begin{equation}
\label{e:Gamma}
    \Gamma^{(T')}_T (\kappa_0,\nu_0)
  = \bigcup_{ \alpha\in (0,T) }
    \left[
    \mathfrak{V}^{\left( \kappa_0\cdot \min\{ \alpha, T'\} \right)}
    \times \Sigma^{(\alpha)} (\nu_0)
    \right] \subset \CC^2\times \CC
\end{equation}
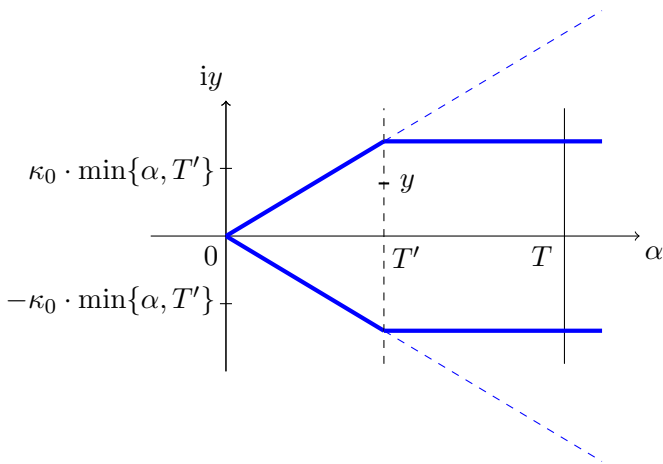
(Figure~4).
Additional properties of this domain will be presented later,
in Section~\ref{s:L^2-bound}, eq.~\eqref{e:Gamma_x}.

In order to get a better picture of the domain
$\Gamma^{(T')}_T (\kappa_0,\nu_0)\subset \CC^2\times \CC$,
it is worth to notice that the mapping
\begin{math}
  (z,\zeta,t) \;\longmapsto\; (z, \log\zeta, \log t)
\end{math}
maps $\Gamma^{(T')}_T (\kappa_0,\nu_0)$ diffeomorphically onto
the set of all complex triples
\begin{align*}
& (z,\zeta',t') = (x + \ii y, \xi' + \ii\eta', \alpha' + \ii\tau')
  \equiv
  (x, \xi',\alpha') + \ii (y,\eta',\tau')
  \in \CC^2\times \CC\simeq \RR^3\times \RR^3 \,,
\\
&   \quad\mbox{ such that }\quad
  0 < \alpha = \RE t = \ee^{\alpha'}\cdot \cos\tau' < T
    \;\mbox{ together with }\;
\\
&   |y| < \kappa_0\alpha \,,\
    |\eta'| < \arctan (\kappa_0\alpha) \,,\;\mbox{ and }\;
    |\tau'| < \arctan (1 / \nu_0) \,.
\end{align*}
In particular, there is no restriction on $x$ and $\xi'$ in the plane
$(x,\xi')\in \RR^2$, while $\alpha' = \log |t|\in \RR$.
These claims follow from simple calculations using
$\zeta = \ee^{\xi'}\cdot \ee^{\ii\eta'}$ and
$t = \ee^{\alpha'}\cdot \ee^{\ii\tau'}$.


\section{Main result}
\label{s:Main}

Our main result, Theorem~\ref{thm-Main}, gives the analyticity
(more precisely, a holomorphic extension to a complex domain)
of a unique weak solution to the homogeneous initial value problem
\eqref{e:Cauchy} with $f\equiv 0$ in $\HH\times (0,T)$.
Such a weak solution exists and is unique by
the following classical result
(Proposition~\ref{prop-Lions})
that summarizes a pair of standard theorems for
abstract parabolic problems due to
{\sc J.-L.\ Lions} \cite[Chapt.~IV]{Lions-61},
Th\'eor\`eme 1.1 ({\S}1, p.~46) and Th\'eor\`eme 2.1 ({\S}2, p.~52).
For alternative proofs, see also e.g.\
{\sc L.~C.\ Evans} \cite[Chapt.~7, {\S}1.2(c)]{Evans-98},
Theorems 3 and~4, pp.\ 356--358,
{\sc J.-L.\ Lions} \cite[Chapt.~III, {\S}1.2]{Lions-71},
Theorem 1.2 (p.~102) and remarks thereafter (p.~103),
{\sc A.\ Friedman} \cite{Friedman-64},
Chapt.~10, Theorem~17, p.~316, or
{\sc H.\ Tanabe} \cite[Chapt.~5, {\S}5.5]{Tanabe},
Theorem 5.5.1, p.~150.
%

\begin{proposition}\label{prop-Lions}
Let\/
$\rho$, $\sigma$, $\theta$, $q_r$, and\/ $\gamma$,
be given constants in $\RR$,
$\rho\in (-1,1)$, $\sigma > 0$, $\theta > 0$, and\/ $\gamma > 0$.
Assume that\/ $\kappa\in \RR$ is sufficiently large,
such that both inequalities,
\eqref{e:Feller} {\rm ({\it Feller's condition\/}) }
and \eqref{ineq:c_1'>0} are satisfied.
Next, let us choose $\beta\in \RR$ such that
$1 < \beta\leq 2\kappa\theta / \sigma^2$.
Set\/
$\mu = (\kappa / \sigma) - \gamma\, |\rho|$ $( > 0)$.
Let\/ $0 < T < \infty$, $f\in L^2((0,T)\to V')$, and\/ $u_0\in H$
be arbitrary.
Then the initial value problem \eqref{e:Cauchy} (with $u_0\in H$)
possesses a unique weak solution
\begin{equation*}
  u\in C([0,T]\to H)\cap L^2((0,T)\to V)
\end{equation*}
in the sense of\/ {\rm Definition~\ref{def-weak_sol}}.
Moreover, this solution satisfies also
$u\in W^{1,2}((0,T)\to V')$ and there exists a constant\/
$C\equiv C(T)\in (0,\infty)$, independent from $f$ and $u_0$, such that
\begin{equation}
\label{est:weak_sol}
\begin{aligned}
    \sup_{t\in [0,T]} \| u(t)\|_H^2
& {}
  + \int_0^T \| u(t)\|_V^2 \,\mathrm{d}t
  + \int_0^T 
    \left\| \genfrac{}{}{}1{\partial u}{\partial t}(t)
    \right\|_{V'}^2 \,\mathrm{d}t
\\
& \leq C
    \left( \| u_0\|_H^2
  + \int_0^T \| f(t)\|_{V'}^2 \,\mathrm{d}t \right) \,.
\end{aligned}
\end{equation}

Finally, if\/ $u_0\colon \HH\to \RR$ defined by
$u_0(x,\xi) = K\, (\mathrm{e}^x - 1)^{+}$, for\/ $(x,\xi)\in \HH$,
should belong to $H$, one needs to take $\gamma > 2$.
\end{proposition}
\par\vskip 10pt

The {\em proof\/} of this proposition is given towards the end of
Section~\ref{s:Heston-real}.
All that we have to do in this proof is to verify
the {\it boundedness\/} and {\it coercivity\/} hypotheses
for the sesquilinear form \eqref{e:Heston-diss:u=w} in $V\times V$
which are assumed in
{\sc J.-L.\ Lions} \cite[Chapt.~IV, {\S}1]{Lions-61},
inequalities (1.1) (p.~43) and (1.9) (p.~46), respectively.

Our main result is the following theorem which provides
an analytic extension of the weak solution $u$
to the initial value problem \eqref{e:Cauchy}
from the real domain $\HH\times [0,T]$ to a complex domain
$\Gamma^{(T')}_T (\kappa_0,\nu_0)$
defined in \eqref{e:Gamma}.

\begin{theorem}\label{thm-Main}
Let\/
$\rho$, $\sigma$, $\theta$, $q_r$, and\/ $\gamma$,
be given constants in $\RR$,
$\rho\in (-1,1)$, $\sigma > 0$, $\theta > 0$, and\/ $\gamma > 0$.
Assume that\/ $\beta$, $\gamma$, $\kappa$, and\/ $\mu$
are chosen as specified in {\rm Proposition~\ref{prop-Lions}} above.
Then the constants $\kappa_0, \nu_0\in (0,\infty)$ and\/ $T'\in (0,T]$
can be chosen sufficiently small and such that
the (unique) weak solution
\begin{equation*}
  u\in C([0,T]\to H)\cap L^2((0,T)\to V)
\end{equation*}
of the homogeneous initial value problem \eqref{e:Cauchy}
(with $f\equiv 0$ and\/ $u_0\in H$)
possesses a unique holomorphic extension
\begin{equation*}
  \tilde{u}\colon \Gamma^{(T')}_T (\kappa_0,\nu_0) \to \CC
\end{equation*}
to the complex domain
$\Gamma^{(T')}_T (\kappa_0,\nu_0) \subset \CC^3$
with the following properties:
There are some constants $C_0, c_0\in \RR_+$ such that\/
\begin{align}
\label{e:u_Hardy^2}
  \int_0^{\infty} \int_{-\infty}^{+\infty}
  \left| \tilde{u}
  \left( x + \ii y ,\, \xi (1 + \ii\omega) ,\, \alpha + \ii\tau
  \right)
  \right|^2\cdot \mathfrak{w}(x,\xi) \,\mathrm{d}x \,\mathrm{d}\xi
  \leq C_0\, \ee^{ c_0\alpha }\cdot \| u_0\|_H^2
\end{align}
for every\/ $\alpha\in (0,T]$ and for all\/
$y, \omega, \tau\in \RR$ satisfying\/
\begin{equation}
\label{dom:u_Hardy^2}
    \max\{ |y| ,\, |\arctan\omega| \}
  < \kappa_0\cdot \min\{ \alpha ,\, T'\}
    \quad\mbox{ and }\quad
  \nu_0 |\tau| < \alpha \,.
\end{equation}

Consequently, for any\/ $T_0\in (0,T']$, the domain
$\Gamma^{(T')}_T (\kappa_0,\nu_0)$ contains the Cartesian product
\begin{equation*}
  \mathfrak{X}^{(\kappa_0 T_0)} \times \Delta_{\kappa_0 T_0}
  \times\left[ (T_0,T)
   + \ii\left( {}- \genfrac{}{}{}1{T_0}{\nu_0} ,\,
                   \genfrac{}{}{}1{T_0}{\nu_0} \right)
        \right]
\end{equation*}
and the estimate in \eqref{e:u_Hardy^2} is valid for every\/
$\alpha\in [T_0,T]$ and for all\/
$y, \omega, \tau\in \RR$ such that, independently from $\alpha$,
\begin{equation}
\label{dom:T_0:u_Hardy^2}
    \max\{ |y| ,\, |\arctan\omega| \}
  < \kappa_0 T_0
    \quad\mbox{ and }\quad
    \nu_0 |\tau| < T_0 \,.
\end{equation}
\end{theorem}
\par\vskip 10pt

The {\em proof\/} of this theorem takes advantage of results from
Sections \ref{s:Heston-compl} and \ref{s:L^2-bound}, and
Appendix~\ref{s:density}.
It is formally completed at the end of Section~\ref{s:End}.


\section{An application to Mathematical Finance}
\label{s:Appl}

This section is concerned with an application of our main result,
Theorem~\ref{thm-Main} (Section~\ref{s:Main}),
to {\sc S.~L.\ Heston}'s {\em stochastic volatility\/} model
\cite{Heston} for {\em European call options\/}
described in Section~\ref{s:reform}.
Our goal will be to provide an affirmative answer to
the problem of {\em market completeness\/}
in Mathematical Finance as described in
{\sc M.~H.~A.\ Davis} and {\sc J.\ Ob{\l}{\'o}j} \cite{Davis-Obloj}.
We recall that the model is defined on a filtered probability space
$(\Omega, \mathcal{F}, (\mathcal{F}_t)_{t\geqslant 0}, \mathbb{P})$,
where $\mathbb{P}$ is the risk neutral probability measure.
Since an equivalent martingale measure exists,
but is \emph{not\/} unique, the market is {\em incomplete\/}.
The reader is referred to
{\sc M.~H.~A.\ Davis} \cite{Davis-royal},
{\sc J.~C.\ Hull} \cite{Hull-book},
{\sc J.\ Hull} and {\sc A.\ White} \cite{Hull-White},
{\sc A.~L.\ Lewis} \cite{Lewis},
{\sc E.~M.\ Stein} and {\sc J.~C.\ Stein} \cite{Stein-Stein}, and
{\sc J.~B.\ Wiggins} \cite{Wiggins}
for additional important work on this subject.
We closely follow the approach in \cite[Sect.~3]{Davis-Obloj} labeled
{\it ``martingale model''\/} for market completeness.
Another interesting paper on market completeness
deserves to be mentioned:
{\sc J.\ Hugonnier}, {\sc S.\ Malamud}, and {\sc E.\ Trubowitz}
\cite{HugoMalaTrubo}.
It is based on the existence of an Arrow\--Debreu equilibrium
and its implementation as a Radner equilibrium.
It is shown or assumed that in this setup,
allocation and prices are analytic functions of
the state and time variables.
The remaining arguments taking advantage of analytic entries
in the parabolic problem are similar to ours.

An extensive account of various stochastic volatility models for
European call options and possible market completion by such options
is given in
{\sc M.~H.~A.\ Davis} and {\sc J.\ Ob{\l}{\'o}j} \cite{Davis-Obloj},
{\sc M.\ Romano} and {\sc N.\ Touzi} \cite{RomanoTouzi}, and
{\sc P.\ Tak\'a\v{c}} \cite[Sect.~8, pp.\ 74--83]{Takac-12}.
Therefore, we restrict the discussion below to
the {\sc Heston\/} model \cite[Sect.~1]{Heston}
which seems to be very popular.
An important basic feature of this model is
the explicit form of its solution
\cite[pp.\ 330--331]{Heston}, eqs.\ (10) -- (18).
We apply our main analyticity result, Theorem~\ref{thm-Main},
to the {\sc Heston} model.
Another frequently used stochastic volatility model is the so\--called
{\it ``\,$3/{2}$ model''\/} investigated in
{\sc S.~L.\ Heston} \cite{Heston-1997},
{\sc P.\ Carr} and {\sc J.\ Sun} \cite{Carr-Sun},
{\sc A.\ Itkin} and {\sc P.\ Carr} \cite{Itkin-Carr},
and in the monographs by
{\sc J.\ Baldeaux} and {\sc E.\ Platen} \cite{Baldeaux-Platen}
and
{\sc A.~L.\ Lewis} \cite{Lewis}.
After a suitable transformation of variables,
it seems to be possible to treat the $3/2$ model by mathematical tools
similar to those we use in our present work.

We will answer the question of {\it\bfseries market completeness\/}
by investigating some qualitative properties
(such as analyticity) of the (unique) weak solution
\begin{equation*}
  u\in C([0,T]\to H)\cap L^2((0,T)\to V)
\end{equation*}
to the initial value problem \eqref{e:Cauchy}
obtained in our Theorem~\ref{thm-Main}.
Let us recall the Heston operator $\mathcal{A}$ defined in formula
\eqref{e:Heston-oper}.
The coefficients of the linear operator $\mathcal{A}$
are independent of time $t$ and $x\in \RR$,
and their dependence on $\xi\in (0,\infty)$ is very simple (linear).
As a natural consequence, the domain
$\Gamma^{(T')}_T (\kappa_0,\nu_0)$
of the holomorphic extension $\tilde{u}$ of the weak solution $u$
obtained in our Theorem~\ref{thm-Main} is simpler than
in the corresponding result obtained in
{\sc P.\ Tak\'a\v{c}} \cite[Theorem 3.3, pp.\ 58--59]{Takac-12}
for uniformly elliptic operators with variable analytic coefficients.

\begin{remark}\label{rem-rho,sigma}\nopagebreak
\begingroup\rm
It seems to be likely that one may allow both,
the correlation coefficient $\rho\equiv \rho(x,\xi,t)$ and
the volatility of volatility $\sigma\equiv \sigma(x,\xi,t)$
to depend on the variables $x$, $\xi$, and $t$, provided
this dependence is analytic, with all partial derivatives bounded,
and both functions $\rho$ and $\sigma$ bounded below and above by
some positive constants.

Last but not least, we would like to mention that
negative values of the correlation coefficient
$\rho\in (-1,1)$ are \emph{not\/} unusual in a volatile market:
asset prices tend to decrease when volatility increases
(\cite[p.~41]{FouqPapaSir}).
\null\hfill\Square
\endgroup
\end{remark}
\par\vskip 10pt

The market completion by a European call option has been obtained in
{\sc M.~H.~A.\ Davis} and {\sc J.\ Ob{\l}{\'o}j}
\cite[Proposition 5.1, p.~56]{Davis-Obloj}
based on the {\it validity\/} of a more general analyticity result
\cite[Theorem 4.1, p.~54]{Davis-Obloj}.
However, the main hypothesis in this theorem is
the {\it analyticity\/} of the solution
$\overline{p}(x,v,t) = p(x,v,T-t)$ of the parabolic problem
\eqref{e:div_Cauchy} in the domain $\HH\times (0,T)$.
(Warning: We use the symbol $\overline{p}$ to denote the function
$(x,v,t)\mapsto p(x,v,T-t)$, not the complex conjugate of $p$.)
Of course, the initial condition
$h(x) = K\, (\ee^x - 1)^{+}$, $x\in \RR$, is {\it not\/} analytic.
Nevertheless, in our Theorem~\ref{thm-Main}
we have established the analyticity result missing in
\cite{Davis-Obloj} (Theorem 4.1, p.~54).
Consequently, all conclusions in \cite{Davis-Obloj}
on market completion, that are based on the validity of
Theorem 4.1 (\cite[p.~54]{Davis-Obloj}),
are valid for the Heston model.
In Heston's model with a European call option, the notion of
a {\em complete market\/} is rigorously defined in
\cite[Definition 3.1, p.~52]{Davis-Obloj} as follows
(in probabilistic and measure\--theoretic terms):
\begingroup\it
Every contingent claim can be replicated by
a self\--financing trading strategy in the stock and bond\/
\endgroup
(contingent claims can be perfectly hedged against risks).
This is the case for Heston's model supplemented by
a European call option, by
Corollary 4.2 (p.~54) and Proposition 5.1 (p.~56)
in \cite{Davis-Obloj}.
We now briefly sketch how the analyticity of the solution
$u(x,\xi,t)$ in $\HH\times (0,T)$ facilitates market completion.
We keep the notation $u(x,\xi,t)$ for
a weak solution to problem \eqref{e:Cauchy}
which is the specific form of problem \eqref{e:div_Cauchy}
for Heston's model.
The relation between the solution
$\overline{p}(x,v,t) = p(x,v,T-t)$ of the parabolic problem
\eqref{e:div_Cauchy} and the weak solution $u(x,\xi,t)$
to the initial value problem \eqref{e:Cauchy} is obvious, i.e.,
\begin{math}
  \overline{p}(x,v,t) = u(x,\xi,t) = u(x, v / \sigma ,t) ,
\end{math}
by means of the substitutions
$v = \sigma\xi$ with the new independent variable $\xi\in \RR_+$ and
$\theta_{\sigma} = \theta / \sigma\in \RR$,
and by replacing the constants $\kappa$ and $\theta$, respectively, by
$\kappa^{*} = \kappa + \lambda > 0$ and
\begin{math}
  \theta^{*} = \frac{\kappa\theta}{\kappa + \lambda} > 0 .
\end{math}
Hence, we may set $r = \lambda = 0$ in eq.~\eqref{e:div_Cauchy}.
Conversely, let
\begin{math}
  p\colon \HH\times (0,T)\to \RR\colon (x,v,t)\mapsto p(x,v,t)
\end{math}
denote the unique solution of
the ({\it terminal\/} value) Cauchy problem \eqref{e:IVP.p}.
We set $u(x,\xi,t) = p(x, \sigma\xi, T-t)$ for all
$(x,\xi)\in \HH$ and $t\in (0,T)$, so that
$u\colon [0,T]\to H$ is the (unique) weak solution of
the initial value problem \eqref{e:Cauchy}
used in Section~\ref{s:Main}, Theorem~\ref{thm-Main}.
By the main result of this article, Theorem~\ref{thm-Main},
the function
$u\colon \HH\times (0,T)\to \RR$
can be (uniquely) extended to a holomorphic function in the domain
$\Gamma^{(T')}_T (\kappa_0,\nu_0) \subset \CC^2\times \CC$.
Consequently, the Jacobian matrix
\begin{equation*}
  G(x,\xi,t) = \left(
  \begin{array}{cc}
  1 \,, & 0 \\
  \frac{\partial u}{\partial x} (x,\xi,t) \,,
& \frac{\partial u}{\partial\xi}(x,\xi,t) \\
  \end{array}
  \right)
\end{equation*}
of the mapping
\begin{math}
  (x,\xi)\mapsto \left( x, u(x,\xi,t) \right)\colon
  \HH\subset \RR^2\to \RR^2
\end{math}
possesses determinant
$\det G(x,\xi,t)$ $= \frac{\partial u}{\partial\xi} (x,\xi,t)$
with a holomorphic extension to
$\Gamma^{(T')}_T (\kappa_0,\nu_0)$.
The determinant $\det G$ being (real) analytic in all of
$\HH\times (0,T)$, its set of zeros is
either Lebesgue negligible (i.e., of zero Lebesgue measure)
or else it is the whole domain $\HH\times (0,T)$
(cf.\
 {\sc S.~G.\ Krantz} and {\sc H.~R.\ Parks} \cite[p.~83]{Krantz-Parks}).
Hence, it suffices to examine $\det G$ in
an arbitrarily small neighborhood of a single ``central'' point.
An analogous result may be obtained in case when
analyticity can be obtained only in time $t$; see
\cite{Anders-Raimond, Davis-Obloj, HugoMalaTrubo, Kramkov-1, Kramkov-Pred}.
This case requires smoother terminal data, cf.\
Remark~\ref{rem-complete}, Part~{\rm (iii)}, below.

Finally, we can apply Proposition 5.1 (and its proof)
from \cite[p.~56]{Davis-Obloj}
to conclude that a European call option
in Heston's model \eqref{e:SVmodel}
{\em completes the market\/}:

\begin{theorem}\label{thm-complete}
Assume that\/ $\kappa > 0$ is sufficiently large, such that at least
the {\em Feller condition\/} \eqref{e:Feller} is satisfied; cf.\
{\rm Proposition~\ref{prop-Lions}}.
Assume that the payoff function
$h(x) = \hat{h}(K\ee^x)$ is {\em not\/} affine, that is,
$h''(x) = 0$ does not hold for every $x\in \RR$.
Then the stochastic volatility model~\eqref{e:SVmodel}
with a European call option yields a {\em complete market\/}.
\end{theorem}
\par\vskip 10pt

Under quite different sufficient conditions,
a related result on market completeness is established in
{\sc M.\ Romano} and {\sc N.\ Touzi}
\cite[Theorem 3.1, p.~406]{RomanoTouzi}:
\begingroup\it
A single European call option completes the market when there is
stochastic volatility driven by one extra Brownian motion\/
\endgroup
(under some additional assumptions; see
 \cite[pp.\ 404--407]{RomanoTouzi}).
The inequality
\begin{math}
  \det G(x,\xi,t) = \frac{\partial u}{\partial\xi} (x,\xi,t)
  \neq 0
\end{math}
(more precisely, $\frac{\partial u}{\partial\xi} (x,\xi,t) > 0$)
plays also there a decisive role.
Unlike in our present work, the inequality
$\frac{\partial u}{\partial\xi} (x,\xi,t) > 0$
in \cite[Theorem 3.1, p.~406]{RomanoTouzi}
is obtained directly from the convexity of the function
$h(x) =  K\, (\ee^x - 1)^{+}$ of $x\in \RR$
combined with the strong maximum principle for linear parabolic problems
which yields
$\frac{\partial^2 u}{\partial x^2} (x,\xi,t) > 0$
and, thus, the strict convexity of the function
$x\mapsto u(x,\xi,t)$ of $x\in \RR$ needed in
\cite[Theorem 3.1]{RomanoTouzi}.
Since we do not impose any convexity hypothesis on
the terminal function $h(x)$,
we are able to valuate much more general contingent claims than just
European call or put options.
An earlier result in
{\sc P.\ Tak\'a\v{c}} \cite[Theorem 8.5, p.~82]{Takac-12}
covers an alternative stochastic volatility model from
{\sc J.-P.\ Fouque}, {\sc G.\ Papanicolaou}, and {\sc K.~R.\ Sircar}
\cite[{\S}2.5, p.~47]{FouqPapaSir}, eqs.\ (2.18) -- (2.19).
The parabolic partial differential operator
(i.e., the It\^{o} operator)
in this model is uniformly parabolic and, consequently,
mathematically entirely different from
the degenerate It\^{o} operator in the Heston model.
Our main analyticity result,
Theorem~\ref{thm-Main} (Section~\ref{s:Main}),
is specialized to cover {\sc Heston}'s model and, consequently,
does not seem to be directly applicable to
the stochastic volatility models in
\cite{FouqPapaSir, Hull-White, Lewis, Stein-Stein, Wiggins}.

Based on the result in Theorem~\ref{thm-complete} above,
combined with those in
{\sc I.\ Bajeux\--Besnainou} and {\sc J.-Ch.\ Rochet}
\cite[p.~12]{Bajeux-Rochet},
we suggest the following (alternative)
{\it\bfseries analytic\/} definition of a {\em complete market\/},
at least in the case of Heston's model:

\begin{definition}\label{def-complete}\nopagebreak
\begingroup\rm
There is a set\/
$N\subset \HH\times (0,\infty)\subset \RR^2\times \RR$
of zero Lebesgue measure such that the mapping\/
\begin{math}
  \pi_t\colon (x,v)\mapsto \left( x, \overline{p}(x,v,t)\right)
  \colon \HH\subset \RR^2\to \RR^2
\end{math}
is a local diffeomorphism at every point\/
\begin{math}
  (x_0,v_0,t)\in \left[ \HH\times (0,\infty)\right] \setminus N .
\end{math}

Equivalently, for every $t\in (0,\infty)$, the set
\begin{math}
  N_t = \{ (x,v)\in \HH\colon (x,v,t)\in N\}\subset \RR^2
\end{math}
has zero Lebesgue measure and, at the point
$(x_0,v_0)\in \HH\setminus N_t$, the Jacobian matrix
\begin{equation*}
  J(x_0,v_0,t) = \left(
  \begin{array}{cc}
  1 \,, & 0 \\
  \frac{ \partial\overline{p} }{\partial x} (x,v,t) \,,
& \frac{ \partial\overline{p} }{\partial v} (x,v,t) \\
  \end{array}
  \right) \Bigg\vert_{ (x,v) = (x_0,v_0) }
\end{equation*}
of the mapping $\pi_t$ is regular which means that
\begin{math}
  \det J(x_0,v_0,t) =
  \frac{ \partial\overline{p} }{\partial v} (x,v,t)
    \big\vert_{ (x,v) = (x_0,v_0) }
  \neq 0 .
\end{math}
\endgroup
\end{definition}
\par\vskip 10pt

The property
\begin{math}
  \frac{ \partial\overline{p} }{\partial v} (x_0,v_0,t) \neq 0
\end{math}
allows us to apply the local {\em implicit function theorem\/}
to conclude that, by fixing $(x_0,t)$, we obtain an open neighborhood
$(v_0 - \delta, v_0 + \delta)$ of $v_0\in (0,\infty)$
($0 < \delta < \infty$ small enough)
such that either
\begin{math}
  \frac{ \partial\overline{p} }{\partial v} (x_0, \,\cdot\, ,t) > 0
\end{math}
(which is the case in \cite{Bajeux-Rochet, RomanoTouzi}),
or else
\begin{math}
  \frac{ \partial\overline{p} }{\partial v} (x_0, \,\cdot\, ,t) < 0
\end{math}
holds throughout $(v_0 - \delta, v_0 + \delta)$.
Hence, the function
\begin{math}
  \overline{p} (x_0, \,\cdot\, ,t) \colon
  (v_0 - \delta, v_0 + \delta)\hfil\break \to \RR
\end{math}
is either strictly monotone increasing or else
strictly monotone decreasing.
This means that, in a small (open) neighborhood of $v_0$,
one can perfectly hedge against small volatility fluctuations,
expressed through the variance $v = (\mathrm{volatility})^2$
satisfying $|v - v_0| < \delta$,
by a European call option $\overline{p} (x_0,v,t)$
priced near the value of $\overline{p} (x_0,v_0,t)$.
Merely the local {\em implicit function theorem\/} has to be envoked.

Our Definition~\ref{def-complete} is tailored for the completion of
the Heston model of a market with only a pair of random variables,
$(X_t,V_t)_{t\geqslant 0}$, as it appears also in
{\sc I.\ Bajeux\--Besnainou} and {\sc J.-Ch.\ Rochet}
\cite[p.~12]{Bajeux-Rochet}.
However, their market completion result in
\cite[Proposition 5.2, p.~12]{Bajeux-Rochet}
does not cover the Heston model.
A closely related definition of a complete market with
multiple random variables is given in 
{\sc M.~H.~A.\ Davis} and {\sc J.\ Ob{\l}{\'o}j}
\cite[Definition 3.1, p.~52]{Davis-Obloj}.
Their two main results in \cite{Davis-Obloj},
Theorem 3.2 (p.~52) which implies Theorem 4.1 (p.~54),
show that our Definition~\ref{def-complete} implies
that also the classical definition of a complete market from
{\sc J.~M.\ Harrison} and {\sc S.~R. Pliska}
\cite[{\S}3.4, pp.\ 241--242]{HarrisPliska-1}
and
\cite[p.~314]{HarrisPliska-2} is fulfilled
(see also {\sc I.\ Karatzas} and {\sc S.~E. Shreve}
 \cite[Chapt.~1, Def.\ 6.1, p.~21]{Karatzas-Shreve}).
For the market completion by a European call or put option,
another definition closely related to ours
(Definition~\ref{def-complete})
can be found in {\sc M.\ Romano} and {\sc N.\ Touzi}
\cite[Definition 3.1, p.~404]{RomanoTouzi}.

\begin{remark}\label{rem-complete}\nopagebreak
\begingroup\rm
{\rm (i)}$\;$
We stress that our {\rm Theorem~\ref{thm-Main}} (Section~\ref{s:Main})
allows to consider any payoff function $h\in H$,
$h(x,v)\equiv h(x) = \hat{h}(K\ee^x)$ for $x\in \RR$, in particular.
This is a significant advantage over the corresponding result in
{\sc P.\ Tak\'a\v{c}} \cite[Theorem 3.3, p.~59]{Takac-12}
which allows only for a payoff function $h\in L^2(\RR)$.
The hypothesis that the payoff function $h\colon \RR\to \RR$
is {\em not\/} affine is technical and comes from the proof of
Proposition 5.1 in \cite[Eq.\ (5.2), p.~57]{Davis-Obloj}.
It excludes a solution $u(x,\xi,t)$ with the partial derivative
\begin{math}
  \frac{\partial u}{\partial x} (x,\xi,t)
  \equiv \mathrm{const}(\xi,t)\in \RR
\end{math}
independent from $x\in \RR$.

{\rm (ii)}$\;$
The {\it Feller condition\/} \eqref{e:Feller}
(cf.\ \cite{Feller, Guo-Grzel-Ooster})
is needed to guarantee the unique solvability and well\--posedness
of the initial value problem \eqref{e:Cauchy}.
This condition was discovered in {\sc W.\ Feller} \cite{Feller}
for the corresponding parabolic problem in the variables
$(\xi,t)\in (0,\infty)^2$ only.
If this condition is violated, a suitable boundary condition on
the behavior of the solution $u(\xi,t)$ needs to be imposed as
$\xi\to 0+$.
{\sc Feller}'s result \cite{Feller} explains why we are able to prove
the well\--posedness of problem \eqref{e:Cauchy} with practically
{\em no\/} boundary conditions as $\xi\to 0+$ or $\xi\to \infty$,
except for \eqref{bc_xi:bound_u} and \eqref{bc:bound_w}
and the requirement that
$u(\,\cdot\,, \,\cdot\,, t)\in H$ together with
\eqref{bc_x:bound_u} and \eqref{e:trace:x=+-infty}
for every $t\in [0,T]$.
Notice that the last three conditions are easily satisfied by
a regular solution, thanks to $\beta > 1$ and $\gamma > 2$.
Our additional {\em condition on the size\/} of $\kappa > 0$,
i.e., $\kappa$ large enough, comes from the facts that
we have to deal with a solution $u(x,\xi,t)$ depending also
on the additional space variable $x\in \RR$ and
our underlying function space $H$ is the Hilbert space
$H = L^2(\mathbb{H};\mathfrak{w})$ with a special weight
$\mathfrak{w}(x,\xi)$.
the initial value $u(0) = u_0$ in $H$;

{\rm (iii)}$\;$
A number of recent articles concerned with endogenous completeness of
a market including stocks and options
(\cite{Anders-Raimond, Davis-Obloj, HugoMalaTrubo, Kramkov-1,
       Kramkov-Pred})
deal with solutions of
a Black\--Scholes\--It\^{o}\--type parabolic problem
that are analytic only in the time variable $t$.
As a result, these works need to impose more restrictive hypotheses
on the coefficients in the equation and the terminal data of
the parabolic problem, while no space analyticity is required for
the coefficients.
In contrast, the articles using a solution that is analytic
in both, the space and time variables $x$ and $t$
(\cite{Davis-Obloj, Takac-12}), need much less restrictive hypotheses
on the coefficients in the equation and the terminal data,
while space and time analyticity is required for the coefficients.
We refer to \cite[{\S}2 and {\S}5]{Davis-Obloj} and
\cite[Remark 3.3, p.~7]{Kramkov-1} for further details.
\endgroup
\end{remark}

\begin{remark}\label{rem-compl_other}\nopagebreak
\begingroup\rm
The {\it ``\,$3/2$ stochastic volatility model''\/}
\cite{Baldeaux-Platen, Carr-Sun, Heston-1997, Itkin-Carr, Lewis}
mentioned at the beginning of this section requires
some major changes in technical details used in our present work,
although we believe that
similar mathematical tools can still be applied.
For instance, the weight function
$\mathfrak{w}(x,\xi)$ defined in \eqref{def:w}
and the sesquilinear form
$(\mathcal{A}u, w)_H$ defined in \eqref{e:Heston-diss:u=w}
will have to be changed significantly.
\endgroup
\end{remark}
\par\vskip 10pt


\section{The Heston operator in the real domain}
\label{s:Heston-real}

At the end of this section {\bf we prove\/} Proposition~\ref{prop-Lions}
by verifying the boundedness and coercivity hypotheses
(in {\S}\ref{ss:bound-R} and {\S}\ref{ss:coerce-R}, respectively)
for the sesquilinear form \eqref{e:Heston-diss:u=w} in $V\times V$
assumed in
{\sc J.-L.\ Lions} \cite[Chapt.~IV, {\S}1]{Lions-61},
inequalities (1.1) (p.~43) and (1.9) (p.~46),
respectively.

Our boundedness and coercivity results for the Heston operator
$\mathcal{A}\colon V\to V'$
make use of five lemmas stated and proved in the Appendix
(Appendix~\ref{s:Trace,Sobolev}).
Recall that $\beta > 0$, $\gamma > 0$, and $\mu > 0$
are constants in the weight $\mathfrak{w}(x,\xi)$
which is defined in eq.~\eqref{def:w}.


\subsection{Boundedness of the Heston operator}
\label{ss:bound-R}

In this paragraph we verify the boundedness of
the sesquilinear form \eqref{e:Heston-diss:u=w} in $V\times V$.
This property is equivalent to
$\mathcal{A}$ being {\it bounded\/} as a linear operator
from $V$ to $V'$.

\begin{proposition}\label{prop-bound-R}
{\rm (Boundedness.)}$\;$
Let\/
$\beta$, $\gamma$, $\mu$, $\rho$, $\sigma$, $\theta$, $q_r$, and\/
$\kappa$ be given constants in $\RR$,
$\beta > 1$, $\gamma > 0$, $\mu > 0$, $-1 < \rho < 1$, $\sigma > 0$,
and\/ $\theta > 0$.
Then there exists a constant\/ $C\in (0,\infty)$, such that,
for all pairs $u,w\in V$, we have
\begin{equation}
\label{est:Heston-diss}
  \left| (\mathcal{A}u, w)_H\right|
  \leq C\cdot \| u\|_V\cdot \| w\|_V \,.
\end{equation}
\end{proposition}

\proof
For any given $u,w\in V$,
we apply Cauchy's inequality to the right\--hand side of
eq.~\eqref{e:Heston-diss:u=w} to estimate the inner product
\begin{equation*}
\begin{aligned}
& \left| (\mathcal{A}u, w)_H\right| \leq
\\
& \frac{\sigma}{2}\int_{\HH}
  \bigl[
    ( |u_x| + 2 |\rho|\, |u_{\xi}| )\cdot |\bar{w}_x|
  + |u_{\xi}|\cdot |\bar{w}_{\xi}|
  \bigr] \cdot \xi\cdot \mathfrak{w}(x,\xi)
    \,\mathrm{d}x \,\mathrm{d}\xi
\\
& {} + \frac{1}{2}\int_{\HH}
  \bigl[
    (1+\gamma)\sigma\, |u_x|
  + \left( | 2\kappa - \mu\sigma | + 2\gamma\rho\sigma \right)\, |u_{\xi}|
  \bigr] \cdot |\bar{w}|\cdot \xi\cdot \mathfrak{w}(x,\xi)
    \,\mathrm{d}x \,\mathrm{d}\xi
\\
& \begin{aligned}
  {} + \int_{\HH}
&   \left(
      |q_r|\, |u_x|
    + \left| \genfrac{}{}{}1{1}{2} \beta\sigma
           - \kappa\theta_{\sigma} \right|\, |u_{\xi}|
    \right) \cdot |\bar{w}|\cdot \mathfrak{w}(x,\xi)
    \,\mathrm{d}x \,\mathrm{d}\xi \,.
  \end{aligned}
\end{aligned}
\end{equation*}
(We abbreviate
 $\theta_{\sigma}\eqdef \theta / \sigma\in \RR$.)

With the abbreviations of the five integrals below,
\begin{align*}
  A_1 &= \int_{\HH}
    ( |u_x| + 2 |\rho|\, |u_{\xi}| )^2\cdot \xi\cdot \mathfrak{w}(x,\xi)
    \,\mathrm{d}x \,\mathrm{d}\xi \,,
\\
  B_1 &= \int_{\HH}
      |w_x|^2\cdot \xi\cdot \mathfrak{w}(x,\xi)
    \,\mathrm{d}x \,\mathrm{d}\xi \,,
\\
  A_2 &= \int_{\HH}
      |u_{\xi}|^2\cdot \xi\cdot \mathfrak{w}(x,\xi)
    \,\mathrm{d}x \,\mathrm{d}\xi \,,\quad
  B_2  = \int_{\HH}
      |w_{\xi}|^2\cdot \xi\cdot \mathfrak{w}(x,\xi)
    \,\mathrm{d}x \,\mathrm{d}\xi \,,
\\
  J   &= \int_{\HH}
    ( |u_x| + |u_{\xi}| )^2\cdot \xi\cdot \mathfrak{w}(x,\xi)
    \,\mathrm{d}x \,\mathrm{d}\xi
\\
      &\leq 2\int_{\HH}
    ( |u_x|^2 + |u_{\xi}|^2 )\cdot \xi\cdot \mathfrak{w}(x,\xi)
    \,\mathrm{d}x \,\mathrm{d}\xi \,,
\end{align*}
we thus obtain
\begin{align*}
& \left| (\mathcal{A}u, w)_H\right| \leq
  \frac{\sigma}{2} \left[ (A_1 B_1)^{1/2} + (A_2 B_2)^{1/2} \right]
\\
& {}
  + \frac{1}{2}\cdot
    \max\bigl\{
    (1+\gamma)\sigma ,\, | 2\kappa - \mu\sigma | + 2\gamma\rho\sigma
        \bigr\}\cdot J^{1/2}
  \left( \int_{\HH} |w|^2\cdot \xi\cdot \mathfrak{w}(x,\xi)
    \,\mathrm{d}x \,\mathrm{d}\xi
  \right)^{1/2}
\\
& {}
  + 
  {}\max\bigl\{ |q_r| ,\,
    \left| \genfrac{}{}{}1{1}{2} \beta\sigma - \kappa\theta_{\sigma}
    \right|
        \bigr\}\cdot J^{1/2}
  \left( \int_{\HH} \frac{|w(x,\xi)|^2}{\xi}\cdot \mathfrak{w}(x,\xi)
    \,\mathrm{d}x \,\mathrm{d}\xi
  \right)^{1/2} \,.
\end{align*}
With the help of these abbreviations and
the Cauchy\--type elementary inequality
\begin{equation*}
  (A_1 B_1)^{1/2} + (A_2 B_2)^{1/2}
  \leq (A_1 + A_2)^{1/2}\cdot (B_1 + B_2)^{1/2} \,,
\end{equation*}
which is equivalent with
\begin{math}
  \left[ (A_1 B_2)^{1/2} - (A_2 B_1)^{1/2} \right]^2 \geq 0 \,,
\end{math}
%
the last inequality above yields
\begin{align*}
& \left| (\mathcal{A}u, w)_H\right| \leq
    \frac{\sigma}{2}\, (A_1 + A_2)^{1/2}\cdot (B_1 + B_2)^{1/2}
\\
&
\begin{aligned}
  {}+ M_1
& \left( \int_{\HH}
    \left( |u_x|^2 + |u_{\xi}|^2 \right)
    \cdot \xi\cdot \mathfrak{w}(x,\xi)
    \,\mathrm{d}x \,\mathrm{d}\xi
  \right)^{1/2}
\\
  \times
& \left[ \int_{\HH}
    \left(
    \genfrac{|}{|}{}0{w(x,\xi)}{\xi}^2 + |w|^2
    \right)
    \cdot \xi\cdot \mathfrak{w}(x,\xi)
    \,\mathrm{d}x \,\mathrm{d}\xi
  \right]^{1/2} \,,
\end{aligned}
\end{align*}
with the constant
\begin{equation*}
  M_1\eqdef 2\cdot
  \max\left\{
    \genfrac{}{}{}1{1}{2} (1+\gamma)\sigma ,\,
    \left| \kappa - \genfrac{}{}{}1{1}{2} \mu\sigma \right|
  + \gamma\rho\sigma ,\,
    |q_r| ,\,
    \left| \genfrac{}{}{}1{1}{2} \beta\sigma - \kappa\theta_{\sigma}
    \right|
    \right\} > 0 \,.
\end{equation*}

With the help of the Cauchy inequality
\begin{equation*}
  4 |\rho|\, |u_x|\cdot |u_{\xi}|
  \leq 4 |u_x|^2 + |\rho|^2\, |u_{\xi}|^2 \,,
\end{equation*}
whence
\begin{align*}
    ( |u_x| + 2 |\rho|\, |u_{\xi}| )^2 + |u_{\xi}|^2
& = |u_x|^2 + 4 |\rho|\, |u_x|\cdot |u_{\xi}|
  + (1 + 4 |\rho|^2)\, |u_{\xi}|^2
\\
& \leq 5 |u_x|^2 + (1 + 5\rho^2) |u_{\xi}|^2
  \leq 6 \left( |u_x|^2 + |u_{\xi}|^2\right) \,,
\end{align*}
by $|\rho| < 1$, this inequality yields
\begin{equation*}
  A_1 + A_2\leq 6\int_{\HH}
    \left( |u_x|^2 + |u_{\xi}|^2 \right)
    \cdot \xi\cdot \mathfrak{w}(x,\xi)
    \,\mathrm{d}x \,\mathrm{d}\xi
\end{equation*}
and, consequently, also
\begin{align*}
  \left| (\mathcal{A}u, w)_H\right|
& {}\leq
  \left( \int_{\HH}
    \left( |u_x|^2 + |u_{\xi}|^2 \right)
    \cdot \xi\cdot \mathfrak{w}(x,\xi)
    \,\mathrm{d}x \,\mathrm{d}\xi
  \right)^{1/2}
\\
&
\begin{aligned}
  \times
& \left\{
    \frac{\sigma}{2}\, \sqrt{6}
  \left( \int_{\HH}
    \left( |w_x|^2 + |w_{\xi}|^2 \right)
    \cdot \xi\cdot \mathfrak{w}(x,\xi)
    \,\mathrm{d}x \,\mathrm{d}\xi
  \right)^{1/2}
  \right.
\\
& {}+ M_1
  \left.
  \left[ \int_{\HH}
    \left(
    \genfrac{|}{|}{}0{w(x,\xi)}{\xi}^2 + |w|^2
    \right)
    \cdot \xi\cdot \mathfrak{w}(x,\xi)
    \,\mathrm{d}x \,\mathrm{d}\xi
  \right]^{1/2}
  \right\} \,.
\end{aligned}
\end{align*}
Applying the Sobolev and Hardy inequalities
\eqref{ineq:Heston:Sobol} and \eqref{e:Hardy-beta}
to this estimate we deduce that there exists a constant
$C\in (0,\infty)$, such that the estimate in \eqref{est:Heston-diss}
holds for all pairs $u,w\in V$.
Here, we recall that, by Remark~\ref{rem-equiv_norm}, the norm
$\| w\|_V^{\sharp}$ defined in the Hilbert space $V$
by eq.~\eqref{def_eq:w_prod-H^1}
is equivalent with the original norm $\| w\|_V$
defined by eq.~\eqref{def:w_prod-H^1}.

Proposition~\ref{prop-bound-R} is proved.
\qed
\par\vskip 10pt


\subsection{Coercivity in the real domain}
\label{ss:coerce-R}

We wish to investigate the Heston operator $\mathcal{A}$
as a densely defined, closed linear operator in
the weighted Lebesgue space $H = L^2(\mathbb{H};\mathfrak{w})$.

We investigate the {\it coercivity\/} of
the linear operator $\mathcal{A}$ in
$V = H^1(\mathbb{H};\mathfrak{w})$.
In fact, we will show that the coercivity property holds for
$\mathcal{A} + \genfrac{}{}{}1{1}{2} c_2'\, I$
in place of $\mathcal{A}$, where
$c_2' > 0$ is a suitable constant (large enough) specified
at the end of this paragraph.
As a trivial consequence, the linear operator
${}- \left( \mathcal{A} + \genfrac{}{}{}1{1}{2} c_2'\, I\right)$
is {\it dissipative\/} in~$H$.
For establishing the coercivity, hypotheses
\eqref{e:Feller} and \eqref{ineq:c_1'>0} described in
Remark~\ref{rem-prop-Feller} are crucial.

We use the sesquilinear form from eq.~\eqref{e:Heston-diss:u=w}
to verify the coercivity of the linear operator $\mathcal{A}$
in the Hilbert space $V$:
\begin{align}
\nonumber
& 2\cdot \RE (\mathcal{A}u, u)_H
  = J_1 + J_2 + \dots + J_5\equiv
\\
\label{Re:Heston-diss:u=w}
&
\begin{aligned}
& \sigma\int_{\HH}
  \left[
    u_x\cdot \bar{u}_x
  + \rho\, ( u_{\xi}\cdot \bar{u}_x + u_x\cdot \bar{u}_{\xi} )
  + u_{\xi}\cdot \bar{u}_{\xi}
  \right] \cdot \xi\cdot \mathfrak{w}(x,\xi)
    \,\mathrm{d}x \,\mathrm{d}\xi
\\
& {}+ \frac{\sigma}{2}\int_{\HH}
  (1 - \gamma\, \Sgn x)\, ( u_x\cdot \bar{u} + \bar{u}_x\cdot u )
    \cdot \xi\cdot \mathfrak{w}(x,\xi)
    \,\mathrm{d}x \,\mathrm{d}\xi
\\
& {}+ \int_{\HH}
  \left( \kappa - \gamma\rho\sigma\, \Sgn x
       - \genfrac{}{}{}1{1}{2} \mu\sigma
  \right) ( u_{\xi}\cdot \bar{u} + \bar{u}_{\xi}\cdot u )
    \cdot \xi\cdot \mathfrak{w}(x,\xi)
    \,\mathrm{d}x \,\mathrm{d}\xi
\end{aligned}
\\
\nonumber
&
\begin{aligned}
& {}+ q_r
    \int_{\HH} ( u_x\cdot \bar{u} + \bar{u}_x\cdot u )
    \cdot \mathfrak{w}(x,\xi)
    \,\mathrm{d}x \,\mathrm{d}\xi
\\
& {}+ \left( \genfrac{}{}{}1{1}{2} \beta\sigma - \kappa\theta_{\sigma}
      \right)
    \int_{\HH} ( u_{\xi}\cdot \bar{u} + \bar{u}_{\xi}\cdot u )
    \cdot \mathfrak{w}(x,\xi)
    \,\mathrm{d}x \,\mathrm{d}\xi \,.
\end{aligned}
\end{align}
All integrals on the right\--hand side converge absolutely for any
$u\in V$, by the proof of Proposition~\ref{prop-bound-R} above.

\begin{proposition}\label{prop-coerc-R}
{\rm (Coercivity.)}$\;$
Let\/
$\rho$, $\sigma$, $\theta$, $q_r$, and\/ $\gamma$
be given constants in $\RR$,
$\rho\in (-1,1)$, $\sigma > 0$, $\theta > 0$, and\/ $\gamma > 0$.
Assume that\/ $\beta$, $\gamma$, $\kappa$, and\/ $\mu$
are chosen as specified in {\rm Proposition~\ref{prop-Lions}}.
Then there exists a constant\/ $c_2'\in (0,\infty)$
such that the following {\em G\r{a}rding inequality\/}
\begin{equation}
\label{ineq:Heston-diss:RE}
  2\cdot \RE (\mathcal{A}u, u)_H
  \geq \sigma\, (1 - |\rho|)\cdot \| u\|_V^2 - c_2'\cdot \| u\|_H^2
\end{equation}
is valid for all\/ $u\in V$.
\end{proposition}

\proof
Let us consider eq.~\eqref{Re:Heston-diss:u=w}
with an arbitrary $u\in V$.
The first integral on the right\--hand side of
eq.~\eqref{Re:Heston-diss:u=w}
is estimated from below by Cauchy's inequality
\begin{equation*}
    u_{\xi}\cdot \bar{u}_x + u_x\cdot \bar{u}_{\xi}
  = 2\cdot\RE ( u_{\xi}\cdot \bar{u}_x )
  \leq 2 |u_{\xi}|\cdot |\bar{u}_x| \leq |u_x|^2 + |u_{\xi}|^2 ,
\end{equation*}
\begin{align}
\nonumber
& \frac{J_1}{\sigma} \equiv \int_{\HH}
  \left[
    u_x\cdot \bar{u}_x
  + \rho\, ( u_{\xi}\cdot \bar{u}_x + u_x\cdot \bar{u}_{\xi} )
  + u_{\xi}\cdot \bar{u}_{\xi}
  \right] \cdot \xi\cdot \mathfrak{w}(x,\xi)
    \,\mathrm{d}x \,\mathrm{d}\xi
\\
\label{e:A-coerce}
& \geq 
  \int_{\HH}
  \left[ |u_x|^2 - |\rho|\, ( |u_x|^2 + |u_{\xi}|^2 ) + |u_{\xi}|^2
  \right] \cdot \xi\cdot \mathfrak{w}(x,\xi)
    \,\mathrm{d}x \,\mathrm{d}\xi
\\
\nonumber
& = (1 - |\rho|)
  \int_{\HH} ( |u_x|^2 + |u_{\xi}|^2 )
    \cdot \xi\cdot \mathfrak{w}(x,\xi)
    \,\mathrm{d}x \,\mathrm{d}\xi
\\
\nonumber
& = (1 - |\rho|) \left( \| u\|_V^2 - \| u\|_H^2 \right) \,.
\end{align}

The second integral in eq.~\eqref{Re:Heston-diss:u=w},
$J_2$,
consists of two different parts that we treat by
integration\--by\--parts as follows,
using the following simple formulas,
\begin{align*}
& \frac{\partial}{\partial x}\, \mathfrak{w}(x,\xi)
  = {}- \gamma\, \xi^{\beta - 1}\, \ee^{ - \gamma |x| - \mu\xi }
    \cdot \Sgn x
  = {}- \gamma\cdot \Sgn x\cdot \mathfrak{w}(x,\xi) \,,
\\
&
\begin{aligned}
& \frac{\partial}{\partial\xi}\, \mathfrak{w}(x,\xi)
  = (\beta - 1)\, \xi^{\beta - 2}\, \ee^{ - \gamma |x| - \mu\xi }
  {}- \mu\, \xi^{\beta - 1}\, \ee^{ - \gamma |x| - \mu\xi }
\\
& = ( \beta - 1 - \mu\xi )\, \xi^{\beta - 2}\,
    \ee^{ - \gamma |x| - \mu\xi }
  = \left( \frac{\beta - 1}{\xi} - \mu \right)
    \cdot \mathfrak{w}(x,\xi) \,,
\end{aligned}
\\
&
\begin{aligned}
& \frac{\partial}{\partial\xi}
  \left( \xi\cdot \mathfrak{w}(x,\xi) \right)
  = \frac{\partial}{\partial\xi}
    \left( \xi^{\beta}\, \ee^{ - \gamma |x| - \mu\xi } \right)
\\
& = \beta\cdot \xi^{\beta - 1}\, \ee^{ - \gamma |x| - \mu\xi }
  - \mu\, \xi^{\beta}\, \ee^{ - \gamma |x| - \mu\xi }
  = (\beta - \mu\xi)\cdot \mathfrak{w}(x,\xi) \,.
\end{aligned}
\end{align*}
Consequently, the first part of the integral in
${2 J_2}/{\sigma}$ in eq.~\eqref{Re:Heston-diss:u=w}, becomes
\begin{equation*}
\begin{aligned}
& \int_{\RR} ( u_x\, \bar{u} + \bar{u}_x\, u )
    \cdot \ee^{- \gamma |x|} \,\mathrm{d}x
  = \int_{\RR} (|u|^2)_x\cdot \ee^{- \gamma |x|} \,\mathrm{d}x
\\
& = |u(x,\xi)|^2\cdot \ee^{- \gamma |x|}
    \Big\vert_{x = -\infty}^{x = +\infty}
  + \gamma\int_{\RR}
    |u(x,\xi)|^2\cdot \Sgn x\cdot \ee^{- \gamma |x|} \,\mathrm{d}x
\\
& = \gamma\int_{\RR}
    |u(x,\xi)|^2\cdot \Sgn x\cdot \ee^{- \gamma |x|} \,\mathrm{d}x
\end{aligned}
\end{equation*}
for almost every $\xi\in (0,\infty)$,
with a help from Lemma~\ref{lem-Trace_x}.
Integrating this equality with respect to $\xi\in (0,\infty)$
and the measure
$\xi^{\beta}\, \ee^{- \mu\xi} \,\mathrm{d}\xi$, we arrive at
\begin{equation}
\label{int_x:u_x.vw}
\begin{aligned}
&   \int_{\HH} ( u_x\, \bar{u} + \bar{u}_x\, u )
    \cdot \xi\cdot \mathfrak{w}(x,\xi) \,\mathrm{d}x \,\mathrm{d}\xi
\\
& = \gamma\int_{\HH}
    |u(x,\xi)|^2\cdot \Sgn x\cdot \xi\cdot \mathfrak{w}(x,\xi)
    \,\mathrm{d}x \,\mathrm{d}\xi \,.
\end{aligned}
\end{equation}
Recall that
\begin{math}
  \mathfrak{w}(x,\xi) =
  \xi^{\beta - 1}\, \ee^{ - \gamma |x| - \mu\xi } .
\end{math}
Similarly, we get
\begin{align*}
&
\begin{aligned}
& \int_{\RR}
    ( u_x\, \bar{u} + \bar{u}_x\, u )
    \cdot \Sgn x\cdot \ee^{- \gamma |x|}
    \,\mathrm{d}x
\\
& =
{}- \int_{-\infty}^0
    ( u_x\, \bar{u} + u\, \bar{u}_x )\, \ee^{\gamma x}
    \,\mathrm{d}x
  + \int_0^{\infty}
    ( u_x\, \bar{u} + u\, \bar{u}_x )\, \ee^{- \gamma x}
    \,\mathrm{d}x
\\
& =
{}- \int_{-\infty}^0
    (|u|^2)_x\cdot \ee^{\gamma x} \,\mathrm{d}x
  + \int_0^{\infty}
    (|u|^2)_x\cdot \ee^{- \gamma x} \,\mathrm{d}x
\end{aligned}
\\
\nonumber
&
\begin{aligned}
  = {} &
{}- |u(x,\xi)|^2\, \ee^{\gamma x}\,
    \bigg\vert_{-\infty}^0
  + \gamma\int_{-\infty}^0
    |u(x,\xi)|^2\, \ee^{\gamma x} \,\mathrm{d}x
\\
& + |u(x,\xi)|^2\, \ee^{- \gamma x}\,
    \bigg\vert_0^{\infty}
  + \gamma\int_0^{\infty}
    |u(x,\xi)|^2\, \ee^{- \gamma x} \,\mathrm{d}x
\end{aligned}
\\
\nonumber
&
\begin{aligned}
  = {} &
{}- 2 |u(0,\xi)|^2
  + \gamma\int_{-\infty}^{\infty}
    |u(x,\xi)|^2\, \ee^{- \gamma |x|} \,\mathrm{d}x \,.
\end{aligned}
\end{align*}
Integrating this equality with respect to $\xi\in (0,\infty)$
and the measure
$\xi^{\beta}\, \ee^{- \mu\xi} \,\mathrm{d}\xi$, we arrive at
\begin{equation}
\label{int_x:sign_x.vw}
\begin{aligned}
& \int_{\HH}
    ( u_x\, \bar{u} + u\, \bar{u}_x )\cdot \Sgn x\cdot \xi\cdot
    \mathfrak{w}(x,\xi)
    \,\mathrm{d}x \,\mathrm{d}\xi
\\
& =
{}- 2\int_0^{\infty} |u(0,\xi)|^2\, \xi^{\beta}\, \ee^{- \mu\xi}
    \,\mathrm{d}\xi
  + \gamma\int_{\HH}
    |u(x,\xi)|^2\cdot \xi\cdot \mathfrak{w}(x,\xi)
    \,\mathrm{d}x \,\mathrm{d}\xi \,.
\end{aligned}
\end{equation}
Finally, we combine the identities in
\eqref{int_x:u_x.vw} and \eqref{int_x:sign_x.vw} to obtain
\begin{align}
\nonumber
& \frac{2 J_2}{\sigma} \equiv \int_{\HH}
  (1 - \gamma\, \Sgn x)\, ( u_x\cdot \bar{u} + \bar{u}_x\cdot u )
    \cdot \xi\cdot \mathfrak{w}(x,\xi)
    \,\mathrm{d}x \,\mathrm{d}\xi
\\
\label{int_v:u_x.u}
& {}
  = 2\gamma \int_0^{\infty}
    |u(0,\xi)|^2\, \xi^{\beta}\, \ee^{- \mu\xi} \,\mathrm{d}\xi
  - \gamma^2 \int_{\HH}
    |u(x,\xi)|^2\cdot \xi\cdot \mathfrak{w}(x,\xi)
    \,\mathrm{d}x \,\mathrm{d}\xi
\\
\nonumber
& {}
  + \gamma\int_{\HH}
    |u(x,\xi)|^2\cdot \Sgn x\cdot \xi\cdot \mathfrak{w}(x,\xi)
    \,\mathrm{d}x \,\mathrm{d}\xi \,.
\end{align}

In order to treat the third integral in eq.~\eqref{Re:Heston-diss:u=w},
we need to calculate
\begin{align*}
& \int_0^{\infty}
  ( u_{\xi}\cdot \bar{u} + \bar{u}_{\xi}\cdot u )
    \cdot \xi^{\beta}\, \ee^{- \mu\xi} \,\mathrm{d}\xi
  = \int_0^{\infty}
  (|u|^2)_{\xi}\cdot \xi^{\beta}\, \ee^{- \mu\xi} \,\mathrm{d}\xi
\\
& {}
  = |u(x,\xi)|^2\cdot \xi^{\beta}\, \ee^{- \mu\xi}
    \Big\vert_{\xi = 0}^{\xi = \infty}
  - \int_0^{\infty} |u(x,\xi)|^2\cdot
    (\beta - \mu\xi)\, \xi^{\beta - 1}\, \ee^{- \mu\xi}
    \,\mathrm{d}\xi \,.
\end{align*}
Integrating first this equality with respect to $x\in (-\infty, \infty)$
and the measure
$\ee^{- \gamma |x|} \,\mathrm{d}x$,
then applying the vanishing trace results
\eqref{e:trace:v=0} and \eqref{e:trace:v=infty}, we arrive at
\begin{align}
\nonumber
  J_3
& {}
  \equiv \int_{\HH}
  \left( \kappa - \gamma\rho\sigma\, \Sgn x
       - \genfrac{}{}{}1{1}{2} \mu\sigma
  \right) ( u_{\xi}\cdot \bar{u} + \bar{u}_{\xi}\cdot u )
    \cdot \xi\cdot \mathfrak{w}(x,\xi)
    \,\mathrm{d}x \,\mathrm{d}\xi
\\
\label{int_v:u_v.u}
& {}
  = {}- \left( \kappa - \genfrac{}{}{}1{1}{2} \mu\sigma \right)
    \int_{\HH} |u(x,\xi)|^2\cdot
    (\beta - \mu\xi)\, \mathfrak{w}(x,\xi)
    \,\mathrm{d}x \,\mathrm{d}\xi
\\
\nonumber
& {}
  + \gamma\rho\sigma \int_{\HH}
    |u(x,\xi)|^2\cdot \Sgn x\cdot (\beta - \mu\xi)\,
    \mathfrak{w}(x,\xi) \,\mathrm{d}x \,\mathrm{d}\xi \,.
\end{align}

The fourth integral in eq.~\eqref{Re:Heston-diss:u=w}
is treated analogously to the second one.
It suffices to replace $\beta$ by $\beta - 1$
in the equality \eqref{int_x:u_x.vw} which then yields
\begin{equation}
\label{int_v:u_v.u.vw}
\begin{aligned}
  \frac{J_4}{q_r}
& {}
  \equiv \int_{\HH}
  ( u_x\, \bar{u} + \bar{u}_x\, u )
    \cdot \mathfrak{w}(x,\xi) \,\mathrm{d}x \,\mathrm{d}\xi
\\
& = \gamma\int_{\HH}
    |u(x,\xi)|^2\cdot \Sgn x\cdot \mathfrak{w}(x,\xi)
    \,\mathrm{d}x \,\mathrm{d}\xi \,.
\end{aligned}
\end{equation}

Finally, the last integral in eq.~\eqref{Re:Heston-diss:u=w}
is treated analogously to the third one,
\begin{equation}
\label{int_v:u_v.u.w}
\begin{aligned}
  \frac{J_5}{ \genfrac{}{}{}1{1}{2} \beta\sigma - \kappa\theta_{\sigma} }
& {}
  \equiv \int_{\HH}
  ( u_{\xi}\cdot \bar{u} + \bar{u}_{\xi}\cdot u )
    \cdot \mathfrak{w}(x,\xi) \,\mathrm{d}x \,\mathrm{d}\xi
\\
& =
{}- \int_{\HH} |u(x,\xi)|^2\cdot
    \left( \frac{\beta - 1}{\xi} - \mu \right)
    \cdot \mathfrak{w}(x,\xi) \,\mathrm{d}x \,\mathrm{d}\xi \,.
\end{aligned}
\end{equation}

We collect the second through fifth integrals,
cf.\ eq.~\eqref{Re:Heston-diss:u=w},
\begin{align*}
& J_2 + \dots J_5 =
    \gamma\sigma \int_0^{\infty}
    |u(0,\xi)|^2\, \xi^{\beta}\, \ee^{- \mu\xi} \,\mathrm{d}\xi
\\
& {}
  + \left[
  {}- \genfrac{}{}{}1{1}{2} \sigma\gamma^2
  + \mu \left( \kappa - \genfrac{}{}{}1{1}{2} \mu\sigma \right)
    \right]
  \int_{\HH} |u(x,\xi)|^2\cdot \xi\cdot \mathfrak{w}(x,\xi)
    \,\mathrm{d}x \,\mathrm{d}\xi
\\
& {}
  + \left[
    \genfrac{}{}{}1{1}{2} \sigma\gamma
  - \mu\gamma\rho\sigma
    \right]
  \int_{\HH} |u(x,\xi)|^2\cdot \Sgn x\cdot \xi\cdot \mathfrak{w}(x,\xi)
    \,\mathrm{d}x \,\mathrm{d}\xi
\\
& {}
  + \left[
  {}- \beta \left( \kappa - \genfrac{}{}{}1{1}{2} \mu\sigma \right)
  + \mu
    \left( \genfrac{}{}{}1{1}{2} \beta\sigma - \kappa\theta_{\sigma}
    \right)
    \right]
  \int_{\HH} |u(x,\xi)|^2\cdot \mathfrak{w}(x,\xi)
    \,\mathrm{d}x \,\mathrm{d}\xi
\\
& {}
  + \left[
    \beta\gamma\rho\sigma
    + \gamma q_r
    \right]
  \int_{\HH} |u(x,\xi)|^2\cdot \Sgn x\cdot \mathfrak{w}(x,\xi)
    \,\mathrm{d}x \,\mathrm{d}\xi
\\
& {}- (\beta - 1)
    \left( \genfrac{}{}{}1{1}{2} \beta\sigma - \kappa\theta_{\sigma}
    \right)
    \int_{\HH} \frac{|u(x,\xi)|^2}{\xi}
    \cdot \mathfrak{w}(x,\xi) \,\mathrm{d}x \,\mathrm{d}\xi \,,
\end{align*}
whence
\begin{align}
\nonumber
& J_2 + \dots J_5\geq
\\
\label{e:J_2+...+J_5}
& {}
    \left\{
    \left[
    \mu\kappa - \genfrac{}{}{}1{1}{2} \sigma (\gamma^2 + \mu^2)
    \right]
  - \sigma\gamma
    \left| \genfrac{}{}{}1{1}{2} - \mu\rho \right|
    \right\}
  \int_{\HH} |u(x,\xi)|^2\cdot \xi\cdot \mathfrak{w}(x,\xi)
    \,\mathrm{d}x \,\mathrm{d}\xi
\\
\nonumber
& {}
  + \left\{
    \left[
    \beta\mu\sigma - \kappa (\beta + \mu\theta_{\sigma})
    \right]
  - \gamma
    \left| \beta\rho\sigma + q_r \right|
    \right\} \| u\|_H^2
\\
\nonumber
& {}+ (\beta - 1)
    \left( \kappa\theta_{\sigma} - \genfrac{}{}{}1{1}{2} \beta\sigma
    \right)
    \int_{\HH} \frac{|u(x,\xi)|^2}{\xi}
    \cdot \mathfrak{w}(x,\xi) \,\mathrm{d}x \,\mathrm{d}\xi
\\
\nonumber
& \equiv c_1
  \int_{\HH} |u(x,\xi)|^2\cdot \xi\cdot \mathfrak{w}(x,\xi)
    \,\mathrm{d}x \,\mathrm{d}\xi
  + c_2\cdot \| u\|_H^2
\\
\nonumber
& {}+ c_3
    \int_{\HH} \frac{|u(x,\xi)|^2}{\xi}
    \cdot \mathfrak{w}(x,\xi) \,\mathrm{d}x \,\mathrm{d}\xi \,,
\end{align}
where the constants
\begin{align*}
& c_1\eqdef
    \left[
    \mu\kappa - \genfrac{}{}{}1{1}{2} \sigma (\gamma^2 + \mu^2)
    \right]
  - \sigma\gamma
    \left| \genfrac{}{}{}1{1}{2} - \mu\rho \right| \,,
\\
& c_2\eqdef
    \left[
    \beta\mu\sigma - \kappa (\beta + \mu\theta_{\sigma})
    \right]
  - \gamma
    \left| \beta\rho\sigma + q_r\right| \,,
\\
& c_3\eqdef (\beta - 1)
    \left( \kappa\theta_{\sigma} - \genfrac{}{}{}1{1}{2} \beta\sigma
    \right) \,,
\end{align*}
are estimated from below as follows:
\begin{align}
\label{est:v^1}
& c_1\geq c_1'\eqdef
    \mu\kappa - \genfrac{}{}{}1{1}{2} \sigma (\gamma^2 + \mu^2)
  - \sigma\gamma
    \left( \genfrac{}{}{}1{1}{2} + \mu\, |\rho| \right) \,,
\\
\label{est:v^0}
& c_2 > -\infty \,,
\\
\label{est:u/v}
& c_3 = \frac{\beta - 1}{\sigma}
    \left( \kappa\theta - \genfrac{}{}{}1{1}{2} \beta\sigma^2
    \right)\geq 0 \,.
\end{align}
The constant $c_3\in \RR$ is nonnegative thanks to
{\it Feller's condition\/},
$\genfrac{}{}{}1{1}{2} \sigma^2 - \kappa\theta < 0$,
provided we choose $\beta\in \RR$ such that
\begin{math}
  1 < \beta\leq 2\kappa\theta / \sigma^2 .
\end{math}
The sign of the constant $c_2$ does not matter as it stands as
a coefficient with the norm $\| u\|_H$.
Finally, in order to guarantee $c_1'\geq 0$, we first choose
$\mu > 0$ such that this value of $\mu$ maximizes the function
\begin{align*}
  \mu\mapsto c_1'
& \equiv c_1'(\mu)
  = \mu\kappa - \genfrac{}{}{}1{1}{2} \sigma (\gamma^2 + \mu^2)
  - \sigma\gamma
    \left( \genfrac{}{}{}1{1}{2} + \mu\, |\rho| \right)
\\
& = \genfrac{}{}{}1{1}{2} \sigma
    \left[
  {}-
    \left( \mu - \frac{\kappa}{\sigma} + \gamma\, |\rho| \right)^2
  + \left( \frac{\kappa}{\sigma} - \gamma\, |\rho| \right)^2
  - \gamma (1+\gamma)
    \right] \,,
\end{align*}
that is,
$\mu = (\kappa / \sigma) - \gamma\, |\rho|$, provided
$\kappa > \sigma\gamma |\rho|$.
With this value of $\mu$, we have to satisfy
\begin{equation*}
  c_1' = \genfrac{}{}{}1{1}{2} \sigma
    \left[
    \left( \frac{\kappa}{\sigma} - \gamma\, |\rho| \right)^2
  - \gamma (1+\gamma)
    \right] \geq 0 \,,
\end{equation*}
that is, ineq.~\eqref{ineq:c_1'>0}.

Finally, applying
inequalities \eqref{est:v^1}, \eqref{est:v^0}, and \eqref{est:u/v}
to the right\--hand side of eq.~\eqref{e:J_2+...+J_5},
and inequality \eqref{e:A-coerce}
to eq.~\eqref{Re:Heston-diss:u=w}, we obtain
\begin{align}
\nonumber
& 2\cdot \RE (\mathcal{A}u, u)_H
  \geq \sigma\, (1 - |\rho|) \left( \| u\|_V^2 - \| u\|_H^2 \right)
\\
\label{e:Heston-diss:RE}
& {}+ c_1'
    \int_{\HH} |u(x,\xi)|^2\cdot \xi\cdot \mathfrak{w}(x,\xi)
    \,\mathrm{d}x \,\mathrm{d}\xi
  {}+ c_2\, \| u\|_H^2
\\
\nonumber
& {}+ c_3
    \int_{\HH} \frac{|u(x,\xi)|^2}{\xi}
    \cdot \mathfrak{w}(x,\xi) \,\mathrm{d}x \,\mathrm{d}\xi
  \geq \sigma\, (1 - |\rho|)\, \| u\|_V^2
     - c_2'\, \| u\|_H^2 \,,
\end{align}
where
$c_2'= \sigma\, (1 - |\rho|) + |c_2| > 0$ is a constant.

Consequently, the linear operator
$\mathcal{A} + \genfrac{}{}{}1{1}{2} c_2'\, I$
is {\em coercive\/} in $V$ and
$- \left( \mathcal{A} + \genfrac{}{}{}1{1}{2} c_2'\, I\right)$
is {\em dissipative\/} in~$H$.
More precisely, ineq.~\eqref{e:Heston-diss:RE},
when combined with our definitions of equivalent norms in
$V = H^1(\mathbb{H};\mathfrak{w})$, yields
the {\em G\r{a}rding inequality\/} in \eqref{ineq:Heston-diss:RE}.

The proof of Proposition~\ref{prop-coerc-R} is complete.
\qed
\par\vskip 10pt

\begin{remark}\label{rem-prop-coerc-R}\nopagebreak
{\rm (Feller's condition.)}$\;$
\begingroup\rm
{\it Feller's condition\/}
$\genfrac{}{}{}1{1}{2} \sigma^2 - \kappa\theta < 0$
and our choice of $\beta\in \RR$ such that
\begin{math}
  1 < \beta\leq 2\kappa\theta / \sigma^2
\end{math}
guarantee $c_3\geq 0$
in the proof of Proposition~\ref{prop-coerc-R} above.
In addition, to guarantee also
\begin{equation*}
  c_1' = \genfrac{}{}{}1{1}{2} \sigma
    \left[
    \left( \frac{\kappa}{\sigma} - \gamma\, |\rho| \right)^2
  - \gamma (1+\gamma)
    \right] \geq 0 \,,
\end{equation*}
we need to assume ineq.~\eqref{ineq:c_1'>0}.
\hfill\Square
\endgroup
\end{remark}
\par\vskip 10pt

{\it Proof of\/} {\bf Proposition~\ref{prop-Lions}.}$\;$
In Propositions \ref{prop-bound-R} and \ref{prop-coerc-R}
above we have verified the boundedness and coercivity hypotheses
for the linear operator $\mathcal{A}\colon V\to V'$
required in
{\sc J.-L.\ Lions} \cite[Chapt.~IV]{Lions-61},
Th\'eor\`eme 1.1 ({\S}1, p.~46) and Th\'eor\`eme 2.1 ({\S}2, p.~52).
Consequently, these well\--known results from
\cite[Chapt.~IV]{Lions-61}
yield the desired conclusion of Proposition~\ref{prop-Lions}
on the existence and uniqueness of a weak solution to
the initial value problem \eqref{e:Cauchy}.
Finally, the energy estimate \eqref{est:weak_sol}
can be found in {\sc L.~C.\ Evans} \cite[Chapt.~7, {\S}1.2(b)]{Evans-98},
Theorem~2, p.~354.
\qed
\par\vskip 10pt


\section{The Heston operator in the complex domain}
\label{s:Heston-compl}

In the first paragraph of this section, {\S}\ref{ss:Anal-compl_t},
we apply the classical theory of {\em sectorial operators\/} as
{\em infinitesimal generators\/} of {\em holomorphic semigroups\/}
of bounded linear operators in the complex Hilbert space
$H = L^2(\mathbb{H};\mathfrak{w})$.
This theory provides a (unique) holomorphic extension of
the unique weak solution
$u\colon \HH\times [0,T]\to \RR$ of the initial value problem
\eqref{e:Cauchy} with $f\equiv 0$,
obtained in Proposition~\ref{prop-Lions},
to the complex domain
$\HH\times \Delta_{\vartheta'}$ that is holomorphic in the time variable
$t\in \Delta_{\vartheta'}$.
To obtain a holomorphic extension of $u$ to the complex domain
\begin{math}
    \mathfrak{V}^{(r)}
  = \mathfrak{X}^{(r)}\times \Delta_{\arctan r} \subset \CC^2
\end{math}
in the space variables $(x,\xi)$,
that has been defined in eq.~\eqref{x,v:strip} for $r\in (0,\infty)$,
we first replace the (possibly nonsmooth) initial data $u_0\in H$
by an entire function
$u_{0,n}\colon \CC^2\to \CC$; $n=1,2,3,\dots$,
constructed in {\S}\ref{ss:Cauchy-compl}, such that
$u_{0,n}\vert_{\HH}\in H$, ineq.~\eqref{est:H_n-L_n} is valid,
and the sequence
\begin{math}
  \| u_{0,n}\vert_{\HH} - u_0\|_H \,\rightarrow\, 0
\end{math}
as $n\to \infty$.
Given such initial data $u_0\vert_{\HH}\in H$, where
$u_0\colon \CC^2\to \CC$ is an entire function satisfying
ineq.~\eqref{est:H_n-L_n},
the main result of the entire section,
Proposition~\ref{prop-Main} proved in {\S}\ref{ss:Cauchy-compl},
provides a (unique) holomorphic extension of the solution $u$
to the complex domain
\begin{math}
    \mathfrak{X}^{(r)}\times \Delta_{\arctan r}
    \times \Delta_{\vartheta'}\subset \CC^3 \,;
\end{math}
hence, in all its variables $(x,\xi,t)$,
provided the initial values (at $t=0$) are holomorphic in
the complex domain
\begin{math}
    \mathfrak{V}^{(r)}
  = \mathfrak{X}^{(r)}\times \Delta_{\arctan r} \subset \CC^2 .
\end{math}
The case of general initial data $u_0\in H$ will be postponed until
Section~\ref{s:End} where we let the analytic initial data
$u_{0,n}\vert_{\HH}$ converge to arbitrary initial data $u_0$
in~$H$ as $n\to \infty$.
Finally, the convergence of the (unique) holomorphic extensions
to a smaller domain
\begin{equation*}
  \Gamma^{(T')}_T (\kappa_0,\nu_0) \subset
  \mathfrak{V}^{(r)}\times \Delta_{\vartheta'}
\end{equation*}
of the corresponding weak solutions
$u_n\colon \HH\times [0,T]\to \RR$ of the initial value problem
\eqref{e:Cauchy} with $f\equiv 0$ and the initial data
$u_{0,n}\vert_{\HH} \in H$,
obtained in Proposition~\ref{prop-Lions},
to a holomorphic function
$u\colon \Gamma^{(T')}_T (\kappa_0,\nu_0)\colon \CC$
will be established in the next section
(Section~\ref{s:L^2-bound}).
This argument will help us to complete the proof of our main result
(Theorem~\ref{thm-Main}).

Next, we define a few function spaces for functions on
$\mathfrak{V}^{(r)}\subset \CC^2$.
We denote by
$\mathcal{L}^{2,\infty}(\mathfrak{V}^{(r)})$
the Banach space of all complex\--valued, Lebesgue measurable functions
$u\colon \mathfrak{V}^{(r)}\to \CC$, such that,
for each pair $y, \omega\in \RR$ with $|y| < r$ and $|\omega| < r$,
the following integral converges,
\begin{equation}
\label{e:L^2(y,omega)}
  \int_0^{\infty} \int_{-\infty}^{+\infty}
  \left| u
  \left( x + \ii y ,\, \xi (1 + \ii\omega) \right)
  \right|^2\cdot \mathfrak{w}(x,\xi) \,\mathrm{d}x \,\mathrm{d}\xi
  < \infty \,,
\end{equation}
and the norm
\begin{equation}
\label{e:L^infty(V^(r)}
\begin{aligned}
& \| u\|_{ \mathcal{L}^{2,\infty}(\mathfrak{V}^{(r)}) } \eqdef
\\
& \esssup_{ |y| < r ,\ |\omega| < r }
  \left(
  \int_0^{\infty} \int_{-\infty}^{+\infty}
  \left| u
  \left( x + \ii y ,\, \xi (1 + \ii\omega) \right)
  \right|^2\cdot \mathfrak{w}(x,\xi) \,\mathrm{d}x \,\mathrm{d}\xi
  \right)^{1/2}
  < \infty \,.
\end{aligned}
\end{equation}
It is well\--known that
$\mathcal{L}^{2,\infty}(\mathfrak{V}^{(r)})$
is a vector space and
$\|\cdot\|_{ \mathcal{L}^{2,\infty}(\mathfrak{V}^{(r)}) }$
defines a norm on it; cf.\
{\sc P.\ Tak\'a\v{c}} \cite[Sect.~5]{Takac-12}.
It is easy to verify that
$\mathcal{L}^{2,\infty}(\mathfrak{V}^{(r)})$
is a Banach space.
We denote by
$\mathcal{H}^2(\mathfrak{V}^{(r)})$
the Hardy space of all holomorphic functions
$u\colon \mathfrak{V}^{(r)}\to \CC$ such that
$u\in \mathcal{L}^{2,\infty}(\mathfrak{V}^{(r)})$.
It is well\--known that
$\mathcal{H}^2(\mathfrak{V}^{(r)})$
is a closed vector subspace of
$\mathcal{L}^{2,\infty}(\mathfrak{V}^{(r)})$.
We refer to
{\sc E.~M.\ Stein} and {\sc G.\ Weiss} \cite[Chapt.~III]{SteinWeiss}
for basic theory of Hardy spaces;
the most relevant results about $H^2(\mathfrak{V}^{(r)})$ can be found in
\cite[Chapt.~III]{SteinWeiss},
{\S}2, pp.\ 91--101, and {\S}6.12, pp.\ 127--128.

The problem of analyticity (holomorphic extension)
of a weak solution to
the homogeneous Cauchy problem \eqref{e:Cauchy} (with $f\equiv 0$)
can be split into two parts,
{\it analyticity in time\/} and {\it analyticity in space\/};
see {\S}\ref{ss:Anal-compl_t} and {\S}\ref{ss:Cauchy-compl} below,
respectively.
Since the partial differential operator
$\mathcal{A}\colon V\to V'$ in eq.~\eqref{e:Cauchy}
is independent from time $t$,
analyticity in the time variable $t$ follows from
the well\--known theory of analytic $C^0$-semigroups as described below.


\subsection{Analyticity in the complex time variable $t$}
\label{ss:Anal-compl_t}

Our results from the previous section (Section~\ref{s:Heston-real})
on the boundedness and coercivity of the linear operator
$\mathcal{A}\colon V\to V'$ in eq.~\eqref{e:Cauchy}
show that $\mathcal{A}$ is a {\em sectorial operator\/}
in the complex Hilbert space~$H$.
More precisely, the linear operator
${}- \left( \mathcal{A} + \genfrac{}{}{}1{1}{2} c_2'\, I\right)$
in $H$ possesses a bounded inverse, by the Lax\--Milgram theorem, and
ineq.~\eqref{ineq:Heston-diss:RE} implies that
there are constants $\vartheta\in (0, \pi / 2)$ and
$M_{\vartheta}\in (0,\infty)$, such that
\begin{align}
\label{e:Hille-Yosida}
& \| \left( \lambda\, I + \genfrac{}{}{}1{1}{2} c_2' + \mathcal{A}
     \right)^{-1}
  \|_{ \mathcal{L}(H\to H) }
  \leq M_{\vartheta} / |\lambda|
\\
\nonumber
& \;\mbox{ holds for all }\;
    \lambda = \varrho\ee^{\ii\theta} \in \CC
  \;\mbox{ with $\varrho > 0$ and }\,
    \theta\in
    \left( {}- \genfrac{}{}{}1{1}{2}\pi - \vartheta ,\,
               \genfrac{}{}{}1{1}{2}\pi + \vartheta
    \right) \,.
\end{align}
Consequently,
${}- \left( \mathcal{A} + \genfrac{}{}{}1{1}{2} c_2'\, I\right)$
is the {\em infinitesimal generator\/}
of a {\em holomorphic semigroup\/} of uniformly bounded linear operators
\begin{math}
  \left\{
     \ee^{ {}- c_2' t/2 }\,
     \ee^{-t\mathcal{A}}\colon t\in \RR_+
  \right\}
\end{math}
in $H$, i.e.,
\begin{equation}
\label{eq:Hille-Yosida}
  \| \ee^{-t\mathcal{A}} \|_{ \mathcal{L}(H\to H) }
  \leq M_{\vartheta'}^{\prime}\, \ee^{ (c_2'/2)\cdot \RE t }
  \quad\mbox{ holds for all }\; t\in \Delta_{\vartheta'} \,,
\end{equation}
where $\vartheta'\in (0,\vartheta)$ is arbitrary and
$M_{\vartheta'}^{\prime}\in (0,\infty)$ is a suitable constant
depending on $\vartheta'$;
see, e.g.,
Theorem 5.7.2 in
{\sc H.~Tanabe} \cite{Tanabe}, {\S}5.7, p.~161, combined with
\cite[Theorem 5.7.6]{Tanabe}, {\S}5.7.4, p.~179.
This means that the strongly continuous mapping
$t\mapsto \ee^{ {}- c_2' t/2 }\, \ee^{-t\mathcal{A}}$
of $\RR_+$ into the Banach algebra of all bounded linear operators on $H$
(endowed with the operator norm
 $\|\cdot\|_{ \mathcal{L}(H\to H) }$)
can be extended uniquely to a holomorphic mapping in
a {\em complex angle\/} $\Delta_{\vartheta'}$ of angular width
$2\vartheta'$, defined in \eqref{e:angle},
$\vartheta'\in (0, \pi / 2)$ small enough,
$0 < \vartheta' < \vartheta < \pi / 2$.

Hence, the unique weak solution
$u\colon \HH\times [0,T]\to \RR$ of the initial value problem
\eqref{e:Cauchy} with $f\equiv 0$,
obtained in Proposition~\ref{prop-Lions},
extends uniquely to the complex domain
$\HH\times \Delta_{\vartheta'}$ and is holomorphic in the time variable
$t\in \Delta_{\vartheta'}$.
Furthermore, by ineq.~\eqref{eq:Hille-Yosida} above,
the following estimate holds for any initial condition $u_0\in H$,
\begin{equation}
\label{ineq:Hille-Yosida}
    \| u(\,\cdot\,, \,\cdot\,, t) \|_H
  = \| \ee^{-t\mathcal{A}} u_0 \|_H
  \leq M_{\vartheta'}^{\prime}\, \ee^{ (c_2'/2)\cdot \RE t }\, \| u_0\|_H
  \quad\mbox{ for all }\; t\in \Delta_{\vartheta'} \,.
\end{equation}
%


\subsection{The Cauchy problem in the complex domain}
\label{ss:Cauchy-compl}

Given an initial condition $u_0\in H$,
in the Appendix (Appendix~\ref{s:density})
there is a sequence of entire functions
$u_{0,n}\colon \CC^2\to \CC$; $n=1,2,3,\dots$, with
$u_{0,n}\vert_{\HH}\in H$, constructed such that
\begin{equation*}
  \| u_{0,n}\vert_{\HH} - u_0\|_H \;\longrightarrow\; 0
    \quad\mbox{ as }\, n\to \infty \,.
\end{equation*}
An important property of each function
$u_{0,n}\colon \CC^2\to \CC$ is the following decay inequality:
Given any numbers $r\in (0,\infty)$ and $\vartheta\in (0, \pi / 2)$,
for each $n=1,2,3,\dots$, there exists a constant
$A_n\equiv A_n(r,\vartheta)\in (0,\infty)$ such that
\begin{align}
\label{est:H_n-L_n}
& | u_{0,n}(x + \ii y, \xi + \ii\eta) |
  \leq A_n\, \ee^{ {}- (x^2 + \xi) / 4 }
\\
\nonumber
& \;\mbox{ whenever }\;
  z = x + \ii y\in \mathfrak{X}^{(r)} \,\mbox{ and }\,
  \zeta = \xi + \ii\eta\in \Delta_{\vartheta} \,,
\end{align}
where the right\--hand side is in
$H = L^2(\mathbb{H};\mathfrak{w})$.

To begin with, let us fix an arbitrary index $n\in \NN$;
$\NN\eqdef \{ 1,2,3,\dots\}$,
for which we abbreviate $u_0\equiv u_{0,n}$ with
$u_0\vert_{\HH}\in H$.
Hence, throughout this paragraph we assume that either
$u_0\colon \CC^2\to \CC$ is an entire function or at least
$u_0\colon \mathfrak{X}^{(r)}\times \Delta_{\vartheta} \to \CC$
is a holomorphic function that satisfies
an analogue of \eqref{est:H_n-L_n} with a constant
$A_0\equiv A_0(r,\vartheta)\in (0,\infty)$:
\begin{align}
\label{est:H_0-L_0}
& | u_0(x + \ii y, \xi + \ii\eta) |
  \leq A_0\, \ee^{ {}- (x^2 + \xi) / 4 }
\\
\nonumber
& \;\mbox{ whenever }\;
  z = x + \ii y\in \mathfrak{X}^{(r)} \,\mbox{ and }\,
  \zeta = \xi + \ii\eta\in \Delta_{\vartheta} \,.
\end{align}
To simplify our hypotheses and notation, we take
$r\in (0,\infty)$ arbitrary and
$\vartheta = \arctan r\in (0, \pi / 2)$, so that
\begin{math}
    \mathfrak{X}^{(r)}\times \Delta_{\vartheta}
  = \mathfrak{V}^{(r)}\subset \CC^2
\end{math}
is the complex domain
\begin{math}
  \mathfrak{V}^{(r)} =
  \mathfrak{X}^{(r)} \times \Delta_{\arctan r}\subset \CC^2
\end{math}
that has been defined in eq.~\eqref{x,v:strip}.
The general case of $u_0\in H$ will be treated in the next section
(Section~\ref{s:L^2-bound}).

We formulate the corresponding analyticity result for
such an initial condition $u_0$ as the following special case of
Theorem~\ref{thm-Main}:

\begin{proposition}\label{prop-Main}
Let\/
$\rho$, $\sigma$, $\theta$, $q_r$, and\/ $\gamma$
be given constants in $\RR$,
$\rho\in (-1,1)$, $\sigma > 0$, $\theta > 0$, and\/ $\gamma > 0$.
Assume that\/ $\beta$, $\gamma$, $\kappa$, and\/ $\mu$
are chosen as specified in {\rm Proposition~\ref{prop-Lions}}.
Finally, let us assume that\/
$u_0\colon \mathfrak{V}^{(r)}\to \CC$
is a holomorphic function that satisfies
a bound similar to \eqref{est:H_0-L_0},
\begin{align}
\label{est:bound-u_0}
& | u_0(x + \ii y, \xi + \ii\eta) |
  \leq A_0\, \ee^{ {}- (x^2 + \xi) / 4 }
\\
\nonumber
& \;\mbox{ whenever }\;
  z = x + \ii y\in \mathfrak{X}^{(r)} \,\mbox{ and }\,
  \zeta = \xi + \ii\eta\in \Delta_{\arctan r} \,,
\end{align}
where $r\in (0,\infty)$ is some number and\/
$A_0\equiv A_0(r)\in (0,\infty)$ is a constant.

Then the (unique) weak solution
\begin{equation*}
  u\in C([0,T]\to H)\cap L^2((0,T)\to V)
\end{equation*}
of the homogeneous initial value problem \eqref{e:Cauchy}
(with $f\equiv 0$ and this $u_0$)
possesses a unique holomorphic extension
\begin{math}
  \tilde{u}\colon
  \mathfrak{V}^{(r')}\times \Delta_{\vartheta'} \to \CC
\end{math}
to the complex domain
$\mathfrak{V}^{(r')}\times \Delta_{\vartheta'} \subset \CC^3$,
where $r'\in (0,r]$ and\/
$\vartheta'\in (0, \pi / 2)$ are some constants.
Furthermore, there are additional constants
$C_0, c_0\in \RR_+$ such that\/
\begin{align}
\label{eq:u_Hardy^2}
\begin{aligned}
& \int_0^{\infty} \int_{-\infty}^{+\infty}
  \left| \tilde{u}
  \left( x + \ii y ,\, \xi (1 + \ii\omega) ,\, t\right)
  \right|^2\cdot \mathfrak{w}(x,\xi) \,\mathrm{d}x \,\mathrm{d}\xi
\\
& \leq C_0\, \ee^{ c_0\cdot \RE t }\cdot
  \int_0^{\infty} \int_{-\infty}^{+\infty}
  \left| u_0\left( x + \ii y ,\, \xi (1 + \ii\omega) \right)
  \right|^2\cdot \mathfrak{w}(x,\xi) \,\mathrm{d}x \,\mathrm{d}\xi
\end{aligned}
\end{align}
for every\/ $t\in \Delta_{\vartheta'}$ and for all\/
$y, \omega\in \RR$ such that\/
$|y| < r'$ and\/ $|\omega| < r'$.
\end{proposition}
\par\vskip 10pt

Before giving the {\em proof\/} of this proposition,
we make a few important remarks:
The proof hinges upon the fact that if the holomorphic extension
\begin{math}
  \tilde{u}\colon
  \mathfrak{V}^{(r')}\times \Delta_{\vartheta'} \to \CC
\end{math}
of a weak solution
\begin{equation*}
  u\in C([0,T]\to H)\cap L^2((0,T)\to V)
\end{equation*}
of the homogeneous initial value problem \eqref{e:Cauchy} exists,
then it must satisfy the following initial value problem
with complex partial derivatives:
\begin{equation}
\label{e:Cauchy-compl}
\left\{
\begin{alignedat}{2}
  \frac{\partial\tilde{u}}{\partial t}
  + ( \tilde{\mathcal{A}}\tilde{u} ) (z,\zeta,t) &= 0
  &&\quad\mbox{ in }\,
    \mathfrak{V}^{(r')}\times \Delta_{\vartheta'} \,;
\\
  \tilde{u}(z,\zeta,0) &= u_0(z,\zeta)
  &&\quad\mbox{ for }\, (z,\zeta)\in \mathfrak{V}^{(r')} \,,
\end{alignedat}
\right.
\end{equation}
where the complex partial differential operator
$\tilde{\mathcal{A}}$ is given by
\begin{align}
\label{eq_C:Heston-oper}
&
\begin{alignedat}{2}
& (\tilde{\mathcal{A}}\tilde{u})(z,\zeta) =
&&{} - \frac{1}{2}\, \sigma\zeta\cdot
  \left[ \frac{\partial}{\partial z}
  \left(
    \frac{\partial\tilde{u}}{\partial z}(z,\zeta)
  + 2\rho\, \frac{\partial\tilde{u}}{\partial\zeta}(z,\zeta)
  \right)
  + \frac{\partial^2\tilde{u}}{\partial\zeta^2}(z,\zeta)
  \right]
\\
&
&&{}
  + \left( q_r + \genfrac{}{}{}1{1}{2} \sigma\zeta\right)
    \cdot \frac{\partial\tilde{u}}{\partial z}(z,\zeta)
  - \kappa (\theta_{\sigma} - \zeta)
    \cdot \frac{\partial\tilde{u}}{\partial\zeta}(z,\zeta)
\end{alignedat}
\\
\nonumber
&
\begin{alignedat}{2}
& \equiv\;
&&{} - \frac{1}{2}\, \sigma\zeta\cdot
  \left[
  \left( \tilde{u}_z + 2\rho\, \tilde{u}_{\zeta}\right)_z
       + \tilde{u}_{\zeta\zeta}
  \right]
  + \left( q_r + \genfrac{}{}{}1{1}{2} \sigma\zeta\right)
    \cdot \tilde{u}_z
  - \kappa (\theta_{\sigma} - \zeta) \cdot \tilde{u}_{\zeta}
\\
&
&&\quad\mbox{ for }\;
  (z,\zeta)\in \mathfrak{V}^{(r')}
             = \mathfrak{X}^{(r')} \times \Delta_{\arctan r'} \,.
\end{alignedat}
\end{align}
This operator has been obtained from
the Heston operator \eqref{eq:Heston-oper}
by the natural complexification of the variables $x$ and $\xi$ as
$z = x + \ii y$ and $\zeta = \xi + \ii\eta$, respectively,
with the imaginary parts $y, \eta$ $\in \RR$.
However, to establish the conclusion of
Proposition~\ref{prop-Main},
we need to choose the imaginary parts $y,\eta\in \RR$ such that
$|y| < r'$ and $\eta = \xi\omega$ with $|\omega| < r'$,
where $y$ and $\omega$ are fixed, while
$x$ and $\xi$ are the independent variables, $(x,\xi)\in \HH$.
Hence, we have to investigate the function
\begin{equation}
\label{e:u(x,xi,t)}
\begin{aligned}
& v\colon
  (x,\xi,t) \;\longmapsto\; v(x,\xi,t) \equiv
  v_{(\ii y + z^{*})}^{(\ii\omega + \omega^{*})} (x,\xi,t)
\\
& \eqdef \tilde{u}
    \Bigl(
    x + \ii y + z^{*} ,\, \xi (1 + \ii\omega + \omega^{*}) ,\, t
    \Bigr)
  \colon \HH\times \Delta_{\vartheta'}\to \CC
\end{aligned}
\end{equation}
with the complexified space variables
\begin{equation}
\label{e_compl:u(x,xi,alpha)}
\begin{aligned}
    z + z^{*}
& = x + \ii y + z^{*}
  = x + x^{*} + \ii (y + y^{*}) \,,
\\
    \zeta + \zeta^{*}
& = \xi (1 + \ii\omega) + \zeta^{*}
  = \xi (1 + \ii\omega + \omega^{*}) \,.
\end{aligned}
\end{equation}
Here,
$z^{*}, \omega^{*}\in \CC$
are complex numbers with sufficiently small absolute values, such that
\begin{equation}
\label{e:z^*,zeta^*,t^*}
  \ii y + z^{*}\in \mathfrak{X}^{(r')}
    \quad\mbox{ and }\quad
  1 + \ii\omega + \omega^{*}\in \Delta_{\arctan r'} \,,
\end{equation}
which guarantees that the argument of the function $\tilde{u}$
in eq.~\eqref{e:u(x,xi,t)} above stays in
$\mathfrak{V}^{(r')}\times \Delta_{\vartheta'}$
for all
$(x,\xi,t)\in \HH\times \Delta_{\vartheta'}$.
Small complex perturbations
$(z^{*}, \omega^{*})\in \CC^2$
are needed to calculate partial derivatives
of the function $\tilde{u}(z,\zeta,t)$ with respect to
the real and imaginary parts of its arguments
$(z,\zeta)\in \mathfrak{V}^{(r')}$.
The complex differentiability (yielding the holomorphy)
with respect to the time variable $t\in \Delta_{\vartheta'}$
has been treated in the previous paragraph
({\S}\ref{ss:Anal-compl_t}).

A simple application of the chain rule,
\begin{equation*}
    \frac{\partial v}{\partial x} (x,\xi,t)
  = \frac{\partial\tilde{u}}{\partial z}
    \left( z + z^{*} ,\, \zeta + \zeta^{*} ,\, t\right)
    \quad\mbox{ and }\quad
    \frac{\partial v}{\partial\xi}
  = (1 + \ii\omega + \omega^{*})\,
    \frac{\partial\tilde{u}}{\partial\zeta} \,,
\end{equation*}
shows that the function
$v\colon \HH\times \Delta_{\vartheta'}\to \CC$
defined in eq.~\eqref{e:u(x,xi,t)} must be a weak solution to
the following initial value problem with real partial derivatives:
\begin{equation}
\label{e:Cauchy-real}
\left\{
\begin{alignedat}{2}
  \frac{\partial v}{\partial t}
& {}
  + \left( \mathcal{A}^{(\ii\omega + \omega^{*})} v\right) (x,\xi,t)
  = 0
  &&\quad\mbox{ in }\, \HH\times \Delta_{\vartheta'} \,;
\\
  v(x,\xi,0) &{} = u_0
    \left( x + \ii y + z^{*} ,\, \xi (1 + \ii\omega + \omega^{*})
    \right)
  &&\quad\mbox{ for }\, (x,\xi)\in \HH \,,
\end{alignedat}
\right.
\end{equation}
where the real partial differential operator
$\mathcal{A}^{(\ii\omega + \omega^{*})}$ is given by
\begin{align*}
& \left( \mathcal{A}^{(\ii\omega + \omega^{*})} v\right) (x,\xi) =
\\
&
\begin{aligned}
  {} - \frac{1}{2}\, (1 + \ii\omega + \omega^{*}) \sigma\xi\cdot
& \left[ \frac{\partial}{\partial x}
  \left(
    \frac{\partial v}{\partial x}(x,\xi)
  + \frac{2\rho}{1 + \ii\omega + \omega^{*}}\cdot
            \frac{\partial v}{\partial\xi}(x,\xi)
  \right)
  \right.
\\
& \left.
{}+ \frac{1}{ (1 + \ii\omega + \omega^{*})^2 }\cdot
    \frac{\partial^2 v}{\partial\xi^2}(x,\xi)
  \right]
\end{aligned}
\\
&
\begin{aligned}
& {}
  + \left[ q_r
         + \genfrac{}{}{}1{1}{2} (1 + \ii\omega + \omega^{*}) \sigma\xi
    \right]
    \cdot \frac{\partial v}{\partial x}(x,\xi)
\\
& {}
  - \frac{\kappa}{1 + \ii\omega + \omega^{*}}\,
      [ \theta_{\sigma} - (1 + \ii\omega + \omega^{*}) \xi ]
    \cdot \frac{\partial v}{\partial\xi}(x,\xi)
\end{aligned}
\\
&
\begin{alignedat}{2}
& \equiv\;
&&{} - \frac{1}{2}\, \sigma\xi\cdot
  \left[
  \left( (1 + \ii\omega + \omega^{*}) v_x + 2\rho\, v_{\xi}\right)_x
       + (1 + \ii\omega + \omega^{*})^{-1} v_{\xi\xi}
  \right]
\\
&
&&{} + \left[ q_r
     + \genfrac{}{}{}1{1}{2} (1 + \ii\omega + \omega^{*}) \sigma\xi
       \right] \cdot v_x
  {} - \kappa
       \left[ (1 + \ii\omega + \omega^{*})^{-1} \theta_{\sigma} - \xi
       \right] \cdot v_{\xi}
\\
&
&&\quad\mbox{ for $(x,\xi)\in \HH$. }
\end{alignedat}
\end{align*}
Consequently, recalling the definition of $\mathcal{A}$
in eq.~\eqref{eq:Heston-oper}, we have
\begin{equation}
\label{eq_R-A:Heston-oper}
\begin{aligned}
& \left( \mathcal{A}^{(\ii\omega + \omega^{*})} v\right) (x,\xi)
  = (\mathcal{A} v) (x,\xi)
\\
& {} - \frac{\sigma}{2}\, (\ii\omega + \omega^{*}) \xi\cdot
  \left( v_{xx} - (1 + \ii\omega + \omega^{*})^{-1} v_{\xi\xi}
  \right)
\\
& {} + \frac{\sigma}{2}\, (\ii\omega + \omega^{*}) \xi\cdot v_x
     + \frac{ \ii\omega + \omega^{*} }{ 1 + \ii\omega + \omega^{*} }\,
       \kappa\theta_{\sigma}\cdot v_{\xi}
  \quad\mbox{ for $(x,\xi)\in \HH$. }
\end{aligned}
\end{equation}
It is important to note that the linear operator
$\mathcal{A}^{(\ii\omega + \omega^{*})} \colon V\to V'$
does {\em not \/} depend on $y\in \RR$ or $z^{*}\in \CC$.
However, it does depend on $\omega\in \RR$ and $\omega^{*}\in \CC$;
more precisely, it depends on the sum $\ii\omega + \omega^{*}$.

To derive the sesquilinear form associated to
$\mathcal{A}^{(\ii\omega + \omega^{*})}$,
\begin{equation}
\label{eq_R:Heston-bilin}
  (v,w)\mapsto
  \left( \mathcal{A}^{(\ii\omega + \omega^{*})} v, w\right)_H \,,
\end{equation}
we apply the same methods as for obtaining eq.~\eqref{e:Heston-diss:u=w}
associated to $\mathcal{A}$.
We thus arrive at
\begin{align*}
\nonumber
& \left( \mathcal{A}^{(\ii\omega + \omega^{*})} v, w\right)_H
  = (\mathcal{A}v, w)_H
\\
&
\begin{aligned}
& {}+ \frac{\sigma}{2}\, (\ii\omega + \omega^{*}) \int_{\HH}
  \left(
    v_x\cdot \bar{w}_x
  - (1 + \ii\omega + \omega^{*})^{-1} v_{\xi}\cdot \bar{w}_{\xi}
  \right) \cdot \xi\cdot \mathfrak{w}(x,\xi)
    \,\mathrm{d}x \,\mathrm{d}\xi
\end{aligned}
\\
\nonumber
&
\begin{aligned}
  {}- \frac{\sigma}{2}\, (\ii\omega + \omega^{*}) \int_{\HH}
& \bigl[
      \gamma\, \Sgn x\cdot v_x\, \bar{w}\cdot \xi
\\
&
{}+ (1 + \ii\omega + \omega^{*})^{-1}
      (\beta - \mu\xi)\, v_{\xi}\cdot \bar{w}
  \bigr] \mathfrak{w}(x,\xi)
    \,\mathrm{d}x \,\mathrm{d}\xi
\end{aligned}
\\
\nonumber
&
\begin{aligned}
& {}+ \frac{\sigma}{2}\, (\ii\omega + \omega^{*}) \int_{\HH}
      v_x\, \bar{w}\cdot \xi\cdot
    \mathfrak{w}(x,\xi) \,\mathrm{d}x \,\mathrm{d}\xi
\\
& {}+ \frac{ \ii\omega + \omega^{*} }{ 1 + \ii\omega + \omega^{*} }\,
      \kappa\theta_{\sigma}
      \int_{\HH} v_{\xi}\, \bar{w}\cdot \mathfrak{w}(x,\xi)
    \,\mathrm{d}x \,\mathrm{d}\xi \,,
\end{aligned}
\end{align*}
where we have taken advantage of the vanishing boundary conditions
\eqref{bc_xi:Heston-bilin} and \eqref{bc_x:Heston-bilin}
with the pair of functions $(v,w)$ in place of $(u,w)$,
while performing integration\--by\--parts
on the second summand on the right\--hand side of
eq.~\eqref{eq_R-A:Heston-oper}; cf.\ also
eqs.\ \eqref{e:trace:v=0}, \eqref{e:trace:v=infty}, and 
\eqref{e:trace:x=+-infty}.

Finally, the sesquilinear form \eqref{eq_R:Heston-bilin} becomes
\begin{align}
\nonumber
& \left( \mathcal{A}^{(\ii\omega + \omega^{*})} v, w\right)_H
  = (\mathcal{A}v, w)_H
\\
\label{eq_R:Heston-diss:u=w}
&
\begin{aligned}
& {}+ \frac{\sigma}{2}\, (\ii\omega + \omega^{*}) \int_{\HH}
  \left(
    v_x\cdot \bar{w}_x
  - (1 + \ii\omega + \omega^{*})^{-1} v_{\xi}\cdot \bar{w}_{\xi}
  \right) \cdot \xi\cdot \mathfrak{w}(x,\xi)
    \,\mathrm{d}x \,\mathrm{d}\xi
\end{aligned}
\\
\nonumber
&
\begin{aligned}
& {}+ \frac{\sigma}{2}\, (\ii\omega + \omega^{*}) \int_{\HH}
  (1 - \gamma\, \Sgn x)\, v_x\cdot \bar{w}
    \cdot \xi\cdot \mathfrak{w}(x,\xi)
    \,\mathrm{d}x \,\mathrm{d}\xi
\\
& {}
    + \frac{\sigma}{2}\cdot
      \frac{ \ii\omega + \omega^{*} }{ 1 + \ii\omega + \omega^{*} }\,
      \mu \int_{\HH}
  v_{\xi}\cdot \bar{w}\cdot \xi\cdot \mathfrak{w}(x,\xi)
    \,\mathrm{d}x \,\mathrm{d}\xi
\\
& {}- \frac{ \ii\omega + \omega^{*} }{ 1 + \ii\omega + \omega^{*} }\,
      \left( \genfrac{}{}{}1{1}{2} \beta\sigma - \kappa\theta_{\sigma}
      \right)
    \int_{\HH} v_{\xi}\cdot \bar{w}\cdot \mathfrak{w}(x,\xi)
    \,\mathrm{d}x \,\mathrm{d}\xi \,.
\end{aligned}
\end{align}
All integrals on the right\--hand side converge absolutely for any pair
$u,w\in V$, in analogy with eq.~\eqref{e:Heston-diss:u=w}.
In what follows we use the last formula,
eq.~\eqref{eq_R:Heston-diss:u=w},
to define the sesquilinear form \eqref{eq_R:Heston-bilin}
in $V\times V$.

The following two results, respectively, are analogues of
Propositions \ref{prop-bound-R} and \ref{prop-coerc-R}
with similar proofs.
Here, the sesquilinear form from eq.~\eqref{eq_R:Heston-diss:u=w}
replaces that from \eqref{e:Heston-diss:u=w}.
We use the former to verify the boundedness and coercivity of
the linear operator
$\mathcal{A}^{(\ii\omega + \omega^{*})} \colon V\to V'$
in the Hilbert space $V = H^1(\mathbb{H};\mathfrak{w})$.
The details of these proofs are left to an interested reader.

\begin{proposition}\label{prop-bound-C}
{\rm (Boundedness.)}$\;$
Let\/
$\beta$, $\gamma$, $\mu$, $\rho$, $\sigma$, $\theta$, $q_r$, and\/
$\kappa$ be given constants in $\RR$,
$\beta > 1$, $\gamma > 0$, $\mu > 0$, $-1 < \rho < 1$, $\sigma > 0$,
and\/ $\theta > 0$.
Then, given any number\/ $r\in (0,\infty)$,
there exists a constant\/
$C^{*}\in (0,\infty)$, such that, for all numbers\/
$\omega\in (-r,r)$ and\/ $\omega^{*}\in \CC$ with
$|\omega^{*}|\leq 1/2$, and for all pairs $u,w\in V$, we have
\begin{equation}
\label{est_C:Heston-diss}
  \left|
  \left( \mathcal{A}^{(\ii\omega + \omega^{*})} u, w\right)_H
  \right| \leq C^{*}\cdot \| u\|_V\cdot \| w\|_V \,.
\end{equation}
\end{proposition}

In our next proposition, the number $r\in (0,\infty)$
has to be sufficiently small, unlike in the analogous
Proposition~\ref{prop-coerc-R} where it is arbitrary.

\begin{proposition}\label{prop-coerc-C}
{\rm (Coercivity.)}$\;$
Let\/
$\rho$, $\sigma$, $\theta$, $q_r$, and\/
$\gamma$ be given constants in $\RR$,
$\rho\in (-1,1)$, $\sigma > 0$, $\theta > 0$, and\/ $\gamma > 0$.
Assume that\/ $\beta$, $\gamma$, $\kappa$, and\/ $\mu$
are chosen as specified in {\rm Proposition~\ref{prop-Lions}}.
Then there exist constants\/
$r\in \bigl( 0 ,\, \frac{1}{2}\bigr]$
and\/ $c_2''\in (0,\infty)$
such that the following {\em G\r{a}rding inequality\/}
\begin{equation}
\label{ineq:Heston-C:RE}
  2\cdot \RE
  \left( \mathcal{A}^{(\ii\omega + \omega^{*})} u, u\right)_H
  \geq \frac{\sigma}{2}\, (1 - |\rho|)\cdot
       \| u\|_V^2 - c_2''\cdot \| u\|_H^2
\end{equation}
is valid for all\/
$\omega\in (-r,r)$ and\/ $\omega^{*}\in \CC$ with
$|\omega^{*}|\leq r$, and for all\/ $u\in V$.
\end{proposition}
\par\vskip 10pt

Now we are ready to prove Proposition~\ref{prop-Main}.

\par\vskip 10pt
{\it Proof of\/} {\bf Proposition~\ref{prop-Main}.}$\;$
It is obvious that we must find a method how to solve
the initial value problem \eqref{e:Cauchy-real}
with a conclusion similar to that provided in paragraph
{\S}\ref{ss:Anal-compl_t}
for the initial value problem
\eqref{e:Cauchy} with $f\equiv 0$, thanks to
Propositions \ref{prop-bound-R} and \ref{prop-coerc-R}
for the linear operator $\mathcal{A}\colon V\to V'$.
Notice that the initial condition in problem \eqref{e:Cauchy-real}
reads
\begin{equation}
\label{e:Cauchy-real:t=0}
  v(x,\xi,0) = v_0(x,\xi)\eqdef u_0
    \left( x + \ii y + z^{*} ,\, \xi (1 + \ii\omega + \omega^{*})
    \right)
    \quad\mbox{ for }\, (x,\xi)\in \HH \,.
\end{equation}
Thus, we must first adapt these two propositions to the linear operator
$\mathcal{A}^{(\ii\omega + \omega^{*})} \colon V\to V'$
for any fixed numbers $y,\omega\in \RR$ with
$|y| < r'$ and $|\omega| < r'$, and for any fixed complex numbers
$z^{*}, \omega^{*}\in \CC$
with sufficiently small absolute values, such that
\eqref{e:z^*,zeta^*,t^*} holds.
It suffices to do this for some $r'\in (0,r]$ small enough.
Hence, the couple
\begin{math}
  \left( z + z^{*} ,\, \zeta + \zeta^{*} \right)
\end{math}
from eq.~\eqref{e_compl:u(x,xi,alpha)}
that appears also as the argument of the function $u_0$
in eq.~\eqref{e:Cauchy-real:t=0} above stays in
$\mathfrak{V}^{(r')} \subset \mathfrak{V}^{(r)}$
for all $(x,\xi)\in \HH$, thanks to $0 < r'\leq r$.

In analogy with
Propositions \ref{prop-bound-R} and~\ref{prop-coerc-R}
(boundedness and coercivity, respectively)
for the operator $\mathcal{A}\colon V\to V'$,
Propositions \ref{prop-bound-C} and~\ref{prop-coerc-C}
(Appendix~\ref{s:Trace,Sobolev})
for the operator
$\mathcal{A}^{(\ii\omega + \omega^{*})} \colon V\to V$
guarantee that
$\mathcal{A}^{(\ii\omega + \omega^{*})}$
is a {\em sectorial operator\/}
in the Hilbert space $H$, provided
$|\omega| < r'$ and $|\omega^{*}|$ is small enough.
Hence,
${}- \mathcal{A}^{(\ii\omega + \omega^{*})}$
is the {\em infinitesimal generator\/}
of a {\em holomorphic semigroup\/} of bounded linear operators
\begin{math}
  \left\{
     \ee^{ -t\mathcal{A}^{(\ii\omega + \omega^{*})} } \colon t\in \RR_+
  \right\}
\end{math}
in $H$, i.e.,
\begin{equation}
\label{eq_C:Hille-Yosida}
  \| \ee^{ -t\mathcal{A}^{(\ii\omega + \omega^{*})} }
  \|_{ \mathcal{L}(H\to H) }
  \leq M_{\vartheta''}^{\prime\prime}\, \ee^{ (c_2''/2)\cdot \RE t }
  \quad\mbox{ holds for all }\; t\in \Delta_{\vartheta''} \,,
\end{equation}
where $\vartheta''\in (0,\vartheta)$ is arbitrary and
$M_{\vartheta''}^{\prime\prime} ,\, c_2''\in (0,\infty)$
are suitable constants depending on $\vartheta''$,
but independent from the particular choice of
$\omega\in \RR$ or $\omega^{*}\in \CC$ such that
$|\omega| < r'$ and $|\omega^{*}|$ is small enough.
This semigroup provides the (unique) holomorphic extension
$v\colon \Delta_{\vartheta''}\to H$
of the (unique) weak solution
\begin{equation*}
  v\equiv v_{(\ii y + z^{*})}^{(\ii\omega + \omega^{*})}
  \in C([0,T]\to H)\cap L^2((0,T)\to V)
\end{equation*}
to the initial value problem \eqref{e:Cauchy-real}.
The uniqueness guarantees that this solution depends on
the fixed data $y, \omega\in \RR$ and $z^{*}, \omega^{*}\in \CC$
only through the sums
$\ii y + z^{*}$ and $\ii\omega + \omega^{*}$,
as so do the operator
$\mathcal{A}^{(\ii\omega + \omega^{*})}$
(which, in fact, is independent from $y$ and $z^{*}$)
and the initial condition \eqref{e:Cauchy-real:t=0}.
Indeed, let
$y_j, \omega_j\in \RR$ and $z_j^{*}, \omega_j^{*}\in \CC$
satisfy \eqref{e:z^*,zeta^*,t^*} for both $j=1,2$, i.e.,
\begin{equation}
\label{e_j:z^*,zeta^*,t^*}
  \ii y_j + z_j^{*}\in \mathfrak{X}^{(r')}
    \quad\mbox{ and }\quad
  1 + \ii\omega_j + \omega_j^{*}\in \Delta_{\arctan r'} \,.
\end{equation}
Consider the corresponding (unique) weak solution
\begin{equation*}
  v^{(j)}\equiv v_{(\ii y_j + z_j^{*})}^{(\ii\omega_j + \omega_j^{*})}
  \in C([0,T]\to H)\cap L^2((0,T)\to V)
\end{equation*}
to the initial value problem \eqref{e:Cauchy-real}
together with its (unique) holomorphic extension
$v^{(j)}\colon \Delta_{\vartheta''}$ $\to H$; $j=1,2$.
The initial condition \eqref{e:Cauchy-real:t=0} is given by
\begin{align}
\label{e_j:Cauchy-real:t=0}
  v^{(j)}(x,\xi,0) = v_0^{(j)}(x,\xi)\eqdef u_0
    \left( x + \ii y_j + z_j^{*} ,\,
           \xi (1 + \ii\omega_j + \omega_j^{*})
    \right)
\\
\nonumber
    \quad\mbox{ for }\, (x,\xi)\in \HH \,.
\end{align}
Consequently, if
\begin{equation*}
  \ii y_1 + z_1^{*} = \ii y_2 + z_2^{*}
    \quad\mbox{ and }\quad
  \ii\omega_1 + \omega_1^{*} = \ii\omega_2 + \omega_2^{*} \,,
\end{equation*}
then $v_0^{(1)} = v_0^{(2)}$ in $H$ and, therefore,
the uniqueness for problem \eqref{e:Cauchy-real} forces
$v^{(1)}(x,\xi,t)\equiv v^{(2)}(x,\xi,t)$ for
$(x,\xi,t)\in \HH\times \Delta_{\vartheta''}$.
This uniqueness result allows us to give
the following (correct) definition of a function
\begin{math}
  \tilde{u}\colon \mathfrak{V}^{(r')}\times \Delta_{\vartheta''} \to \CC
\end{math}
by the formula
\begin{align}
\label{def:u(x,xi,t)}
  \tilde{u}
    \Bigl(
    x + \ii y + z^{*} ,\, \xi (1 + \ii\omega + \omega^{*}) ,\, t
    \Bigr)
  \eqdef v_{(\ii y + z^{*})}^{(\ii\omega + \omega^{*})} (x,\xi,t)
\\
\nonumber
  \quad\mbox{ for all $(x,\xi)\in \HH$ and for all }\,
              t\in \Delta_{\vartheta''} \,.
\end{align}
Notice that it suffices to take $z^{*} = \omega^{*} = 0$
and arbitrary numbers $y, \omega\in \RR$ with
$|y| < r'$ and $|\omega| < r'$ to define $\tilde{u}$.

The function
\begin{equation*}
  t\mapsto v_{(\ii y + z^{*})}^{(\ii\omega + \omega^{*})} (x,\xi,t)
  \colon \Delta_{\vartheta''} \to \CC
\end{equation*}
being holomorphic, by {\S}\ref{ss:Anal-compl_t},
it is obvious that also
\begin{math}
  \tilde{u}\colon \mathfrak{V}^{(r')}\times \Delta_{\vartheta''}
  \to \CC
\end{math}
is holomorphic in the time variable $t\in \Delta_{\vartheta''}$.
Furthermore, the estimate in \eqref{eq:u_Hardy^2}
follows immediately from inequality \eqref{eq_C:Hille-Yosida}
by taking
$C_0 = M_{\vartheta''}^{\prime\prime} > 0$ and
$c_0 = c_2''/2 > 0$.

Taking advantage of the differentiability
of the coefficients of the partial differential operator
$\mathcal{A}^{(\ii\omega + \omega^{*})}$
in eq.~\eqref{eq_R-A:Heston-oper},
we observe that if the initial data
$u_0\in \mathcal{L}^{2,\infty}(\mathfrak{V}^{(r)})$
are $C^{\infty}$-smooth (in the real\--variable sense)
then also the (unique) solution
\begin{math}
  \tilde{u}(\,\cdot\,,\,\cdot\, t) \colon
  \mathfrak{V}^{(r')}\to \CC
\end{math}
to the initial value problem \eqref{e:Cauchy-real}
is $C^{\infty}$-smooth in $\HH$, by
Theorem~19 and Corollary (to Theorem~19) in
{\sc A.\ Friedman} \cite[Chapt.~10]{Friedman-64},
on p.~321 and p.~322, respectively.

Now we take advantage of the holomorphic data $v_0$
in the initial condition \eqref{e:Cauchy-real:t=0}
with respect to the small complex parameters
$(z^{*},\omega^{*})\in \CC^2$
in order to show that, for each fixed
$t\in \Delta_{\vartheta'}$, the function
\begin{math}
  \tilde{u}(\,\cdot\,,\,\cdot\, t) \colon
  \mathfrak{V}^{(r')}\to \CC
\end{math}
is holomorphic.
To this end we first realize that the initial data
$v_0$ in \eqref{e:Cauchy-real:t=0},
which depend on the real parameters
$x^{*} = \RE z^{*}$, $y^{*} = \IM z^{*}$,
$\alpha^{*} = \RE\omega^{*}$, and $\beta^{*} = \IM\omega^{*}$,
are continuously differentiable 
(i.e., $C^1$-smooth in the real\--variable sense)
with respect to these parameters.
We wish to prove that the same is true of each function
$v_{(\ii y + z^{*})}^{(\ii\omega + \omega^{*})}$
with respect to the parameters
$x^{*}, y^{*}, \alpha^{*}, \beta^{*}\in \RR$.

In order to be able to apply well\--known results from
{\sc D.~Henry} \cite[Chapt.~3, {\S}4]{Henry}
on the continuous dependence and differentiability of the solution
$v_{(\ii y + z^{*})}^{(\ii\omega + \omega^{*})}$
with respect to parameters, we rewrite
the initial value problem \eqref{e:Cauchy-real} equivalently as
\begin{equation}
\label{e:Cauchy:u-u_0}
\left\{
\begin{alignedat}{2}
  \frac{\partial w}{\partial t}
{}+ \left( \mathcal{A}^{(\ii\omega + \omega^{*})} w\right) (x,\xi,t)
  &= {}
  - \left( \mathcal{A}^{(\ii\omega + \omega^{*})} v_0\right) (x,\xi)
  &&\quad\mbox{ in }\, \HH\times \Delta_{\vartheta'} \,;
\\
  w(x,\xi,0) &{} = 0
  &&\quad\mbox{ for }\, (x,\xi)\in \HH \,,
\end{alignedat}
\right.
\end{equation}
where
\begin{align}
\label{def:w=u-u_0}
    w(x,\xi,t)\equiv 
    w_{(\ii y + z^{*})}^{(\ii\omega + \omega^{*})} (x,\xi,t)
  \eqdef 
    v_{(\ii y + z^{*})}^{(\ii\omega + \omega^{*})} (x,\xi,t)
  - v_0(x,\xi,t) \equiv
\\
\nonumber
  \tilde{u}
    \Bigl(
    x + \ii y + z^{*} ,\, \xi (1 + \ii\omega + \omega^{*}) ,\, t
    \Bigr)
  - u_0
    \Bigl(
    x + \ii y + z^{*} ,\, \xi (1 + \ii\omega + \omega^{*})
    \Bigr)
\end{align}
is the new unknown function of
$(x,\xi,t)\in \HH\times \Delta_{\vartheta'}$.
It is easy to see that the function
\begin{equation*}
  {}-
  \left( \mathcal{A}^{(\ii\omega + \omega^{*})} v_0\right) (x,\xi)
  = {}- (\tilde{A} u_0)
    \Bigl(
    x + \ii y + z^{*} ,\, \xi (1 + \ii\omega + \omega^{*})
    \Bigr)
\end{equation*}
of $(z^{*},\omega^{*})\in \CC$ is holomorphic, for
$|z^{*}|$ and $|\zeta^{*}|$ small enough; hence,
$C^1$-smooth with respect to the real parameters
$x^{*} = \RE z^{*}$, $y^{*} = \IM z^{*}$,
$\alpha^{*} = \RE\omega^{*}$, and $\beta^{*} = \IM\omega^{*}$.
By {\sc Henry}'s theorem \cite[Theorem 3.4.4, pp.\ 64--65]{Henry},
the unknown function
$w_{(\ii y + z^{*})}^{(\ii\omega + \omega^{*})} (x,\xi,t)$
possesses the same $C^1$-smoothness property, for every fixed
$t\in \Delta_{\vartheta'}$.
Next, we apply the Cauchy\--Riemann operators
\begin{align*}
  \frac{\partial}{\partial\bar{z}^{*}}\eqdef
  \frac{1}{2}
    \left( \frac{\partial}{\partial x^{*}}
    + \ii\,\frac{\partial}{\partial y^{*}}
    \right) \quad\mbox{ and }\quad
  \frac{\partial}{\partial\bar{\omega}^{*}}\eqdef
  \frac{1}{2}
    \left( \frac{\partial}{\partial\alpha^{*}}
    + \ii\,\frac{\partial}{\partial\beta^{*}}
    \right)
\end{align*}
to both sides of eq.~\eqref{e:Cauchy:u-u_0}
(differentiation with respect to parameters),
thus concluding that both derivatives,
\begin{equation*}
  \frac{\partial}{\partial\bar{z}^{*}}
    w_{(\ii y + z^{*})}^{(\ii\omega + \omega^{*})} (x,\xi,t)
      \quad\mbox{ and }\quad
  \frac{\partial}{\partial\bar{\omega}^{*}}
    w_{(\ii y + z^{*})}^{(\ii\omega + \omega^{*})} (x,\xi,t) \,,
\end{equation*}
are the (unique) weak solutions of
the initial value problem \eqref{e:Cauchy:u-u_0}
with the zero initial data.
Thus, both derivatives must vanish identically for all
$(z^{*},\omega^{*})\in \CC$
with $|z^{*}|$ and $|\zeta^{*}|$ small enough.
Consequently, the difference
\begin{math}
  \tilde{u}(\,\cdot\,,\,\cdot\, t) - u_0\colon
  \mathfrak{V}^{(r')}\to \CC
\end{math}
is holomorphic, and so is the function
\begin{math}
  \tilde{u}(\,\cdot\,,\,\cdot\, t) \colon
  \mathfrak{V}^{(r')}\to \CC \,,
\end{math}
as claimed.
{\sc D.~Henry} provides an alternative proof of analyticity in his
\cite[Corollary 3.4.5, p.~65]{Henry}
that employs an analytic implicit function theorem via
Lemmas 3.4.2 and 3.4.3 in \cite[pp.\ 63--64]{Henry}.

To complete our proof of Proposition~\ref{prop-Main},
we apply the classical Hartogs's theorem on separate analyticity
(see, e.g.,
 {\sc S.~G.\ Krantz} \cite[Theorem 1.2.5, p.~32]{Krantz}
 and remarks around)
to conclude that the function
\begin{math}
  \tilde{u}\colon \mathfrak{V}^{(r)}\times \Delta_{\vartheta''}
  \to \CC \,,
\end{math}
defined by the formula in eq.~\eqref{def:u(x,xi,t)},
is holomorphic not only separately in the variables
$(z,\zeta)\in \mathfrak{V}^{(r')}$ and $t\in \Delta_{\vartheta''}$,
but also jointly in $(z,\zeta,t)$ in its entire domain.
\qed
\par\vskip 10pt


\section{$L^2$-bounds in the complex domain}
\label{s:L^2-bound}

In order to give a plausible lower estimate on
the space\--time domain of holomorphy
(i.e., the domain of complex analyticity)
of a weak solution $u$ to
the homogeneous initial value problem \eqref{e:Cauchy}
with $f\equiv 0$,
we introduce a few more subsets of $\CC^2\times \CC$
(cf.\ {\sc P.\ Tak\'a\v{c}} et al.\ \cite[p.~428]{TakacBoller}
 or   {\sc P.\ Tak\'a\v{c}} \cite[pp.\ 58--59]{Takac-12}):

The two constants $\kappa_0, \nu_0\in (0,\infty)$
used below will be specified later
(in the proof of Theorem~\ref{thm-Main});
$0\leq \alpha < \infty$ is an arbitrary number.
First, we recall the definitions of the complex sets
$\mathfrak{V}^{(\kappa_0\alpha)}\subset \CC^2$,
$\Sigma^{(\alpha)} (\nu_0)\subset \CC$, and
$\Gamma^{(T')}_T (\kappa_0,\nu_0) \subset \CC^2\times \CC$
given in Section~\ref{s:prelim}, eqs.\
\eqref{e:Pi}, \eqref{e:Sigma}, and \eqref{e:Gamma}, respectively.

Let us introduce the function
$\chi(s)\eqdef \min\{ s,\, 1\}$ for
$s\in \RR_+\eqdef [0,\infty)$; hence, it's derivative is given by
$\chi'(s) = 1$ for $0\leq s\leq 1$ and
$\chi'(s) = 0$ for $1 < s < \infty$.
Since the $x$\--section of
$\Gamma^{(T')}_T (\kappa_0,\nu_0)$ is independent from $x\in \RR$,
if $\kappa_0 T' < \pi / 2$, setting
\begin{equation}
\label{e:Gamma_x}
\begin{split}
  \hat{\Gamma}^{(T')}_T (\kappa_0,\nu_0)
  \eqdef
    \Bigl\{
&   (y,\zeta,t) = (y, \xi + \ii\eta, \alpha + \ii\tau)
    \in \RR\times \CC\times \CC\colon
\\
&   0 < \alpha < T \;\mbox{ together with }\;
    |y| < \kappa_0 T'\chi\genfrac{(}{)}{}1{\alpha}{T'} \,,\ \xi > 0 \,,
\\
&   |\arctan (\eta / \xi)|
    < \kappa_0 T'\chi\genfrac{(}{)}{}1{\alpha}{T'} \,,
  \;\mbox{ and\/ }\; \nu_0 |\tau| < T'\chi\genfrac{(}{)}{}1{\alpha}{T'}
    \Bigr\} \,,
\end{split}
\end{equation}
we may identify
\begin{math}
    \Gamma^{(T')}_T (\kappa_0,\nu_0)
  \simeq \RR\times \hat{\Gamma}^{(T')}_T (\kappa_0,\nu_0) .
\end{math}

The most important part of the proof of Theorem~\ref{thm-Main}
is the a~priori estimate in \eqref{e:u_Hardy^2}.
It is proved in the following proposition.
An example of a holomorphic extension
\begin{math}
  \tilde{u}\colon
  \mathfrak{V}^{(r)}\times \Delta_{\vartheta'} \to \CC
\end{math}
to a complex domain containing
$\Gamma^{(T')}_T (\kappa_0,\nu_0) \subset \CC^3$
is given in Proposition~\ref{prop-Main},
provided $\kappa_0$, $\nu_0^{-1}$, and $T'\in (0,T]$ are small enough.

\begin{proposition}\label{prop-Main-est}
Let\/
$\rho$, $\sigma$, $\theta$, $q_r$, and\/ $\gamma$
be given constants in $\RR$,
$\rho\in (-1,1)$, $\sigma > 0$, $\theta > 0$, and\/ $\gamma > 0$.
Assume that\/ $\beta$, $\gamma$, $\kappa$, and\/ $\mu$
are chosen as specified in {\rm Proposition~\ref{prop-Lions}}.
Then, given any numbers $r\in (0,\infty)$ and\/
$\vartheta'\in (0, \pi / 2)$,
the constants $\kappa_0, \nu_0^{-1}\in (0,\infty)$ and\/
$T'\in (0,T]$ can be chosen sufficiently small, such that
\begin{equation*}
  \Gamma^{(T')}_T (\kappa_0,\nu_0) \subset
  \mathfrak{V}^{(r)}\times \Delta_{\vartheta'}
\end{equation*}
and there exist some constants $C_0, c_0\in \RR_+$
with the following property:

If\/
$u_0\colon \mathfrak{V}^{(r)}\to \CC$
is a holomorphic function that satisfies
the bound \eqref{est:bound-u_0} in {\rm Proposition~\ref{prop-Main}}
and if\/
\begin{math}
  \tilde{u}\colon
  \mathfrak{V}^{(r)}\times \Delta_{\vartheta'} \to \CC
\end{math}
is the holomorphic extension of the (unique) weak solution
\begin{equation*}
  u\in C([0,T]\to H)\cap L^2((0,T)\to V)
\end{equation*} 
of the homogeneous initial value problem \eqref{e:Cauchy}
(with $f\equiv 0$ and this $u_0$)
that has been obtained in {\rm Proposition~\ref{prop-Main}},
then the estimate in \eqref{e:u_Hardy^2} holds
with the constants $C_0 = 1$ and $c_0 = c_2'\in \RR_+$ from
{\rm Proposition~\ref{prop-coerc-R}\/},
for every\/ $\alpha\in (0,T]$ and for all\/
$y, \omega, \tau\in \RR$ satisfying\/
\eqref{dom:u_Hardy^2}, depending on $\alpha$.
depending on $\alpha$.
\end{proposition}
\par\vskip 10pt

Before giving the {\em proof\/} of this proposition,
we first observe that the holomorphic extension
$\tilde{u}(z,\zeta,t)$ must be unique,
by uniqueness of the holomorphic extension in each of the variables
$z, \zeta, t\in \CC$.
Consequently, the remarks following the statement of
Proposition~\ref{prop-Main} apply also in the setting of our
Proposition~\ref{prop-Main-est}.
The holomorphic extension
\begin{math}
  \tilde{u}\colon \Gamma^{(T')}_T (\kappa_0,\nu_0) \to \CC
\end{math}
of a weak solution
\begin{equation*}
  u\in C([0,T]\to H)\cap L^2((0,T)\to V)
\end{equation*}
of the homogeneous initial value problem \eqref{e:Cauchy}
must satisfy the following initial value problem
with complex partial derivatives; cf.~\eqref{e:Cauchy-compl}:
\begin{equation}
\label{e_t:Cauchy-compl}
\left\{
\begin{alignedat}{2}
  \frac{\partial\tilde{u}}{\partial t}
  + ( \tilde{\mathcal{A}}\tilde{u} ) (z,\zeta,t) &= 0
  &&\quad\mbox{ in }\,
    \Gamma^{(T')}_T (\kappa_0,\nu_0) \,;
\\
  \tilde{u}(z,\zeta,0) &= u_0(z,\zeta)
  &&\quad\mbox{ for }\, (z,\zeta) = (x,\xi)\in \HH \,,
\end{alignedat}
\right.
\end{equation}
where the complex partial differential operator
$\tilde{\mathcal{A}}$ is given by eq.~\eqref{eq_C:Heston-oper} and
$\tilde{u}\in \mathcal{H}^2(\mathfrak{V}^{(r)})$.

\par\vskip 10pt
{\it Proof of\/} {\bf Proposition~\ref{prop-Main-est}.}$\;$
In order to establish the estimate in \eqref{e:u_Hardy^2},
we need to control the behavior of the holomorphic extension
$\tilde{u}(z,\zeta,t)$ of the solution $u(x,\xi,t)$
at every point
\begin{equation*}
  (z,\zeta,t) =
  \left( x + \ii y ,\, \xi (1 + \ii\omega) ,\, \alpha + \ii\tau
  \right) \in \Gamma^{(T')}_T (\kappa_0,\nu_0)
\end{equation*}
by the initial condition
$u_0\colon \HH\to \CC$ defined only at points
\begin{math}
  (x,\xi,0)\in \HH\times \{ 0\} = \RR\times (0,\infty)\times \{ 0\} .
\end{math}
Given any such two points, $(x,\xi,0)$ and $(z,\zeta,t)$,
we connect them by the following piecewise linear path parametrized by
the real time $s\in [ 0,\, \RE t ]$, i.e., by $0\leq s\leq \alpha$:

Given any point
\begin{equation*}
  (z,\zeta,t) =
  \left( x + \ii y ,\, \xi (1 + \ii\omega) ,\, \alpha + \ii\tau
  \right) \in \Gamma^{(T')}_T (\kappa_0,\nu_0) \,,
\end{equation*}
we set
\begin{equation*}
  y_0 = \frac{T'}{ \min\{ \alpha ,\, T'\} }\, y \,,\quad
  \omega_0 = \tan
    \left( \frac{T'}{ \min\{ \alpha ,\, T'\} }\, \arctan\omega
    \right) \,, \quad\mbox{ and }\quad
  \phi = \frac{\tau}{\alpha} \,.
\end{equation*}
Thus, conditions \eqref{dom:u_Hardy^2} are equivalent with
\begin{equation}
\label{dom_0:u_Hardy^2}
    \max\{ |y_0| ,\, |\arctan\omega_0| \}
  < \kappa_0 T'
    \quad\mbox{ and }\quad
  |\phi| < \nu_0^{-1} \,.
\end{equation}
Fixing $(y_0, \omega_0, \phi)\in \RR^3$
as in \eqref{dom_0:u_Hardy^2} above, we recall
$\chi(s)\eqdef \min\{ s,\, 1\}$ for
$s\in \RR_+\eqdef [0,\infty)$ and define the path
\begin{align}
\nonumber
& \varsigma\equiv \varsigma_{x,\xi}\colon
  [0,T]\to \{ (x,\xi,0) \}\cup \Gamma^{(T')}_T (\kappa_0,\nu_0)
  \colon
\\
\label{e:z,zeta,t}
& s \,\longmapsto\,
  \Bigl( x + \ii\chi(s/T') y_0 ,\,
  \xi \left( 1 + \ii\chi(s/T')\omega_0 \right) ,\, (1 + \ii\phi) s
  \Bigr) \,.
\\
\nonumber
& {}= (x,\xi,s)
    + \ii\bigl( \chi(s/T') y_0 ,\, \chi(s/T')\omega_0 ,\, \phi s\bigr) \,.
\end{align}
The numbers $y, \omega, \phi\in \RR$ are related to $(z,\zeta,t)$ by
$\phi = \frac{\tau}{\alpha}$, $y = \IM z$, and
$\omega = \frac{\IM\zeta}{\RE\zeta}$.
For $s = 0$ and $s = \alpha = \RE t$ we get the points
$(x,\xi,0)$ and $(z,\zeta,t)$, respectively.

Next, we define the function
$v\colon \HH\times [0,T]\to \CC$ by the values of $\tilde{u}$
on the image of the path $\varsigma$,
\begin{equation}
\label{e:v(x,xi,s)}
  v(x,\xi,s)
  \eqdef \tilde{u}
    \Bigl(
    x + \ii\chi\genfrac{(}{)}{}1{s}{T'} y_0 ,\,
    \xi \left( 1 + \ii\chi\genfrac{(}{)}{}1{s}{T'} \omega_0 \right) ,\,
    (1 + \ii\phi) s
    \Bigr) \,,\quad (x,\xi,s)\in \HH\times [0,T] \,.
\end{equation}
We calculate
\begin{align}
\label{e:dv/ds}
    \frac{\partial v}{\partial s} (x,\xi,s)
& = (1 + \ii\phi)\,
    \frac{\partial\tilde{u}}{\partial t}
  + \frac{\ii}{T'}\cdot \chi' \genfrac{(}{)}{}1{s}{T'}
    \left(
      \frac{\partial\tilde{u}}{\partial z}\, y_0
    + \frac{\partial\tilde{u}}{\partial\zeta}\, \xi\omega_0
    \right) \,,
\\
\label{e:dv/dx}
    \frac{\partial v}{\partial x} (x,\xi,s)
& = \frac{\partial\tilde{u}}{\partial z} \,,
\\
\label{e:dv/dxi}
    \frac{\partial v}{\partial\xi} (x,\xi,s)
& = \left( 1 + \ii\chi\genfrac{(}{)}{}1{s}{T'} \omega_0 \right)
    \frac{\partial\tilde{u}}{\partial\zeta} \,.
\end{align}

We prefer to use the complex form \eqref{eq_C:Heston-oper}
of the (time\--independent) Heston operator \eqref{eq:Heston-oper}.
Hence, according to the initial value problem \eqref{e_t:Cauchy-compl},
\begin{equation*}
  v\in C([0,T]\to H)\cap L^2((0,T)\to V)
\end{equation*}
is a weak solution of the following initial value problem,
\begin{equation}
\label{e_t:Cauchy-real}
\left\{
\begin{alignedat}{2}
  \frac{\partial v}{\partial s}
  + ( \hat{\mathcal{A}}(s) v ) (x,\xi,s) &= 0
  &&\quad\mbox{ in }\, \HH\times (0,T) \,;
\\
  v(x,\xi,0) &= u_0(x,\xi)
  &&\quad\mbox{ for }\, (x,\xi)\in \HH \,,
\end{alignedat}
\right.
\end{equation}
where the (time\--dependent) partial differential operator
$\hat{\mathcal{A}}(s)$ is given by
\begin{align*}
& ( \hat{\mathcal{A}}(s) v ) (x,\xi)\eqdef
  (1 + \ii\phi)\,
  ( \tilde{\mathcal{A}}\tilde{u} ) (z,\zeta)
  - \frac{\ii}{T'}\cdot \chi' \genfrac{(}{)}{}1{s}{T'}
    \left(
      \frac{\partial\tilde{u}}{\partial z}\, y_0
    + \frac{\partial\tilde{u}}{\partial\zeta}\, \xi\omega_0
    \right)
\\
&
\begin{alignedat}{2}
  = {}
  - \frac{1}{2}\, (1 + \ii\phi)\sigma\xi\cdot
& \left[
    \left( 1 + \ii\chi\genfrac{(}{)}{}1{s}{T'} \omega_0 \right)
    \frac{\partial^2 v}{\partial x^2}
  + 2\rho\,
    \frac{\partial^2 v}{\partial x\, \partial\xi}(x,\xi)
  \right.
\\
& \left.
{}
  + \left( 1 + \ii\chi\genfrac{(}{)}{}1{s}{T'} \omega_0 \right)^{-1}
    \frac{\partial^2 v}{\partial\xi^2}(x,\xi)
  \right]
\end{alignedat}
\\
& {}
  + (1 + \ii\phi)
    \left[ q_r + \genfrac{}{}{}1{1}{2}
       \left( 1 + \ii\chi\genfrac{(}{)}{}1{s}{T'} \omega_0 \right)
       \,\sigma\xi
    \right]
    \frac{\partial v}{\partial x}(x,\xi)
\\
& {}
  - (1 + \ii\phi)\kappa
    \left[ \theta_{\sigma}
    \left( 1 + \ii\chi\genfrac{(}{)}{}1{s}{T'} \omega_0 \right)^{-1}
  - \xi
    \right]
    \frac{\partial v}{\partial\xi}(x,\xi)
\\
& {}
  - \frac{\ii}{T'}\cdot \chi' \genfrac{(}{)}{}1{s}{T'}
    \left[
      y_0\, \frac{\partial v}{\partial x}
  + \left( 1 + \ii\chi\genfrac{(}{)}{}1{s}{T'} \omega_0 \right)^{-1}
      \xi\omega_0\, \frac{\partial v}{\partial\xi}
    \right]
\\
& {}
  = (1 + \ii\phi)\cdot (\mathcal{A}v) (x,\xi)
\\
& {} - \frac{\ii}{2}\, (1 + \ii\phi)\sigma\xi\cdot
    \chi\genfrac{(}{)}{}1{s}{T'} \omega_0
  \left[
    \frac{\partial^2 v}{\partial x^2}
  - \left( 1 + \ii\chi\genfrac{(}{)}{}1{s}{T'} \omega_0 \right)^{-1}
    \frac{\partial^2 v}{\partial\xi^2}
  \right]
\\
& {}
  + \frac{\ii}{2}\, (1 + \ii\phi)\cdot
    \chi\genfrac{(}{)}{}1{s}{T'} \omega_0
    \left[
    \sigma\xi\, \frac{\partial v}{\partial x}(x,\xi)
  + 2\kappa\theta_{\sigma}
    \left( 1 + \ii\chi\genfrac{(}{)}{}1{s}{T'} \omega_0 \right)^{-1}
    \frac{\partial v}{\partial\xi}(x,\xi)
    \right]
\\
& {}
  - \frac{\ii}{T'}\cdot \chi' \genfrac{(}{)}{}1{s}{T'}
    \left[
      y_0\, \frac{\partial v}{\partial x}
  + \left( 1 + \ii\chi\genfrac{(}{)}{}1{s}{T'} \omega_0 \right)^{-1}
      \xi\omega_0\, \frac{\partial v}{\partial\xi}
    \right]
\end{align*}
which yields the following formula,
\begin{align}
\nonumber
& ( \hat{\mathcal{A}}(s) v ) (x,\xi)
  {}= (1 + \ii\phi)\cdot (\mathcal{A}v) (x,\xi)
\\
\label{eq_R:Heston-oper}
& {}
  - \ii\cdot \frac{y_0}{T'}\cdot ( \mathcal{L}_1(s) v ) (x,\xi)
  - \ii\cdot \frac{\omega_0}{T'}\cdot ( \mathcal{L}_2(s) v ) (x,\xi)
\\
\nonumber
& {}
  + \frac{\ii}{2}\, (1 + \ii\phi)\sigma\, \omega_0
    \cdot ( \mathcal{L}_3(s) v ) (x,\xi)
  + \ii (1 + \ii\phi)\kappa\theta_{\sigma}\, \omega_0
    \cdot ( \mathcal{L}_4(s) v ) (x,\xi) \,,
\end{align}
where we have abbreviated
\begin{align}
\label{e:L_1}
& (\mathcal{L}_1(s) v)(x,\xi)
  \eqdef
    \chi'\genfrac{(}{)}{}1{s}{T'}
    \cdot \frac{\partial v}{\partial x}(x,\xi) \,,
\\
\label{e:L_2}
& (\mathcal{L}_2(s) v)(x,\xi)
  \eqdef
    \chi'\genfrac{(}{)}{}1{s}{T'}
    \left( 1 + \ii\chi\genfrac{(}{)}{}1{s}{T'} \omega_0 \right)^{-1}
      \xi\, \frac{\partial v}{\partial\xi}(x,\xi) \,,
\\
\label{e:L_3}
&
\begin{alignedat}{2}
  (\mathcal{L}_3(s) v)(x,\xi)
  \eqdef
& {} - \chi\genfrac{(}{)}{}1{s}{T'} \xi
  \left[
    \frac{\partial^2 v}{\partial x^2}
  - \left( 1 + \ii\chi\genfrac{(}{)}{}1{s}{T'} \omega_0 \right)^{-1}
    \frac{\partial^2 v}{\partial\xi^2}
  - \frac{\partial v}{\partial x}
  \right] \,, \quad\mbox{ and }\quad
\end{alignedat}
\\
\label{e:L_4}
& (\mathcal{L}_4(s) v)(x,\xi)
  \eqdef
    \chi\genfrac{(}{)}{}1{s}{T'}
    \left( 1 + \ii\chi\genfrac{(}{)}{}1{s}{T'} \omega_0 \right)^{-1}
    \frac{\partial v}{\partial\xi}
  \quad\mbox{ for }\, (x,\xi)\in \HH \,.
\end{align}

We insert eq.~\eqref{eq_R:Heston-oper} into \eqref{e_t:Cauchy-real},
thus arriving at
\begin{align}
\nonumber
&   \frac{\partial v}{\partial s} (x,\xi,s) =
{}- (1 + \ii\phi)\cdot (\mathcal{A}v) (x,\xi)
\\
\label{eq:dv/ds}
& {}
  + \ii\cdot \frac{y_0}{T'}\cdot ( \mathcal{L}_1(s) v ) (x,\xi)
  + \ii\cdot \frac{\omega_0}{T'}\cdot ( \mathcal{L}_2(s) v ) (x,\xi)
\\
\nonumber
& {}
  - \frac{\ii}{2}\, (1 + \ii\phi)\sigma\, \omega_0
    \cdot ( \mathcal{L}_3(s) v ) (x,\xi)
  - \ii (1 + \ii\phi)\kappa\theta_{\sigma}\, \omega_0
    \cdot ( \mathcal{L}_4(s) v ) (x,\xi)
\end{align}
for $(x,\xi,s)\in \HH\times (0,T)$.

In Propositions \ref{prop-bound-R} and \ref{prop-coerc-R}
above we have verified the boundedness and coercivity hypotheses
for the linear operator $\mathcal{A}\colon V\to V'$
defined by sesquilinear form in eq.~\eqref{e:Heston-diss:u=w}.
Estimates analogous to those used in the proof of
Proposition~\ref{prop-bound-R} show that all linear operators
$\mathcal{L}_j(s)\colon V\to V'$; $j=1,2,3,4$,
are uniformly bounded for $s\in [0,T]$ and $\omega_0\in \RR$, i.e.,
there is a constant $L\in (0,\infty)$ such that
\begin{equation}
\label{e:bound_L_j}
  \left\vert
  \left( \mathcal{L}_j(s) v, w\right)_H
  \right\vert \leq L\cdot \| v\|_V\, \| w\|_V
  \quad\mbox{ holds for all }\, v,w\in V
\end{equation}
and for all $s\in [0,T]$ and all $\omega_0\in \RR$; $j=1,2,3,4$.
Here, we have used the definition of
$\chi(s) = \min\{ s,\, 1\}$ and
\begin{math}
  \left| 1 + \ii\chi\genfrac{(}{)}{}1{s}{T'} \omega_0 \right|
  \geq 1 .
\end{math}

In order to obtain the upper bound \eqref{e:u_Hardy^2}
for the integral on the left\--hand side,
\begin{align*}
& \int_0^{\infty} \int_{-\infty}^{+\infty}
  \left| \tilde{u}
  \left( x + \ii y ,\, \xi (1 + \ii\omega) ,\, \alpha + \ii\tau
  \right)
  \right|^2\cdot \mathfrak{w}(x,\xi) \,\mathrm{d}x \,\mathrm{d}\xi
\\
& {}
  = \int_0^{\infty} \int_{-\infty}^{+\infty} |v(x,\xi,s)|^2
    \mathfrak{w}(x,\xi) \,\mathrm{d}x \,\mathrm{d}\xi
  = \| v(\,\cdot\,,\,\cdot\,,s) \|_H^2 \,,
\end{align*}
cf.\ eq.~\eqref{e:v(x,xi,s)},
we first take the time derivative of the second integral above,
then apply eq.~\eqref{eq:dv/ds}:
\begin{align*}
&   \frac{\mathrm{d}}{\mathrm{d}s}\, \| v(\,\cdot\,,\,\cdot\,,s) \|_H^2
  = \int_{\HH}
    \left( \frac{\partial v}{\partial s}\, \bar{v}
         + v\, \frac{\partial\bar{v}}{\partial s}
    \right)
    \mathfrak{w}(x,\xi) \,\mathrm{d}x \,\mathrm{d}\xi
\\
& = {}- \int_{\HH}
    \left( (\mathcal{A}v) (x,\xi)\, \bar{v}
         + v\, \overline{(\mathcal{A}v) (x,\xi)}
    \right) 
    \mathfrak{w}(x,\xi) \,\mathrm{d}x \,\mathrm{d}\xi
\\
& {}- \ii\phi \int_{\HH}
    \left( (\mathcal{A}v) (x,\xi)\, \bar{v}
         - v\, \overline{(\mathcal{A}v) (x,\xi)}
    \right) 
    \mathfrak{w}(x,\xi) \,\mathrm{d}x \,\mathrm{d}\xi
\\
& {}+ \ii\, \frac{y_0}{T'}\int_{\HH}
    \left( (\mathcal{L}_1(s) v) (x,\xi)\, \bar{v}
         - v\, \overline{(\mathcal{L}_1(s) v) (x,\xi)}
    \right) 
    \mathfrak{w}(x,\xi) \,\mathrm{d}x \,\mathrm{d}\xi
\\
& {}+ \ii\, \frac{\omega_0}{T'}\int_{\HH}
    \left( (\mathcal{L}_2(s) v) (x,\xi)\, \bar{v}
         - v\, \overline{(\mathcal{L}_2(s) v) (x,\xi)}
    \right) 
    \mathfrak{w}(x,\xi) \,\mathrm{d}x \,\mathrm{d}\xi
\\
& {}- \frac{\ii}{2}\, \sigma\omega_0 \int_{\HH}
    \left( (1 + \ii\phi) (\mathcal{L}_3(s) v) (x,\xi)\, \bar{v}
         - (1 - \ii\phi) v\, \overline{(\mathcal{L}_3(s) v) (x,\xi)}
    \right) 
    \mathfrak{w}(x,\xi) \,\mathrm{d}x \,\mathrm{d}\xi
\\
& {}- \ii\kappa \theta_{\sigma}\omega_0 \int_{\HH}
    \left( (1 + \ii\phi) (\mathcal{L}_4(s) v) (x,\xi)\, \bar{v}
         - (1 - \ii\phi) v\, \overline{(\mathcal{L}_4(s) v) (x,\xi)}
    \right) 
    \mathfrak{w}(x,\xi) \,\mathrm{d}x \,\mathrm{d}\xi \,.
\end{align*}
We estimate the integrals on the right\--hand side above as follows.
First, we take advantage of the coercivity of
$\mathcal{A}\colon V\to V'$ expressed in terms of
the {\em G\r{a}rding inequality\/} \eqref{ineq:Heston-diss:RE}.
Second, we employ the boundedness of $\mathcal{A}$, i.e.,
ineq.~\eqref{est:Heston-diss}.
Third, we employ the boundedness of $\mathcal{L}_j(s)$, i.e.,
ineq.~\eqref{e:bound_L_j}.
Consequently, we arrive at
\begin{align}
\nonumber
&   \frac{\mathrm{d}}{\mathrm{d}s}\,
    \| v(\,\cdot\,,\,\cdot\,,s) \|_H^2
  = \int_{\HH}
    \left( \frac{\partial v}{\partial s}\, \bar{v}
         + v\, \frac{\partial\bar{v}}{\partial s}
    \right)
    \mathfrak{w}(x,\xi) \,\mathrm{d}x \,\mathrm{d}\xi
\\
\label{est:u_Hardy^2}
& \leq {}- \sigma\, (1 - |\rho|)\cdot \| v\|_V^2 + c_2'\cdot \| v\|_H^2
\\
\nonumber
& {}+ 2 C |\phi|\, \| v\|_V^2
    + 2 L\, \frac{|y_0|}{T'}\, \| v\|_V^2
    + 2 L\, \frac{|\omega_0|}{T'}\, \| v\|_V^2
\\
\nonumber
& {}+ L\, |1 + \ii\phi|\, \sigma |\omega_0|\, \| v\|_V^2
    + 2 L\, |1 + \ii\phi|\, \kappa\theta_{\sigma}\, |\omega_0|\,
      \| v\|_V^2 \,.
\end{align}
To estimate the coefficients on the right\--hand side above,
we recall the conditions on $(y_0, \omega_0, \phi)\in \RR^3$
required in \eqref{dom_0:u_Hardy^2}.
In order to estimate the ratio $\omega_0 / T'$ in a simple way,
let us take the constants $\kappa_0\in (0,\infty)$ and $T'\in (0,T]$
small enough, such that $\kappa_0 T'\leq \pi / 4$.
The function $x\mapsto x^{-1}\tan x$ being
strictly monotone increasing on $(0,\infty)$,
with the limit equal to~$1$ as $x\to 0+$,
we employ condtition \eqref{dom_0:u_Hardy^2} to obtain
\begin{equation*}
    \frac{|\omega_0|}{T'}
  < \frac{\kappa_0}{\kappa_0 T'}\cdot \tan (\kappa_0 T')
  \leq \kappa_0\cdot \frac{ \tan (\pi / 4) }{\pi / 4}
  = \frac{4\kappa_0}{\pi} < 2\kappa_0 \,.
\end{equation*}
Then ineq.~\eqref{est:u_Hardy^2} yields
\begin{align}
\nonumber
&   \frac{\mathrm{d}}{\mathrm{d}s}\,
    \| v(\,\cdot\,,\,\cdot\,,s) \|_H^2
  \leq {}- \sigma\, (1 - |\rho|)\cdot \| v\|_V^2
  + c_2'\cdot \| v\|_H^2
\\
\label{ineq:u_Hardy^2}
& {}+ \left( 2 C \nu_0^{-1} + 2 L\kappa_0 + 4 L\kappa_0
      \right) \| v\|_V^2
\\
\nonumber
& {}+ \left(
      L (1 + \nu_0^{-1})\sigma\cdot 2\kappa_0 T'
    + 2 L (1 + \nu_0^{-1})\kappa \theta_{\sigma}\cdot 2\kappa_0 T'
      \right) \| v\|_V^2
\\
\nonumber
& = {}- \sigma\, (1 - |\rho|)\cdot \| v\|_V^2
  + c_2'\cdot \| v\|_H^2 + \tilde{C}\, \| v\|_V^2 \,,
\end{align}
where $\tilde{C}\in (0,\infty)$ is a constant,
\begin{align*}
  \tilde{C}
& \eqdef
    \left( 2 C \nu_0^{-1} + 2 L\kappa_0 + 4 L\kappa_0 \right)
\\
& {}
  + \left(
      L (1 + \nu_0^{-1})\sigma\cdot 2\kappa_0 T'
    + 2 L (1 + \nu_0^{-1})\kappa \theta_{\sigma}\cdot 2\kappa_0 T'
      \right)
\\
& {}
  = 2 C \nu_0^{-1} + 6 L\kappa_0
  + 2 L (1 + \nu_0^{-1}) (\sigma + 2\kappa \theta_{\sigma})
    \cdot \kappa_0 T' \,.
\end{align*}
Here, the constants
$\kappa_0, \nu_0^{-1}\in (0,\infty)$ and
$T'\in (0,T]$ can be chosen sufficiently small, such that
\begin{equation*}
  \Gamma^{(T')}_T (\kappa_0,\nu_0) \subset
  \mathfrak{V}^{(r)}\times \Delta_{\vartheta'}
\end{equation*}
holds together with
\begin{math}
  0 < \tilde{C}\leq \sigma\, (1 - |\rho|) .
\end{math}

Then ineq.~\eqref{ineq:u_Hardy^2} yields
\begin{equation*}
    \frac{\mathrm{d}}{\mathrm{d}s}\,
    \| v(\,\cdot\,,\,\cdot\,,s) \|_H^2
  \leq c_2'\cdot \| v\|_H^2
  \quad\mbox{ for }\, s\in (0,T) \,.
\end{equation*}
The desired inequality \eqref{e:u_Hardy^2}
now follows by taking
$C_0 = 1$, $c_0 = c_2'$, and
$s = \alpha$.

The proof of Proposition~\ref{prop-Main-est} is complete.
\qed
\par\vskip 10pt


\section{End of the proof of the main result}
\label{s:End}

In this section we finally finish the proof of Theorem~\ref{thm-Main}.
We will make use of the holomorphic approximation and the a~priori
estimates established in the previous two sections,
Sections \ref{s:Heston-compl} and~\ref{s:L^2-bound}.

For a given function
$u_0\in H = L^2(\mathbb{H};\mathfrak{w})$,
a sequence of entire (holomorphic) functions
\begin{equation*}
  \tilde{u}_{0,n}\colon \CC^2\to \CC \,;\quad n=1,2,2,\dots \,,
\end{equation*}
is constructed in Appendix~\ref{s:density} ({\S}~\ref{ss:Galerkin}),
whose restrictions to the complex domain
$\mathfrak{X}^{(r)} \times \Delta_{\vartheta_v}$
belong to
$H^2( \mathfrak{X}^{(r)} \times \Delta_{\vartheta_v} )$
and satify
\begin{equation*}
  \| \tilde{u}_{0,n}\vert_{\HH} - u_0\|_H \;\longrightarrow\; 0
    \quad\mbox{ as }\, n\to \infty \,;
\end{equation*}
cf.\ {\S}~\ref{ss:Galerkin},
properties {\rm (i)}, {\rm (ii)}, and {\rm (iii)}.
In paragraph {\S}\ref{ss:Cauchy-compl},
for every fixed $n = 1,2,3,\dots$, we have used the function
$\tilde{u}_{0,n}$ as the initial data for the initial value problem
\eqref{e:Cauchy-compl},
\begin{equation}
\label{e:Cauchy_u_n}
\left\{
\begin{alignedat}{2}
    \frac{\partial\tilde{u}_n}{\partial t}
  + \tilde{\mathcal{A}} \tilde{u}_n &= 0
    \quad\mbox{ for }\, (x,\xi,s)\in \HH\times (0,T) \,;
\\
  \tilde{u}_n
    \bigl( x + \ii y ,\, \xi (1 + \ii\omega) ,\, 0 \bigr)
& = \tilde{u}_{0,n}
    \bigl( x + \ii y ,\, \xi (1 + \ii\omega) \bigr)
    \quad\mbox{ for }\, (x,\xi)\in \HH \,.
\end{alignedat}
\right.
\end{equation}
Recall that $\tilde{\mathcal{A}}$ stands for
the natural complexification of the Heston operator $\mathcal{A}$
defined in eq.~\eqref{eq_C:Heston-oper}.
More precisely, this initial value problem has been solved by
general theory of holomorphic semigroups for fixed values of
$y,\omega\in \RR$ such that $|y| < r$ and
$|\arctan\omega| < \vartheta_v$.
In paragraph {\S}\ref{ss:Anal-compl_t}
we have proved that the unique weak solution
\begin{equation*}
  t \,\longmapsto\,
  \left[ (x,\xi)\mapsto
  \tilde{u}_n\bigl( x + \ii y ,\, \xi (1 + \ii\omega) ,\, t\bigr)
  \right] \colon [0,T]\to H
\end{equation*}
to problem \eqref{e:Cauchy_u_n} possesses a holomorphic extension
with respect to time $t$ to an angle $\Delta_{\vartheta_t}$,
for some $\vartheta_t\in (0, \pi / 2)$.
Furthermore, in paragraph {\S}\ref{ss:Cauchy-compl}
(Proposition~\ref{prop-Main})
we have proved that, for every $t\in \Delta_{\vartheta_t}$,
the solution
\begin{math}
  \tilde{u}_n (\,\cdot\,, \,\cdot\,,t)
  \colon \mathfrak{X}^{(r)} \times \Delta_{\vartheta_v}
  \,\longrightarrow\, \CC
\end{math}
is a holomorphic function that belongs to
$H^2( \mathfrak{X}^{(r)} \times \Delta_{\vartheta_v} )$.
Consequently, the function
\begin{math}
  \tilde{u}_n\colon
  \mathfrak{X}^{(r)} \times \Delta_{\vartheta_v}
  \times \Delta_{\vartheta_t} \,\longrightarrow\, \CC
\end{math}
is holomorphic in all its variables.

Now let us recall the time\--dependent path $\varsigma$
from \eqref{e:z,zeta,t},
\begin{align*}
& \varsigma\equiv \varsigma_{x,\xi}\colon
  [0,T]\to \{ (x,\xi,0) \}\cup \Gamma^{(T')}_T (\kappa_0,\nu_0)
  \colon
\\
& s \,\longmapsto\,
  \Bigl( x + \ii\chi(s/T') y_0 ,\,
  \xi \left( 1 + \ii\chi(s/T')\omega_0 \right) ,\, (1 + \ii\phi) s
  \Bigr) \,.
\\
& {}= (x,\xi,s)
    + \ii\bigl( \chi(s/T') y_0 ,\, \chi(s/T')\omega_0 ,\, \phi s\bigr) \,,
\end{align*}
where the numbers $y_0, \omega_0, \phi\in \RR$ obey conditions
\eqref{dom_0:u_Hardy^2},
\begin{equation*}
    \max\{ |y_0| ,\, |\arctan\omega_0| \}
  < \kappa_0 T'
    \quad\mbox{ and }\quad
  |\phi| < \nu_0^{-1} \,,
\end{equation*}
with some constants $\kappa_0, \nu_0^{-1}\in (0,\infty)$ and
$T'\in (0,T]$ small enough, such that also
\begin{equation*}
  \kappa_0 T'\leq \min\{ r,\, \vartheta_v\}
    \quad\mbox{ and }\quad 
  \nu_0^{-1}\leq \tan\vartheta_t \,.
\end{equation*}
Here,
$0 < \vartheta_v, \vartheta_t < \pi / 2$
are some given numbers.
In the previous section (Section~\ref{s:L^2-bound}),
Proposition~\ref{prop-Main-est},
we have shown that along this path,
$\varsigma\equiv \varsigma_{x,\xi}$,
whose value at each $s\in [0,T]$ is viewed as
a function of the pair $(x,\xi)\in \HH$,
the $H$-norm of the function
\begin{math}
  (x,\xi) \,\longmapsto\, v_n(x,\xi,s)
  \colon \HH\times [0,T]\to \CC \,,
\end{math}
defined by \eqref{e:v(x,xi,s)},
\begin{align*}
  v_n(x,\xi,s)
  \eqdef \tilde{u}_n
    \Bigl(
    x + \ii\chi\genfrac{(}{)}{}1{s}{T'} y_0 ,\,
    \xi \left( 1 + \ii\chi\genfrac{(}{)}{}1{s}{T'} \omega_0 \right) ,\,
    (1 + \ii\phi) s
    \Bigr) \,,
\\
  \quad (x,\xi,s)\in \HH\times [0,T] \,,
\end{align*}
is uniformly bounded with the bound depending solely on the norm
\begin{math}
  \| \tilde{u}_{0,n}\vert_{\HH} \|_H ,
\end{math}
the time interval length $T>0$, and the constant $c_2'> 0$
in inequality \eqref{ineq:Heston-diss:RE}.

Next, we take advantage of the fact that we treat
homogeneous {\em linear\/} parabolic problems,
\eqref{e:Cauchy} (with $f\equiv 0$) in the real domain
$\HH\times (0,T)$, and its natural complexification
\eqref{e:Cauchy-compl} in the complex domain
$\mathfrak{V}^{(r')}\times \Delta_{\vartheta'}$.
Consequently, given any indices $m,n\in \NN$, the difference
\begin{math}
  \tilde{u}_n - \tilde{u}_m\colon
  \mathfrak{V}^{(r')}\times \Delta_{\vartheta'} \to \CC
\end{math}
is a holomorphic function that obeys the parabolic equation
in problem~\eqref{e:Cauchy-compl}.
Hence, we may apply our crucial a~priori estimate
\eqref{e:u_Hardy^2} in Proposition~\ref{prop-Main-est}
to the difference $\tilde{u}_n - \tilde{u}_m$, thus obtaining
\begin{alignat}{2}
\nonumber
  \int_0^{\infty} \int_{-\infty}^{+\infty}
& \left| \tilde{u}_n
  \left( x + \ii y ,\, \xi (1 + \ii\omega) ,\, \alpha + \ii\tau
  \right)
  \right.
\\
\label{e:(u_n-u_m)_Hardy^2}
& \left.
  {}   - \tilde{u}_m
  \left( x + \ii y ,\, \xi (1 + \ii\omega) ,\, \alpha + \ii\tau
  \right)
  \right|^2\cdot \mathfrak{w}(x,\xi) \,\mathrm{d}x \,\mathrm{d}\xi
\\
\nonumber
& \leq \ee^{ c_2'\alpha }\cdot
  \int_0^{\infty} \int_{-\infty}^{+\infty}
    | \tilde{u}_n(x,\xi,0) - \tilde{u}_m(x,\xi,0) |^2
  \cdot \mathfrak{w}(x,\xi) \,\mathrm{d}x \,\mathrm{d}\xi
\\
\nonumber
& {}
    = \ee^{ c_2'\alpha }\cdot \| u_{0,n} - u_{0,m}\|_H^2
\end{alignat}
for every $\alpha\in (0,T]$ and for all
$y, \omega, \tau\in \RR$ satisfying conditions \eqref{dom:u_Hardy^2},
\begin{equation*}
    \max\{ |y| ,\, |\arctan\omega| \}
  < \kappa_0\cdot \min\{ \alpha ,\, T'\}
    \quad\mbox{ and }\quad
  \nu_0 |\tau| < \alpha \,,
\end{equation*}
depending on $\alpha$.

It follows from
$\tilde{u}_{0,n}\vert_{\HH} \to u_0$ in $H$ as $n\to \infty$, that
$\left\{ \tilde{u}_{0,n}\vert_{\HH} \right\}_{n=1}^{\infty}$
is a Cauchy sequence in~$H$.
By ineq.~\eqref{e:(u_n-u_m)_Hardy^2}, also the functions
\begin{align}
\label{e:w(x,xi)}
  w_n(x,\xi)
  \eqdef \tilde{u}_n
    \Bigl( x + \ii y ,\, \xi (1 + \ii\omega) ,\, \alpha + \ii\tau
    \Bigr) \,,
  \quad (x,\xi)\in \HH \,,
\end{align}
form a Cauchy sequence $\{ w_n\}_{n=1}^{\infty}$ in $H$,
{\em uniformly\/} for all choices of $\alpha + \ii\tau\in \CC$ and
$y, \omega\in \RR$ satisfying $0 < \alpha\leq T$
and conditions \eqref{dom:u_Hardy^2}, that is to say, for
\begin{equation}
\label{cond:u_Hardy^2}
    \max\{ |y| ,\, |\arctan\omega| \}
  < \kappa_0\cdot \min\{ \alpha ,\, T'\}
    \quad\mbox{ and }\quad
  \nu_0 |\tau| < \alpha\leq T \,.
\end{equation}
Such numbers $\alpha + \ii\tau\in \CC$ and $y, \omega\in \RR$
being fixed, let
$w\eqdef \lim_{n\to \infty} w_n$ be the limit in $H$ of
this Cauchy sequence.
In analogy with eq.~\eqref{e:w(x,xi)}, we set
\begin{align}
\label{e:tilde_u=w(x,xi)}
  \tilde{u}
    \Bigl( x + \ii y ,\, \xi (1 + \ii\omega) ,\, \alpha + \ii\tau
    \Bigr)
  \eqdef w(x,\xi) \,,
  \quad (x,\xi)\in \HH \,.
\end{align}
Then
\begin{math}
  \tilde{u}\colon \Gamma^{(T')}_T (\kappa_0,\nu_0) \to \CC
\end{math}
is a complex\--valued, Lebesgue measurable function that satisfies
the following inequality,
by letting $m\to \infty$ in ineq.~\eqref{e:(u_n-u_m)_Hardy^2},
\begin{alignat}{2}
\nonumber
  \int_0^{\infty} \int_{-\infty}^{+\infty}
& \left| \tilde{u}_n
  \left( x + \ii y ,\, \xi (1 + \ii\omega) ,\, \alpha + \ii\tau
  \right)
  \right.
\\
\label{e:(u_n-u)_Hardy^2}
& \left.
  {}   - \tilde{u}
  \left( x + \ii y ,\, \xi (1 + \ii\omega) ,\, \alpha + \ii\tau
  \right)
  \right|^2\cdot \mathfrak{w}(x,\xi) \,\mathrm{d}x \,\mathrm{d}\xi
\\
\nonumber
& \leq \ee^{ c_2'\alpha }\cdot
  \int_0^{\infty} \int_{-\infty}^{+\infty}
    | \tilde{u}_n(x,\xi,0) - u_0(x,\xi) |^2
  \cdot \mathfrak{w}(x,\xi) \,\mathrm{d}x \,\mathrm{d}\xi
\\
\nonumber
& {}
    = \ee^{ c_2'\alpha }\cdot \| u_{0,n} - u_0\|_H^2
\end{alignat}
for all choices of $\alpha + \ii\tau\in \CC$ and
$y, \omega\in \RR$ satisfying conditions \eqref{cond:u_Hardy^2} above.

A trivial consequence of
\eqref{e:(u_n-u)_Hardy^2} and \eqref{cond:u_Hardy^2}
is that the sequence of functions
\begin{math}
  \tilde{u}_n\colon \Gamma^{(T')}_T (\kappa_0,\nu_0)\hfil\break \to \CC ;
\end{math}
$n=1,2,3,\dots$, converges in the complex domain 
$\Gamma^{(T')}_T (\kappa_0,\nu_0)$ to the function
\begin{math}
  \tilde{u}\colon\hfil\break \Gamma^{(T')}_T (\kappa_0,\nu_0) \to \CC
\end{math}
locally in the $L^2$-topology.
Since $\tilde{u}_n$ is holomorphic in
$\Gamma^{(T')}_T (\kappa_0,\nu_0)$,
it can be expressed by
the Cauchy integral formula for polydiscs
({\sc S.~G.\ Krantz} \cite{Krantz}, Theorem 1.2.2 (p.~24), or
 {\sc F.\ John} \cite{John}, Chapt.~3, Sect.\ 3(c), eq.\ (3.22c), p.~71).
From this formula we deduce by standard limiting arguments using
ineq.~\eqref{e:(u_n-u)_Hardy^2}
that also the limit function $\tilde{u}$ is expressed by
the same Cauchy integral formula for polydiscs.
It follows that also $\tilde{u}$ is holomorphic in
$\Gamma^{(T')}_T (\kappa_0,\nu_0)$, as desired.
Obviously, Proposition~\ref{prop-Main-est} guarantees that
$\tilde{u}$ satisfies ineq.~\eqref{e:u_Hardy^2}.

To derive the relation of $\tilde{u}$ to
problem~\eqref{e:Cauchy} (with $f\equiv 0$)
in the real domain $\HH\times (0,T)$, let us take
$y = \omega = \tau = 0$ in ineq.~\eqref{e:(u_n-u)_Hardy^2}.
Letting $n\to \infty$ we observe that the function
\begin{align}
\label{e:u(x,xi)}
  \hat{u}\colon (x,\xi,t) \longmapsto \tilde{u}(x,\xi,t)\colon
    \HH\times (0,T)\to \CC
\end{align}
is a weak solution to
the Cauchy problem~\eqref{e:Cauchy} (with $f\equiv 0$).
However, the initial value problem~\eqref{e:Cauchy}
(with $f\equiv 0$)
possesses a unique weak solution
\begin{equation*}
  u\in C([0,T]\to H)\cap L^2((0,T)\to V) \,,
\end{equation*}
by a pair of standard theorems for abstract parabolic problems due to
{\sc J.-L.\ Lions} \cite[Chapt.~IV]{Lions-61},
Th\'eor\`eme 1.1 ({\S}1, p.~46) and Th\'eor\`eme 2.1 ({\S}2, p.~52)
(for alternative proofs, see also e.g.\
{\sc L.~C.\ Evans} \cite[Chapt.~7, {\S}1.2(c)]{Evans-98},
Theorems 3 and~4, pp.\ 356--358,
{\sc J.-L.\ Lions} \cite[Chapt.~III, {\S}1.2]{Lions-71},
Theorem 1.2 (p.~102) and remarks thereafter (p.~103),
{\sc A.\ Friedman} \cite{Friedman-64},
Chapt.~10, Theorem~17, p.~316, or
{\sc H.\ Tanabe} \cite[Chapt.~5, {\S}5.5]{Tanabe},
Theorem 5.5.1, p.~150).

Hence, we have $\hat{u} = u$ in $\HH\times (0,T)$,
thus proving that 
\begin{math}
  \tilde{u}\colon \Gamma^{(T')}_T (\kappa_0,\nu_0) \to \CC
\end{math}
is a holomorphic extension of~$u$.

The proof of Theorem~\ref{thm-Main} is complete.
\qed
\par\vskip 10pt


\appendix

\section{Appendix: Trace, Sobolev's, and\\ Hardy's inequalities}
\label{s:Trace,Sobolev}

Our boundedness and coercivity results for the Heston operator
$\mathcal{A}\colon V\to V'$
make use of the following five lemmas:
Recall that $V = H^1(\mathbb{H};\mathfrak{w})$
and $\beta > 0$, $\gamma > 0$, and $\mu > 0$
are constants in the weight $\mathfrak{w}(x,\xi)$
which is defined in eq.~\eqref{def:w}.

\begin{lemma}\label{lem-Trace_pt}
{\rm (A pointwise trace inequality.)}$\;$
Let\/
$\beta > 0$, $\gamma > 0$, and\/ $\mu > 0$.
Then the following inequality holds for every function $u\in V$
and at almost every point\/ $x\in \RR$,
\begin{equation}
\label{ineq:trace_pt:v}
  \frac{\partial}{\partial\xi}
  \left( \xi^{\beta}\, \ee^{- \mu\xi}\, |u(x,\xi)|^2 \right)
  \leq \frac{1}{\mu}\, |u_{\xi}(x,\xi)|^2
    \cdot \xi^{\beta}\, \ee^{- \mu\xi}
  + \beta\, |u(x,\xi)|^2\cdot \xi^{\beta - 1}\, \ee^{- \mu\xi}
\end{equation}
for almost every\/ $\xi\in (0,\infty)$.

Furthermore, for a.e.\ $x\in \RR$ we have the limits\/
\begin{align}
\label{lim:trace_pt:v=0}
& \lim_{\xi\to 0+} \left( \xi^{\beta}\cdot |u(x,\xi)|^2 \right)
  = 0 \quad\mbox{ and }
\\
\label{lim:trace_pt:v=infty}
& \lim_{\xi\to \infty}
    \left( \xi^{\beta}\, \ee^{- \mu\xi}\cdot |u(x,\xi)|^2 \right)
  = 0 \,.
\end{align}
\end{lemma}

\proof
The following partial derivatives exist almost everywhere in $\HH$;
we first calculate
\begin{align*}
& \frac{\partial}{\partial\xi}
  \left( \xi^{\beta}\, \ee^{ - \mu\xi} |u(x,\xi)|^2
  \right)
\\
& = \left( u_{\xi}\, \bar{u} + u\, \bar{u}_{\xi} \right)
    \cdot \xi^{\beta}\, \ee^{- \mu\xi}
  + \beta\, |u(x,\xi)|^2\cdot \xi^{\beta - 1}\, \ee^{- \mu\xi}
  - \mu\, |u(x,\xi)|^2\cdot \xi^{\beta}\, \ee^{- \mu\xi} \,,
\end{align*}
then apply the Cauchy inequality
\begin{equation*}
    u_{\xi}\, \bar{u} + u\, \bar{u}_{\xi}
  = 2\cdot \RE ( u_{\xi}\, \bar{u} )
  \leq 2 |u_{\xi}|\cdot |u|
  \leq \mu^{-1}\, |u_{\xi}|^2 + \mu\, |u|^2
\end{equation*}
to estimate
\begin{equation*}
  \frac{\partial}{\partial\xi}
  \left( \xi^{\beta}\, \ee^{- \mu\xi}\, |u(x,\xi)|^2 \right)
  \leq \frac{1}{\mu}\, |u_{\xi}|^2
    \cdot \xi^{\beta}\, \ee^{- \mu\xi}
  + \beta\, |u(x,\xi)|^2\cdot \xi^{\beta - 1}\, \ee^{- \mu\xi} \,.
\end{equation*}
This proves ineq.~\eqref{ineq:trace_pt:v}.

Recall that $u\in V$.
Integrating the right\--hand side of the last inequality
with respect to the measure
$\ee^{ - \gamma |x| - \mu\xi } \,\mathrm{d}x \,\mathrm{d}\xi$
over $\HH = \RR\times (0,\infty)$
we infer that, for a.e.\ $x\in \RR$, both integrals below converge,
\begin{equation}
\label{int:trace:v}
  \int_0^{\infty} |u_{\xi}(x,\xi)|^2\cdot \xi^{\beta}\, \ee^{- \mu\xi}
    \,\mathrm{d}\xi < \infty
    \quad\mbox{ and }\quad
  \int_0^{\infty} |u(x,\xi)|^2\cdot \xi^{\beta - 1}\, \ee^{- \mu\xi}
    \,\mathrm{d}\xi < \infty \,.
\end{equation}
Let $x\in \RR$ be such a point.
The right\--hand side of ineq.~\eqref{ineq:trace_pt:v}
is integrable with respect to the Lebesgue measure
$\mathrm{d}\xi$ over $(0,\infty)$, and so is the positive part
$\phi^{+}(\xi) = \max\{ \phi(\xi),\, 0\}$
of the partial derivative
\begin{equation*}
  \xi \,\longmapsto\, \phi(\xi)\eqdef
    \frac{\partial}{\partial\xi}
  \left( \xi^{\beta}\, \ee^{- \mu\xi}\, |u(x,\xi)|^2 \right) \,.
\end{equation*}
Thus, the existence of the limit in \eqref{lim:trace_pt:v=0},
\begin{equation*}
  \lim_{\xi\to 0+} \left( \xi^{\beta}\cdot |u(x,\xi)|^2 \right)
  = L_0(x)
    \quad\mbox{ for a.e. }\, x\in \RR \,,
\end{equation*}
is deduced from
\begin{equation}
\label{lim-inf:trace_pt:v=0}
  L_0(x)\eqdef \liminf_{\xi\to 0+}
    \left( \xi^{\beta}\cdot |u(x,\xi)|^2 \right)
\end{equation}
and the following inequality, obtained by integrating
ineq.~\eqref{ineq:trace_pt:v}
and valid for all $0 < \xi' < \xi'' < \infty$,
\begin{align}
\label{est:trace:v}
&   (\xi'')^{\beta}\, \ee^{- \mu\xi''}\, |u(x,\xi'')|^2
  - (\xi')^{\beta}\,  \ee^{- \mu\xi'}    |u(x,\xi')|^2
  \eqdef
  \left[ \xi^{\beta}\, \ee^{- \mu\xi}\, |u(x,\xi)|^2
  \right]_{\xi = \xi'}^{\xi = \xi''}
\\
\nonumber
& \leq \frac{1}{\mu} \int_{\xi'}^{\xi''}
    |u_{\xi}(x,\xi)|^2\cdot \xi^{\beta}\, \ee^{- \mu\xi}
    \,\mathrm{d}\xi
  + \beta \int_{\xi'}^{\xi''}
    |u(x,\xi)|^2\cdot \xi^{\beta - 1}\, \ee^{- \mu\xi}
    \,\mathrm{d}\xi \,.
\end{align}
By similar reasoning, one derives
the existence of the limit in \eqref{lim:trace_pt:v=infty},
\begin{equation*}
  \lim_{\xi\to \infty}
    \left( \xi^{\beta}\, \ee^{- \mu\xi}\cdot |u(x,\xi)|^2 \right)
  = L_{\infty}(x)
    \quad\mbox{ for a.e. }\, x\in \RR \,,
\end{equation*}
from
\begin{equation}
\label{lim-inf:trace_pt:v=infty}
  L_{\infty}(x)\eqdef \liminf_{\xi\to \infty}
    \left( \xi^{\beta}\, \ee^{- \mu\xi}\cdot |u(x,\xi)|^2 \right) \,.
\end{equation}

Finally, both limits, $L_0(x)$ and $L_{\infty}(x)$,
are nonnegative and finite, by the integrability properties of
$u_{\xi}(x,\,\cdot\,)$ and $u(x,\,\cdot\,)$
stated in \eqref{int:trace:v}.
Moreover, the second integral in \eqref{int:trace:v} forces
$L_0(x) = L_{\infty}(x) = 0$, thanks to
\begin{math}
  \int_0^{\delta} \xi^{-1} \,\mathrm{d}\xi =
  \int_{1 / \delta}^{\infty} \xi^{-1} \,\mathrm{d}\xi = \infty
\end{math}
for any $\delta > 0$.
%
\qed
\par\vskip 10pt

Lemma~\ref{lem-Trace_pt} has the following global analogue
with a similar proof.

\begin{lemma}\label{lem-Trace}
{\rm (A trace inequality.)}$\;$
Let\/
$\beta > 0$, $\gamma > 0$, and\/ $\mu > 0$.
Then the following inequality holds for every function $u\in V$,
\begin{align}
\label{ineq:trace:v}
& \frac{\partial}{\partial\xi}
  \left( \xi^{\beta}\, \ee^{- \mu\xi} \int_{\RR} |u(x,\xi)|^2
    \cdot \ee^{- \gamma |x|} \,\mathrm{d}x
  \right)
\\
\nonumber
& \leq \frac{1}{\mu}\int_{\RR} |u_{\xi}(x,\xi)|^2
    \cdot \xi^{\beta}\, \ee^{ - \gamma |x| - \mu\xi } \,\mathrm{d}x
  + \beta\int_{\RR} |u(x,\xi)|^2
    \cdot \xi^{\beta - 1}\, \ee^{ - \gamma |x| - \mu\xi } \,\mathrm{d}x
\end{align}
for almost every\/ $\xi\in (0,\infty)$.

Furthermore, the limits in\/
\eqref{e:trace:v=0}  and~\eqref{e:trace:v=infty}
are valid.
\end{lemma}

\proof
We integrate both sides of ineq.~\eqref{ineq:trace_pt:v}
with respect to the measure
$\ee^{- \gamma |x|} \,\mathrm{d}x$ over $\RR$ to obtain
ineq.~\eqref{ineq:trace:v}.

Since $u\in V$, the right\--hand side of ineq.~\eqref{ineq:trace:v}
is integrable with respect to the Lebesgue measure
$\mathrm{d}\xi$ over $(0,\infty)$, and so is the positive part
$\phi^{+}(\xi) = \max\{ \phi(\xi),\, 0\}$
of the partial derivative
\begin{equation*}
  \xi \,\longmapsto\, \phi(\xi)\eqdef
    \frac{\partial}{\partial\xi}
  \left( \xi^{\beta}\, \ee^{- \mu\xi} \int_{\RR} |u(x,\xi)|^2
    \cdot \ee^{- \gamma |x|} \,\mathrm{d}x
  \right) \,.
\end{equation*}
Thus, the existence of the limit in \eqref{e:trace:v=0},
\begin{equation*}
  \lim_{\xi\to 0+}
    \left( \xi^{\beta}\cdot \int_{-\infty}^{+\infty}
      |u(x,\xi)|^2\cdot \ee^{- \gamma |x|} \,\mathrm{d}x
    \right) = L_0 \,,
\end{equation*}
is deduced from
\begin{equation}
\label{lim-inf:trace:v=0}
  L_0\eqdef \liminf_{\xi\to 0+}
    \left( \xi^{\beta}\cdot \int_{-\infty}^{+\infty}
      |u(x,\xi)|^2\cdot \ee^{- \gamma |x|} \,\mathrm{d}x
    \right)
\end{equation}
and the following inequality, obtained by integrating
ineq.~\eqref{ineq:trace:v}
and valid for all $0 < \xi' < \xi''$ $< \infty$,
cf.~\eqref{est:trace:v}:
\begin{align*}
&   (\xi'')^{\beta}\, \ee^{- \mu\xi''} \int_{\RR} |u(x,\xi'')|^2
    \cdot \ee^{- \gamma |x|} \,\mathrm{d}x
  - (\xi')^{\beta}\, \ee^{- \mu\xi'} \int_{\RR} |u(x,\xi')|^2
    \cdot \ee^{- \gamma |x|} \,\mathrm{d}x
\\
& \eqdef
  \left[ \xi^{\beta}\, \ee^{- \mu\xi} \int_{\RR} |u(x,\xi)|^2
    \cdot \ee^{- \gamma |x|} \,\mathrm{d}x
  \right]_{\xi = \xi'}^{\xi = \xi''}
\\
& \leq \frac{1}{\mu} \int_{\xi'}^{\xi''} \int_{\RR}
    |u_{\xi}(x,\xi)|^2\cdot \xi^{\beta}\, \ee^{ - \gamma |x| - \mu\xi }
    \,\mathrm{d}x \,\mathrm{d}\xi
\\
& {}
  + \beta \int_{\xi'}^{\xi''} \int_{\RR}
    |u(x,\xi)|^2\cdot \xi^{\beta - 1}\, \ee^{ - \gamma |x| - \mu\xi }
    \,\mathrm{d}x \,\mathrm{d}\xi \,.
\end{align*}
By similar reasoning, one derives
the existence of the limit in \eqref{e:trace:v=infty},
\begin{equation*}
  \lim_{\xi\to \infty}
    \left( \xi^{\beta}\, \ee^{- \mu\xi}\cdot \int_{-\infty}^{+\infty}
      |u(x,\xi)|^2\cdot \ee^{- \gamma |x|} \,\mathrm{d}x
    \right) = L_{\infty} \,,
\end{equation*}
from
\begin{equation}
\label{lim-inf:trace:v=infty}
  L_{\infty}\eqdef \liminf_{\xi\to \infty}
    \left( \xi^{\beta}\, \ee^{- \mu\xi}\cdot \int_{-\infty}^{+\infty}
      |u(x,\xi)|^2\cdot \ee^{- \gamma |x|} \,\mathrm{d}x
    \right) \,.
\end{equation}

Again, as in our proof of Lemma~\ref{lem-Trace_pt} above,
both limits, $L_0$ and $L_{\infty}$,
are nonnegative and finite, by the integrability properties of $u\in V$.
Moreover, $u\in H$ forces $L_0 = L_{\infty} = 0$, thanks to
\begin{math}
  \int_0^{\delta} \xi^{-1} \,\mathrm{d}\xi =
  \int_{1 / \delta}^{\infty} \xi^{-1} \,\mathrm{d}\xi = \infty
\end{math}
for any $\delta > 0$.
%
\qed
\par\vskip 10pt

Our second trace result, Lemma~\ref{lem-Trace_x} below,
is a simple analogue in the $x$-direction of Lemma~\ref{lem-Trace} above.
Its proof is analogous to that of Lemma~\ref{lem-Trace}
and is left to the reader; cf.\
{\sc A.\ Kufner} \cite{Kufner}.

\begin{lemma}\label{lem-Trace_x}
{\rm (Another trace inequality.)}$\;$
Let\/
$\beta > 0$, $\gamma > 0$, and\/ $\mu > 0$.
Then the limits in\/ \eqref{e:trace:x=+-infty}
hold for every function $u\in V$.
\end{lemma}

We take advantage of the trace results in
Lemmas \ref{lem-Trace_pt} and~\ref{lem-Trace}
to derive the following embedding lemma.

\begin{lemma}\label{lem-Sobolev}
{\rm (A Sobolev\--type inequality.)}$\;$
Let\/
$\beta > 0$, $\gamma > 0$, and\/ $\mu > 0$.
Then the following {\em Sobolev\--type inequality\/}
holds for every function $u\in V$,
\begin{align}
\label{ineq:Heston:Sobol}
&
\begin{aligned}
  \int_{\HH} |u(x,\xi)|^2
  \cdot \xi^{\beta}\, \ee^{ - \gamma |x| - \mu\xi }
    \,\mathrm{d}x \,\mathrm{d}\xi
& \leq
    \genfrac{(}{)}{}0{2}{\mu}^2
    \int_{\HH} |u_{\xi}(x,\xi)|^2
  \cdot \xi^{\beta}\, \ee^{ - \gamma |x| - \mu\xi }
    \,\mathrm{d}x \,\mathrm{d}\xi
\\
& {}
  + \frac{2\beta}{\mu}
    \int_{\HH} |u(x,\xi)|^2
  \cdot \xi^{\beta - 1}\, \ee^{ - \gamma |x| - \mu\xi }
    \,\mathrm{d}x \,\mathrm{d}\xi \,.
\end{aligned}
\end{align}
\end{lemma}

\proof
It suffices to verify the following inequality:
\begin{equation}
\label{ineq:Sobol:xi}
\begin{aligned}
  \int_0^{\infty}
  |u(\xi)|^2\cdot \xi^{\beta}\, \ee^{- \mu\xi} \,\mathrm{d}\xi
& \leq
    \genfrac{(}{)}{}0{2}{\mu}^2
    \int_0^{\infty} |u_{\xi}(\xi)|^2\cdot \xi^{\beta}\, \ee^{- \mu\xi}
    \,\mathrm{d}\xi
\\
& + \frac{2\beta}{\mu}
    \int_0^{\infty} |u(\xi)|^2\cdot \xi^{\beta - 1}\, \ee^{- \mu\xi}
    \,\mathrm{d}\xi
\end{aligned}
\end{equation}
holds for an arbitrary function
$u\in W^{1,2}_{\mathrm{loc}}(0,\infty)$ such that
\begin{align}
\label{e:u(xi):H^1}
& \int_0^{\infty} |u_{\xi}(\xi)|^2\cdot \xi^{\beta}\, \ee^{- \mu\xi}
    \,\mathrm{d}\xi < \infty
    \quad\mbox{ and }\quad
\\
\label{e:u(xi):trace:xi=0}
&   \lim_{\xi\to 0+} \left( \xi^{\beta}\cdot |u(\xi)|^2 \right)
  = \lim_{\xi\to \infty}
    \left( \xi^{\beta}\, \ee^{- \mu\xi}\cdot |u(\xi)|^2 \right)
  = 0 \,.
\end{align}
The boundary conditions in \eqref{e:u(xi):trace:xi=0}
are justified by Lemma~\ref{lem-Trace_pt}.

Indeed, we begin with the identities
\begin{align}
\nonumber
&   \mu\int_0^{\infty}
    |u(\xi)|^2\cdot \xi^{\beta}\, \ee^{- \mu\xi}
    \,\mathrm{d}\xi
  = {}- \int_0^{\infty}
    |u(\xi)|^2\cdot \xi^{\beta}\, (\ee^{- \mu\xi})_{\xi}
    \,\mathrm{d}\xi
\\
\label{e:int_parts}
& = {}
  - |u(\xi)|^2\cdot \xi^{\beta}\, \ee^{- \mu\xi}
    \Big\vert_{\xi = 0}^{\xi = \infty}
  + \int_0^{\infty}
    \left( |u(\xi)|^2\cdot \xi^{\beta} \right)_{\xi} \ee^{- \mu\xi}
    \,\mathrm{d}\xi
\\
\nonumber
& = \int_0^{\infty}
    ( |u(\xi)|^2 )_{\xi}\cdot \xi^{\beta}\, \ee^{- \mu\xi}
    \,\mathrm{d}\xi
  + \beta\int_0^{\infty}
    |u(\xi)|^2\cdot \xi^{\beta - 1}\, \ee^{- \mu\xi}
    \,\mathrm{d}\xi
\\
\nonumber
& = \int_0^{\infty}
  \left( u_{\xi}\, \bar{u} + u\, \bar{u}_{\xi} \right)
    \cdot \xi^{\beta}\, \ee^{- \mu\xi}
    \,\mathrm{d}\xi
  + \beta\int_0^{\infty}
    |u(\xi)|^2\cdot \xi^{\beta - 1}\, \ee^{- \mu\xi}
    \,\mathrm{d}\xi \,,
\end{align}
by the zero trace conditions \eqref{e:u(xi):trace:xi=0}.
We apply Cauchy's inequality,
\begin{equation*}
    u_{\xi}\, \bar{u} + u\, \bar{u}_{\xi}
  = 2\cdot \RE (u_{\xi}\, \bar{u})
  \leq 2\cdot | u_{\xi}\, \bar{u} |
  \leq \genfrac{}{}{}1{2}{\mu}\, |u_{\xi}|^2
     + \genfrac{}{}{}1{\mu}{2}\, |u|^2 \,,
\end{equation*}
to the integral
\begin{align*}
&   \int_0^{\infty}
  \left( u_{\xi}\, \bar{u} + u\, \bar{u}_{\xi} \right)
    \cdot \xi^{\beta}\, \ee^{- \mu\xi}
    \,\mathrm{d}\xi
\\
& \leq
    \frac{2}{\mu}\int_0^{\infty}
    |u_{\xi}(\xi)|^2\cdot \xi^{\beta}\, \ee^{- \mu\xi} \,\mathrm{d}\xi
  + \frac{\mu}{2}\int_0^{\infty}
    |u(\xi)|^2\cdot \xi^{\beta}\, \ee^{- \mu\xi} \,\mathrm{d}\xi \,.
\end{align*}
We estimate the last line in \eqref{e:int_parts} by this inequality,
thus arriving at
\begin{align*}
& \mu\int_0^{\infty}
    |u(\xi)|^2\cdot \xi^{\beta}\, \ee^{- \mu\xi} \,\mathrm{d}\xi
\\
&
\begin{aligned}
& \leq
    \frac{2}{\mu}\int_0^{\infty}
    |u_{\xi}(\xi)|^2\cdot \xi^{\beta}\, \ee^{- \mu\xi} \,\mathrm{d}\xi
  + \frac{\mu}{2}\int_0^{\infty}
    |u(\xi)|^2\cdot \xi^{\beta}\, \ee^{- \mu\xi} \,\mathrm{d}\xi \,,
\\
& {}
  + \beta\int_0^{\infty}
    |u(\xi)|^2\cdot \xi^{\beta - 1}\, \ee^{- \mu\xi} \,\mathrm{d}\xi \,.
\end{aligned}
\end{align*}
The desired inequality \eqref{ineq:Sobol:xi} follows.

Finally, we integrate ineq.~\eqref{ineq:Sobol:xi} with
$u$ replaced by
$\tilde{u}\equiv u(x,\,\cdot\,)\in W^{1,2}_{\mathrm{loc}}(0,\infty)$
(for almost every fixed $x\in \RR$)
with respect to the measure
$\ee^{- \gamma |x|} \,\mathrm{d}x$ over $\RR$ to obtain
ineq.~\eqref{ineq:Heston:Sobol}.
\qed
\par\vskip 10pt

Now we are ready to prove the following Hardy inequality.

\begin{lemma}\label{lem-Hardy}
{\rm (A Hardy\--type inequality.)}$\;$
Let\/
$\beta > 1$, $\gamma > 0$, and\/ $\mu > 0$.
Then the following {\em Hardy\--type inequality\/}
holds for every function $u\in V$,
\begin{align}
\label{e:Hardy-beta}
&
\begin{aligned}
  \int_{\HH} \genfrac{|}{|}{}0{u(x,\xi)}{\xi}^2
  \cdot \xi^{\beta}\, \ee^{ - \gamma |x| - \mu\xi }
    \,\mathrm{d}x \,\mathrm{d}\xi
& \leq \frac{8}{(\beta - 1)^2}
    \int_{\HH} |u_{\xi}(x,\xi)|^2
  \cdot \xi^{\beta}\, \ee^{ - \gamma |x| - \mu\xi }
    \,\mathrm{d}x \,\mathrm{d}\xi
\\
& {}
  + \frac{2\mu^2}{(\beta - 1)^2}
    \int_{\HH} |u(x,\xi)|^2
  \cdot \xi^{\beta}\, \ee^{ - \gamma |x| - \mu\xi }
    \,\mathrm{d}x \,\mathrm{d}\xi \,.
\end{aligned}
\end{align}
\end{lemma}

\proof
It suffices to verify the following inequality:
\begin{equation}
\label{e:Hardy-beta:xi}
\begin{aligned}
  \int_0^{\infty} \genfrac{|}{|}{}0{u(\xi)}{\xi}^2
  \cdot \xi^{\beta} \cdot \ee^{- \mu\xi}
  \,\mathrm{d}\xi
& \leq \frac{8}{(\beta - 1)^2}
    \int_0^{\infty} |u_{\xi}(\xi)|^2\cdot \xi^{\beta} \cdot \ee^{- \mu\xi}
    \,\mathrm{d}\xi
\\
& {}
  + \frac{2\mu^2}{(\beta - 1)^2}
    \int_0^{\infty} |u(\xi)|^2\cdot \xi^{\beta} \cdot \ee^{- \mu\xi}
    \,\mathrm{d}\xi
\end{aligned}
\end{equation}
holds for an arbitrary function
$u\in W^{1,2}_{\mathrm{loc}}(0,\infty)$ such that
\begin{equation}
\label{e:u,u'(xi):H^1}
  \int_0^{\infty} |u_{\xi}(\xi)|^2\cdot \xi^{\beta}\, \ee^{- \mu\xi}
    \,\mathrm{d}\xi < \infty \quad\mbox{ and }\quad
  \int_0^{\infty} |u(\xi)|^2\cdot \xi^{\beta}\, \ee^{- \mu\xi}
    \,\mathrm{d}\xi < \infty \,.
\end{equation}
The integrability hypotheses in \eqref{e:u,u'(xi):H^1}
are valid for $u$ replaced by the restricted function
$\tilde{u}\equiv u(x,\,\cdot\,)\in W^{1,2}_{\mathrm{loc}}(0,\infty)$
for a.e.\ fixed $x\in \RR$;
the first one by $u\in V$ and the second one by the previous lemma,
Lemma~\ref{lem-Sobolev}.

Inequality \eqref{e:Hardy-beta:xi} is obtained easily from
the standard weighted Hardy inequality
(\cite[Theorem 330, pp.\ 245--246]{Hardy-Polya}),
\begin{equation}
\label{eq:Hardy-beta}
  \int_0^{\infty} \genfrac{|}{|}{}0{f(\xi)}{\xi}^2 \cdot \xi^{\beta}
  \,\mathrm{d}\xi
  \leq \genfrac{(}{)}{}0{2}{\beta - 1}^2
  \int_0^{\infty} \genfrac{|}{|}{}0{\mathrm{d}f}{\mathrm{d}\xi}^2
  \cdot \xi^{\beta} \,\mathrm{d}\xi \,,
\end{equation}
where $\beta > 1$ and
$f\in W^{1,2}_{\mathrm{loc}}(0,\infty)$ satisfies
$\lim_{\xi\to \infty} f(\xi) = 0$, as follows:
We first replace the function $f$ by the product
$f(\xi) = u(x,\xi)\cdot \ee^{- \mu \xi/2}$,
then estimate the partial derivative
\begin{equation*}
\begin{aligned}
& f'(\xi)
  = \frac{\partial}{\partial\xi}
    \left( u(x,\xi)\cdot \ee^{- \mu \xi/2} \right)
  = u_{\xi}(x,\xi)\cdot \ee^{- \mu \xi/2}
  - \frac{\mu}{2}\, u(x,\xi)\cdot \ee^{- \mu \xi/2}
\\
& = \left( u_{\xi}(x,\xi) + \frac{\mu}{2}\, u(x,\xi) \right)
    \cdot \ee^{- \mu \xi/2}
\end{aligned}
\end{equation*}
by
\begin{equation*}
  |f'(\xi)|^2
  = \left|
    \frac{\partial}{\partial\xi}
    \left( u(x,\xi\cdot \ee^{- \mu \xi/2} \right)
    \right|^2
  \leq
    2\left[ |u_{\xi}(x,\xi)|^2 + \genfrac{(}{)}{}0{\mu}{2}^2\, |u(x,\xi)|^2
     \right] \cdot \ee^{- \mu\xi}
\end{equation*}
and insert it into ineq.~\eqref{eq:Hardy-beta},
thus arriving at ineq.~\eqref{e:Hardy-beta:xi}.
Here, the hypothesis
$f\in W^{1,2}_{\mathrm{loc}}(0,\infty)$ is satisfied, thanks to
$u\in V$, whence even
\begin{math}
  \int_0^{\infty} |f'(\xi)|^2\cdot \xi^{\beta} \,\mathrm{d}\xi < \infty \,,
\end{math}
with a help from \eqref{e:u,u'(xi):H^1}.
Hypothesis $\lim_{\xi\to \infty} f(\xi) = 0$
follows from the trace result
\eqref{lim:trace_pt:v=infty} in Lemma~\ref{lem-Trace_pt}.

The proof is completed by integrating
ineq.~\eqref{e:Hardy-beta:xi} with
$u$ replaced by
$\tilde{u}\equiv u(x,\,\cdot\,)\in W^{1,2}_{\mathrm{loc}}(0,\infty)$
(for a.e.\ $x\in \RR$)
with respect to the measure
$\ee^{- \gamma |x|} \,\mathrm{d}x$ over $\RR$ to obtain
ineq.~\eqref{e:Hardy-beta}.
\qed
\par\vskip 10pt

Recall that any function
$u\in V = H^1(\mathbb{H};\mathfrak{w})$
satisfies the hypotheses of
Lemmas \ref{lem-Sobolev} and~\ref{lem-Hardy} above.

\begin{remark}\label{rem-equiv_norm}\nopagebreak
\begingroup\rm
Owing to the Sobolev- and Hardy\--type inequalities
\eqref{ineq:Heston:Sobol} and \eqref{e:Hardy-beta}
proved in Lemmas \ref{lem-Sobolev} and~\ref{lem-Hardy},
with $1 < \beta < \infty$,
the following inner product defines an {\em equivalent norm\/}
on the Hilbert space~$V$:
\begin{equation}
\label{def_eq:w_prod-H^1}
  (u,w)_{V}^{\sharp} \eqdef (u,w)_V + (u,w)_{V}^{\flat}
    \quad\mbox{ for }\, u,w\in V \,,
\end{equation}
where
\begin{align}
\label{def_HS:w_prod-H^1}
&
\begin{aligned}
  (u,w)_V^{\flat}
& \eqdef
    \int_{\mathbb{H}} \frac{u(x,\xi)}{\xi}\cdot
                      \frac{\bar{w}(x,\xi)}{\xi}\cdot
                      \xi\cdot \mathfrak{w}(x,\xi)
    \,\mathrm{d}x \,\mathrm{d}\xi
\\
& {}
  + \int_{\mathbb{H}} u\, \bar{w}\cdot \xi\cdot \mathfrak{w}(x,\xi)
    \,\mathrm{d}x \,\mathrm{d}\xi
\end{aligned}
\\
\nonumber
& {}
  = \int_{\mathbb{H}} u\, \bar{w}\left( \xi + \frac{1}{\xi}\right)
    \mathfrak{w}(x,\xi) \,\mathrm{d}x \,\mathrm{d}\xi
    \quad\mbox{ for }\, u,w\in V \,.
\end{align}
This fact is used in paragraphs
{\S}\ref{ss:bound-R} and {\S}\ref{ss:coerce-R}.
\hfill\Square
\endgroup
\end{remark}
\par\vskip 10pt


\section{Appendix: Density of entire functions\\ in
         $H = L^2(\mathbb{H};\mathfrak{w})$}
\label{s:density}

As we have already suggested in paragraph {\S}\ref{ss:Cauchy-compl},
we wish to approximate an arbitrary initial condition
$u_0\in H = L^2(\mathbb{H};\mathfrak{w})$
by a sequence of entire functions,
$u_{0,n}\colon \CC^2\to \CC$; $n=1,2,3,\dots$,
such that their restrictions $u_{0,n}\vert_{\HH}$ to
$\HH = \RR\times (0,\infty)$ satisfy
\begin{equation*}
  \| u_{0,n}\vert_{\HH} - u_0\|_H \;\longrightarrow\; 0
    \quad\mbox{ as }\, n\to \infty \,.
\end{equation*}
Below, we construct rather simple
{\em entire\/} (holomorphic) functions
\begin{math}
  u_{0,n}\colon \CC^2\to \CC ;
\end{math}
$n=1,2,3,\dots$, with this property, by using standard results about
{\em Hermite\/} and {\em Laguerre\/} functions.
The reader is referred to the monographs by
{\sc A.~N.\ Kolmogorov\/} and {\sc S.~V.\ Fomin\/}
\cite[Chapt.~VII, {\S}3.7, pp.\ 395--396]{Kolm-Fomin}
and
{\sc N.~N.\ Lebedev} \cite[Chapt.~4]{Lebedev},
{\S}4.9, pp.\ 60--61 and {\S}4.17, pp.\ 76--78,
for details and proofs.


\subsection{Hermite and Laguerre functions in the complex domain}
\label{ss:Hermite-Laguerre}

In our approximation procedure below, we first take advantage of
the (complex) Hilbert space
$H = L^2(\mathbb{H};\mathfrak{w})$
being the {\it tensor product\/} of the Hilbert spaces
$\mathfrak{H}_1 = L^2(\RR;\mathfrak{w}_1)$ and
$\mathfrak{H}_2 = L^2(\RR_+;\mathfrak{w}_2)$, with the weights
\begin{equation}
\label{def:w_1,2}
  \mathfrak{w}_1(x)\eqdef \ee^{- \gamma |x|} \quad\mbox{ and }\quad
  \mathfrak{w}_2(\xi)\eqdef \xi^{\beta - 1}\, \ee^{- \mu\xi}
  \quad\mbox{ for $(x,\xi)\in \HH$, }
\end{equation}
i.e.,
$H = \mathfrak{H}_1\otimes \mathfrak{H}_2$, as defined in
{\sc M.\ Reed} and {B.\ Simon} \cite[Chapt.~II, {\S}4]{RSimon-I},
pp.\ 49--54.
All general properties of a tensor product of two Hilbert spaces
that we use below can be found there.
Thus, both $\mathfrak{H}_1$ and $\mathfrak{H}_2$
are weighted Lebesgue $L^2$-spaces with the weighted Lebesgue measures
$\mathfrak{w}_1(x) \,\mathrm{d}x$ and
$\mathfrak{w}_2(x) \,\mathrm{d}\xi$, respectively.

In order to keep our approximation procedure simple,
we take advantage of the density of the weighted Lebesgue $L^2$-spaces
as follows:
$L^2(\HH)$ is densely and continuously imbedded into $H$,
$L^2(\RR)$ into $\mathfrak{H}_1$, and $L^2(\RR_+)$ into $\mathfrak{H}_2$.
This claim is an easy consequence of the fact that all weights,
$\mathfrak{w}(x,\xi) = \mathfrak{w}_1(x)\cdot \mathfrak{w}_2(\xi)$,
$\mathfrak{w}_1(x)$, and $\mathfrak{w}_2(\xi)$ are bounded.

We use a standard approximation method in $\mathfrak{H}_1$ by
{\it Hermite functions\/},
\begin{math}
  h(x) =\hfil\break p(x)\, \exp\left( {}- \frac{1}{2} x^2 \right) ,
\end{math}
where $p(x)$ is a polynomial obtained by a linear combination of
{\it Hermite polynomials\/} $H_n(x)$; $n=0,1,2,\dots$.
We refer to
{\sc N.~N.\ Lebedev} \cite[{\S}4.9, pp.\ 60--61]{Lebedev}
for a common definition of Hermite polynomials
and their basic properties.
In particular, $H_n(x)$ is a polynomial of degree $n\geq 0$
and the Hermite functions
\begin{equation*}
  h_n(x) = H_n(x)\,
    \exp\left( {}- \genfrac{}{}{}1{1}{2} x^2 \right)
  \quad\mbox{ of }\, x\in \RR \,;\quad n=0,1,2,\dots \,,
\end{equation*}
form an orthonormal basis in $L^2(\RR)$, by
{\sc N.~N.\ Lebedev} \cite[{\S}4.13, pp.\ 65--66]{Lebedev}.
Furthermore, an arbitrary linear combination of these functions,
\begin{math}
  h(x) = p(x)\, \exp\left( {}- \frac{1}{2} x^2 \right) ,
\end{math}
where $p(x)$ is a polynomial, can be extended uniquely to
an entire function
\begin{math}
  \tilde{h}(z) = p(z)\, \exp\left( {}- \frac{1}{2} z^2 \right)
\end{math}
of the complex variable $z = x + \ii y\in \CC$.
Finally, given any $r > 0$ and $\delta > 0$, there is a constant
$C_{r,\delta,p}\in (0,\infty)$,
depending only on $r$, $\delta$, and the polynomial $p$, such that
the following inequalities hold for all
$z = x + \ii y, z^{*}\in \CC$ with $|y|\leq r$ and $|z^{*}|\leq \delta$:
\begin{align}
\nonumber
& |\tilde{h}(x + \ii y + z^{*})| = |p(x + \ii y + z^{*})|\cdot
    \exp\left( {}- \genfrac{}{}{}1{1}{2}\cdot \RE [ (x + \ii y + z^{*})^2 ]
        \right)
\\
\label{est:x}
& = |p(x + \ii y + z^{*})|\cdot
    \exp\left( {}- \genfrac{}{}{}1{1}{2}\cdot \RE
    \left[ (x + \ii y)^2 + 2\, (x + \ii y) z^{*} + (z^{*})^2 \right]
        \right)
\\
\nonumber
& \leq |p(x + \ii y + z^{*})|\cdot
    \exp\left( {}- \genfrac{}{}{}1{1}{2}\cdot
    \left[ x^2 - y^2 - 2\, (|x| + |y|)\, |z^{*}| - |z^{*}|^2 \right]
        \right)
\\
\nonumber
& \leq C_{r,\delta,p}\cdot
    \exp\left( {}- \genfrac{}{}{}1{1}{2}\, x^2 + 2\delta\, |x| \right) \,.
\end{align}
Consequently, the square of the $L^2(\RR)$-norm of the function
$x\mapsto \tilde{h}(x + \ii y + z^{*})\colon \RR\to \CC$
is uniformly bounded, provided
$|y|\leq r$ and $|z^{*}|\leq \delta$ are satisfied:
\begin{align*}
  \int_{-\infty}^{+\infty} |\tilde{h}(x + \ii y + z^{*})|^2 \,\mathrm{d}x
  \leq C_{r,\delta,p}^2\cdot
  \int_{-\infty}^{+\infty}
    \exp\left( {}- x^2 + 4\delta\, |x| \right) \,\mathrm{d}x
  \equiv \mathrm{const}_{r,\delta,p}^2 < \infty \,.
\end{align*}
A Hermite polynomial based expansion has already been applied to
Black\--Scholes and Merton type models for European option prices,
e.g., in the recent work by
{\sc D.\ Xiu\/} \cite{Xiu}.

Analogously, in $\mathfrak{H}_2$ we use
{\it Laguerre functions\/},
\begin{math}
  \ell(\xi) = q(\xi)\, \exp\left( {}- \frac{1}{2}\xi \right) ,
\end{math}
where $q(\xi)$ is a polynomial obtained by a linear combination of
{\it Laguerre polynomials\/} $L_n(\xi)$; $n=0,1,2,\dots$.
We refer to
{\sc N.~N.\ Lebedev} \cite[{\S}4.17, pp.\ 76--78]{Lebedev}
for a common definition of Laguerre polynomials
and their basic properties.
In particular, $L_n(\xi)$ is a polynomial of degree $n\geq 0$
and the Laguerre functions
\begin{equation*}
  \ell_n(\xi) = L_n(\xi)\,
    \exp\left( {}- \genfrac{}{}{}1{1}{2}\xi \right)
  \quad\mbox{ of }\, \xi\in \RR_+ \,;\quad n=0,1,2,\dots \,,
\end{equation*}
form an orthonormal basis in $L^2(\RR_+)$, by
{\sc N.~N.\ Lebedev} \cite[{\S}4.21, pp.\ 83--84]{Lebedev}.
Furthermore, an arbitrary linear combination of these functions,
\begin{math}
  \ell(\xi) = q(\xi)\, \exp\left( {}- \frac{1}{2}\xi \right) ,
\end{math}
where $q(\xi)$ is a polynomial, can be extended uniquely to
an entire function
\begin{math}
  \tilde{\ell}(\zeta) = q(\zeta)\, \exp\left( {}- \frac{1}{2} \zeta \right)
\end{math}
of the complex variable
$\zeta = \xi (1 + \ii\omega) \in \CC$.
Finally, given any $\vartheta_v > 0$ and $\delta > 0$, there is a constant
$C_{{\vartheta_v},\delta,q}\in (0,\infty)$,
depending only on $\vartheta_v$, $\delta$, and the polynomial $q$,
such that the following inequalities hold for all
$\zeta = \xi (1 + \ii\omega), \zeta^{*}\in \CC$ with $\xi\in \RR_+$,
$|\arctan\omega|\leq \vartheta_v$, and $|\zeta^{*}|\leq \delta$:
\begin{align}
\label{est:xi}
&   |\tilde{\ell}( \xi (1 + \ii\omega) + \zeta^{*} )|
  = |q( \xi (1 + \ii\omega) + \zeta^{*} )|\cdot
    \exp\left( {}- \genfrac{}{}{}1{1}{2}\cdot
    \RE [ \xi (1 + \ii\omega) + \zeta^{*} ]
        \right)
\\
\nonumber
& \leq
    |q( \xi (1 + \ii\omega) + \zeta^{*} )|\cdot
    \exp\left( {}- \genfrac{}{}{}1{1}{2}\cdot ( \xi - |\zeta^{*}| )
        \right)
  \leq C_{{\vartheta_v},\delta,q}\cdot
    \exp\left( {}- \genfrac{}{}{}1{1}{4}\, \xi \right) \,.
\end{align}
Consequently, the square of the $L^2(\RR_+)$-norm of the function
\begin{math}
  \xi\mapsto \tilde{\ell}( \xi (1 + \ii\omega) + \zeta^{*} )
  \colon \RR_+\to \CC
\end{math}
is uniformly bounded, provided
$|\arctan\omega|\leq \vartheta_v$ and $|\zeta^{*}|\leq \delta$
are satisfied:
\begin{align*}
  \int_0^{+\infty} |\tilde{\ell}( \xi (1 + \ii\omega) + \zeta^{*} )|^2
    \,\mathrm{d}\xi
  \leq C_{{\vartheta_v},\delta,q}^2\cdot
  \int_0^{+\infty}
    \exp\left( {}- \genfrac{}{}{}1{1}{2}\, \xi \right) \,\mathrm{d}\xi
  = 2\, C_{{\vartheta_v},\delta,q}^2 < \infty \,.
\end{align*}

Summarizing the properties of the Hermite and Laguerre functions,
we observe that the product functions
\begin{equation*}
  e_{mn}(x,\xi)\eqdef h_m(x)\, \ell_n(\xi)
  \quad\mbox{ of }\, (x,\xi)\in \HH \,;\quad m,n=0,1,2,\dots \,,
\end{equation*}
form an orthonormal basis in $L^2(\HH)$
(\cite[Chapt.~II, {\S}4]{RSimon-I}).


\subsection{Approximation of the initial conditions\\
           (Gal\"erkin's method)}
\label{ss:Galerkin}

We have just shown that, given any initial condition
$u_0\in H = L^2(\mathbb{H};\mathfrak{w})$,
there is a sequence of {\em entire\/} (holomorphic) functions
\begin{equation*}
    u_{0,n}(z,\zeta)
  = P_n(z,\zeta)\, \exp\left( {}- \frac{1}{2}\, (z^2 + \zeta) \right) \,,
    \quad (z,\zeta)\in \CC^2 \,;\quad n=1,2,3,\dots \,,
\end{equation*}
with the restrictions $u_{0,n}\vert_{\HH}$
in the tensor product
$L^2(\HH) = L^2(\RR)\otimes L^2(\RR_+)$ $\hookrightarrow$
$H = \mathfrak{H}_1\otimes \mathfrak{H}_2$, such that:
\begin{itemize}
\item[{\rm (i)}]
$P_n\colon \CC^2\to \CC$ is a polynomial with complex coefficients.
\item[{\rm (ii)}]
The restrictions $u_{0,n}\vert_{\HH}$ of $u_{0,n}$ to
$\HH = \RR\times (0,\infty)$ satisfy
\begin{equation*}
  \| u_{0,n}\vert_{\HH} - u_0\|_H \;\longrightarrow\; 0
    \quad\mbox{ as }\, n\to \infty \,.
\end{equation*}
\item[{\rm (iii)}]
There is a constant
$K_n\equiv K_{P_n}\in (0,\infty)$,
depending on $P_n$, $r$, and $\vartheta_v$,
$0 < r < \infty$ and $0 < \vartheta_v < \pi / 2$,
but independent from $y, \omega\in \RR$ in
$z = x + \ii y ,\, \zeta = \xi (1 + \ii\omega) \in \CC$ and
$z^{*}, \zeta^{*}\in \CC$
with
$|y| < r$, $|\arctan\omega| < \vartheta_v$, and
$\max\{ |z^{*}| ,\, |\zeta^{*}| \} < \delta$, such that
\begin{align*}
& \int_{\HH}
    \left| u_{0,n}
    \left( x + \ii y + z^{*} ,\, \xi (1 + \ii\omega) + \zeta^{*} \right)
    \right|^2 \,\mathrm{d}x \,\mathrm{d}\xi
  \leq K_n\equiv \mathrm{const} < \infty
\\
&   \quad\mbox{ whenever }\,
    |y| < r ,\quad |\arctan\omega| < \vartheta_v ,\ \mbox{ and }\,
    \max\{ |z^{*}| ,\, |\zeta^{*}| \} < \delta \,.
\end{align*}
An analogous estimate remains valid in the weighted Lebesgue space $H$
if the standard Lebesgue measure $\,\mathrm{d}x \,\mathrm{d}v$
is replaced by the weighted Lebesgue measure
$\mathfrak{w}(x,v) \,\mathrm{d}x \,\mathrm{d}v$, thanks to
$0 < \mathfrak{w}(x,v)\leq \mathrm{const} < \infty$.
\end{itemize}

Notice that the estimate in {\rm (iii)} above follows from
\begin{align}
\nonumber
& \int_{\HH}
    \left| u_{0,n}
    \left( x + \ii y + z^{*} ,\, \xi (1 + \ii\omega) + \zeta^{*} \right)
    \right|^2 \,\mathrm{d}x \,\mathrm{d}\xi
\\
\label{e:y,eta_0-fixed}
&
\begin{aligned}
  = \int_0^{\infty} \int_{-\infty}^{\infty}
&   \left| P_n
    \left( x + \ii y + z^{*} ,\, \xi (1 + \ii\omega) + \zeta^{*} \right)
    \right|^2
\\
& {}\times \exp
      \left( {}- \RE
           [ (x + \ii y + z^{*})^2 + \xi (1 + \ii\omega) + \zeta^{*} ]
      \right) \,\mathrm{d}x \,\mathrm{d}\xi
\end{aligned}
\\
\nonumber
&
\begin{aligned}
  \leq \int_0^{\infty} \int_{-\infty}^{\infty}
&   \left| P_n
    \left( x + \ii y + z^{*} ,\, \xi (1 + \ii\omega) + \zeta^{*} \right)
    \right|^2\cdot
      \exp\left( {}- (x^2 - y^2)  - \xi \right)
\\
& {}\times \exp
    \left( 2 |x + \ii y|\cdot |z^{*}| + |z^{*}|^2 + |\zeta^{*}| \right)
    \,\mathrm{d}x \,\mathrm{d}\xi
\end{aligned}
\\
\nonumber
&
\begin{aligned}
  \leq \int_0^{\infty} \int_{-\infty}^{\infty}
&   \left| P_n
    \left( x + \ii y + z^{*} ,\, \xi (1 + \ii\omega) + \zeta^{*} \right)
    \right|^2\cdot
      \exp\left( {}- x^2 - \xi \right)
\\
& {}\times \exp
    \left( r^2 + 2 (|x| + r)\delta + \delta^2 + \delta \right)
    \,\mathrm{d}x \,\mathrm{d}\xi
\end{aligned}
\\
\nonumber
& \leq K_n\equiv \mathrm{const} < \infty
\\
\nonumber
&   \quad\mbox{ whenever }\,
    |y| < r ,\quad |\arctan\omega| < \vartheta_v ,\ \mbox{ and }\,
    \max\{ |z^{*}| ,\, |\zeta^{*}| \} < \delta \,.
\end{align}

As an obvious consequence of properties
{\rm (i)}, {\rm (ii)}, and {\rm (iii)} we obtain that
\begin{math}
  u_{0,n}\colon \mathfrak{X}^{(r)} \times \Delta_{\vartheta_v} \to \CC
\end{math}
is a holomorphic function in both its variables $(z,\zeta)$
and belongs to the Hardy space
$H^2( \mathfrak{X}^{(r)} \times \Delta_{\vartheta_v} )$.


\section*{Acknowledgment.}
This work was supported in part by
le Minist\`ere des Affaires \'Etrang\`eres (France)
and
the German Academic Exchange Service (DAAD, Germany)
within the exchange program ``PROCOPE'' between France and Germany.
It was performed while the second author (P.T.) was
a Visiting Professor (``professeur invit\'e'')
at Toulouse School of Economics, I.M.T.,
Universit\'e de Toulouse -- Capitole, Toulouse, France.
The authors would like to thank Professor Paul M.~N.\ Feehan
(Rutgers University, N.J., U.S.A.)
for his kind advice concerning properties and use of
weighted Lebesgue and Sobolev spaces introduced in
\cite{Daska-Feehan-14, Daska-Feehan-16, Feehan-Pop-14}.


%
%
\makeatletter \renewcommand{\@biblabel}[1]{\hfill#1.} \makeatother
\ifx\undefined\bysame
\newcommand{\bysame}{\leavevmode\hbox to3em{\hrulefill}\,}
\fi
%

%
%
%
\end{document}